%% file: teseleo.tex
\documentclass[12pt]{report}
\usepackage{amstext}

\usepackage[brazil]{babel}
\usepackage[latin1]{inputenc}
\usepackage{amsfonts}
\usepackage{amssymb}
\usepackage{color}
\usepackage{graphicx}

\include{macro}

\begin{document}

\include{capa}

\include{dedic}
 \setlength{\parskip}{0.4cm} \baselineskip 22pt
\include{agrad}

 \setlength{\parskip}{0.4cm} \baselineskip 22pt
\include{abstract}

\include{resumo}
 \setlength{\parskip}{0.0cm}
 \baselineskip 16pt
\tableofcontents
 \setlength{\parskip}{0.4cm}
 \baselineskip 22pt

\include{intro}

\include{rgana}

\include{rgnum}

\include{resul}

\include{homog}

\include{apend}

\nocite{braga03b} \nocite{braga04}
\addcontentsline{toc}{chapter}{Referências Bibliográficas}
 \bibliographystyle{\bibdir mylib}
\bibliography{\bibdir bibliography}

\include{errata}

\end{document}

%% file: macro.tex
\newcommand{\bibdir}{}

\newcommand{\prova}{\noindent \emph{Prova:} }
\newcommand{\qed}{\nopagebreak\hfill\fbox{ }}
\newcommand{\eop}{\qed}

\newcommand{\R}{\mathbb{R}}

\newcommand{\N}{\mathbb{N}}

\newcommand{\B}{{\cal B}}
\newcommand{\F}{{\cal F}}

\newcommand{\oo}{\infty}
 \newcommand{\tot} { \displaystyle \mathop{\longrightarrow}_{ {t\to\oo}} }
 \newcommand{\ton} { \displaystyle \mathop{\longrightarrow}_{ {n\to\oo}} }
 \newcommand{\toe} { \displaystyle \mathop{\longrightarrow}_{{\epsilon\to0}} }

\newcommand{\quando}{,\quad \mbox{quando } t \gg 1}

\newcommand{\A}{{\cal A}}
\newcommand{\CC}{{\cal C}}
\newcommand{\sgn}{{\rm sgn}}

\newcommand{\pot}{{\rm pot}}

\renewcommand{\div}{{\rm div}}
\newcommand{\cont}{\subset}
\newcommand{\D}{\mathrm{d}}
\newcommand{\dx} {{\partial \over \partial x  }}
\newcommand{\dxi}{{\partial \over \partial x_i}}

\newcommand{\LL}{\mathbf{L}}
\newcommand{\grad}{\nabla}

\newcommand{\dlim}{\displaystyle \lim}

\newcommand{\dsup}{\displaystyle \sup}

\newcommand{\av}[1]{\left\langle#1\right\rangle}
\newcommand{\der}[2]{{\partial#1 \over \partial #2}}
\newcommand{\Der}[2]{{\D#1 \over \D #2}}

\definecolor{red}{cmyk}         {0,1,1,0.1}
\definecolor{blue}{cmyk}        {1,1,0,0.2}
\definecolor{pink}{cmyk}        {0.,1,0.,0.2}
\definecolor{lightpink}{cmyk}   {.0,1,.5,.0}
\definecolor{darkblue}{cmyk}    {1,1,0,.6}
\definecolor{black}{cmyk}       {0,0,0,1}
\definecolor{gray}{cmyk}        {0,0,0,.5}
\definecolor{yellow}{cmyk}{0,0,1,0.55}

\newcommand{\textgray}[1]{}

\newcommand{\obs}{\noindent {\bf Observa\c c\~ao: }}

\newcommand{\ba}{\begin{array}}
\newcommand{\ea}{\end{array}}
\newcommand{\be}{\begin{equation}}
\newcommand{\ee}{\end{equation}}
\newcommand{\ben}{\begin{enumerate}}
\newcommand{\een}{\end{enumerate}}

\newtheorem{teo}{Teorema}[chapter]
\newtheorem{conj}[teo]{Conjectura}
\newtheorem{defi}[teo]{Defini\c c\~ao}
\newtheorem{prop}[teo]{Proposi\c c\~ao}
\newtheorem{lema}[teo]{Lema}
\newtheorem{cor}[teo]{Corol\'ario}
\newtheorem{prob}{Problema}

\newcommand{\insertfigure}[3]
{\begin{figure}[htb]
  \begin{center}
    \includegraphics[height=8cm]{#2}
  \end{center}
  \caption{#3}
  \label{#1}
\end{figure}}

\newcommand{\inserttable}[4]
{\begin{table}[htb]
  \begin{center} \small
    \begin{tabular}{#2} #4 \end{tabular}
  \end{center}
  \caption{#3}
  \label{#1}
\end{table}}

%% file: capa.tex
\thispagestyle{empty}
\begin{center}
{\Large UNIVERSIDADE FEDERAL DE MINAS GERAIS}\\{INSTITUTO DE
CIÊNCIAS EXATAS}\\
{\sc Departamento de Matemática}
\end{center}

\hspace*{-2.0cm}

\vspace{3cm}
\begin{center}
{\bf Dissertação de Mestrado}
\end{center}
\vspace{2cm}

\baselineskip 2.3em

\begin{center}
{\LARGE
\textbf{ANÁLISE ASSINTÓTICA DE SOLUÇÕES
DE EQUAÇÕES DIFUSIVAS NÃO-LINEARES VIA
MÉTODOS DE ESCALAS MÚLTIPLAS}
}
\end{center}

\baselineskip 1.3em

\begin{center}

\normalsize \vspace{10mm}

\vspace{15mm}

{\large {\bf Leonardo Trivellato Rolla}}

\vspace{10mm}
\end{center}

\begin{center}

{\bf Orientador: {\large Gastão de Almeida Braga}}
\end{center}
\vspace{10mm}

\begin{center} \small

\normalsize

\vfill \vspace{8mm} \vfill
\textsc{20 de Setembro de 2004}
\end{center}


%% file: dedic.tex
\ \vfill
\begin{flushright}
\textit{Dedico este trabalho
\\ \ 
\\ \ 
\\
à minha mãe e aos meus irmãos,
\\ \
\\ \
\\
meus maiores motivos de alegria.
\\ \
\\ \
}
\end{flushright}

%% file: agrad.tex
\noindent {\Huge \bf Agradecimentos}

\vspace{10mm}

Agradeço ao CNPq, pela bolsa de estudos.

As pessoas, são várias às quais devo agradecimentos. Peço desculpas a
quem tenha sido omitido nesta página, por algum descuido meu. Agradeço:

a todos os meus verdadeiros amigos, pela cumplicidade incondicional e
pela companhia agradável, mesmo que às vezes aparentemente separados pelo
tempo, pela rotina ou pela distância;

aos demais alunos da pós-graduação, pelo coleguismo e pela prazerosa
convivência no nosso querido ``beco''. Em especial, aos colegas
Alexandre, Ivana, Leandro e Rodrigo;

aos camaradas Bruno, Dalgiza, Eliana, Érico, José Raimundo, Leonardo e
Marquinhos, e a todos aqueles que lutam pela construção de uma sociedade
mais justa, pelo companheirismo e por me ajudarem a dar sentido prático à
minha existência enquanto ser político e social;

aos Profs.~Antônio Plascak, Carlos Isnard, Gastão Braga, Mário Jorge
Caneiro, Márcio Soares, Suzana Fornari, Sylvie Kamphorst, Vladas
Sidoravicius e Wagner Meira, por me mostrarem que não existe dicotomia
entre pesquisar com excelência e lecionar com clareza, dedicação e
entusiasmo;

à Profa.~Sylvie Kamphorst, por todo o apoio, preocupa\c c\~ao e aten\c
c\~ao que sempre teve comigo, pelos incentivos e pelos v\'arios conselhos
dados durante nossas longas conversas, que jamais serão esquecidas;

à Jussara Moreira, pela colaboração nos trabalhos publicados e pela ajuda
nos estudos analíticos;

ao Prof.~Paulo Cupertino, pelos valiosos comentários durante os
seminários de homogeneização;

ao Prof.~Vincenzo Isaia, pelas várias discussões, em especial quanto à
simulação do problema periódico;

ao Prof.~Paulo César Carrião, por ter aceitado o convite para participar
da banca e pelo cuidado matemático minuncioso na leitura e crítica deste
texto;

ao Prof.~Luis Felipe Pereira, pelos ensinamentos durante minha estada no
IPRJ, pela singular hospitalidade, por ter aceitado participar da banca
e, principalmente, pelas sugestões e comentários que certamente
contribuiram com a qualidade deste trabalho;

ao Prof.~Frederico Furtado, pelos momentos de trabalho conjunto, com a
seriedade e ao mesmo tempo o bom homor que lhe são peculiares. Pelas
sugestões, críticas, comentários, incentivos e ``palpites'' dados nos
últimos três anos, que sem dúvida alguma foram de enorme valia não só
para esta dissertação como para a minha formação individual;

ao Prof.~Gastão Braga, por ter transmitido, não apenas conhecimento, mas
principalmente sabedoria, com sua paciência incomensurável e, por isso
mesmo, com máxima eficácia. Pela dediação, proteção, compreensão e
altruísmo absolutamente atípicos em uma relação oriendador-aluno. Enfim,
agradeço por ter sido, nos últimos quatro anos, o melhor orientador que
eu poderia esperar.

%% file: abstract.tex
\noindent {\Huge \bf Abstract}

\vspace{1cm}

In the present work we shall describe and apply the techniques of the Renormalization Group - based in data rescaling and operator renormalizing - and of Homogenization - that substitutes, in a  certain limit, a periodically inhomogeneous medium by a homogeneous one - for the study of asymptotic behavior of solutions of PDE's.

One class of problems studied in this dissertation consists of nonlinear, diffusive differential equations with periodic coefficients, that supposedly model the diffusion in a heterogeneous medium. Our study of these problems involves arguments originated from the two techniques described above. This makes, in certain sense, a connection between both of them.

Another class of problems studied here consists of nonlinear, diffusive equations with time dependent coefficients, originated from montecarlo simulation of the asymptotic behavior of the averaged solutions of stochastic differential equations modeling the phenomenon of two-phase flow in porous media.

Using the Renormalization Group we shall classify the problems in accordance with their qualitative behavior. This classification will be verified numerically. We shall also study the curve of phase transition that separates asymptotic regimes predominantly linear and nonlinear. The main results presented in this work are given by a numerical version of the Renormalization Group, which will be described in detail in the body of this dissertation.

%% file: resumo.tex
\noindent {\Huge \bf Resumo}

\vspace{1cm}

Neste trabalho vamos descrever e aplicar as técnicas do Grupo de
Renormalização - baseada em mudança de escalas e renormalização de
operadores - e de Homogeneização -  que substitui, num
certo limite, um meio periodicamente não-homogêneo por um outro homogêneo
 - para o estudo do comportamento assintótico
de soluções de EDP's.

Uma classe de problemas estudados nesta dissertação consiste em equações
diferenciais difusivas não-lineares com coeficiente de difusão periódico,
que supostamente modelam a difusão em um meio heterogêneo. O nosso estudo
desses problemas envolve argumentos originados das duas técnicas descritas
acima, fornecendo, em certo sentido, uma conexão entre ambas.

Outra classe de problemas estudados nesta dissertação consiste em equações
difusivas não-lineares com coeficiente de difusão dependente do tempo, que
têm sua motivação física dada pelo comportamento assintótico médio, com
respeito a simulações estocásticas, do fenômeno de escoamento bifásico em
meios porosos.

Utilizando o grupo de renormalização vamos classificar os problemas de
acordo com seu comportamento qualitativo. Esse comportamento qualitativo
será verificado numericamente. Vamos também estudar a curva de transição
de fase entre os regimes assintóticos predominantemente linear e o
não-linear.  Os principais resultados apresentados neste trabalho são
obtidos via uma versão numérica do grupo de renormalização, que será
descrita em detalhes no corpo desta dissertação.

%% file: intro.tex
\chapter{Introdução}
\label{chap:intro}

\section{Auto-similaridade e universalidade}

A evolução no tempo de sistemas físicos fora do estado de
equilíbrio é muitas vezes bem descrita por soluções de
equações diferenciais parciais.
Freqüentemente é observado que, sob circunstâncias apropriadas
(mas razoavelmente gerais), a evolução temporal de tais sistemas
torna-se
assintoticamente auto-similar\footnote
{O significado de ``auto-similar'' está relacionado com a
invariância com respeito à ação de um certo grupo de simetria,
como será explicado logo adiante.}.
No contexto de EDP's, tal
afirmação traduz-se no fato de a equação em questão ter soluções
que se comportam como
\be
  \label{eq:noprefss}
  u(x,t) \approx \frac{1}{t^{\alpha}}\ \phi\left( \frac{x}{t^{\beta}}\right)
  \,\,\, \mbox{quando}\,\,\, t\rightarrow \infty.
\ee

Ocorre muitas vezes que tal comportamento auto-similar
possui um forte viés de \emph{universalidade}: os expoentes
$\alpha$ e $\beta$ e a maioria dos aspectos da função
$\phi(\cdot)$ independem do problema de valor inicial
associado, isto é, independem dos dados iniciais (dentro
de uma certa classe), dos parâmetros da equação ou até mesmo da
sua forma.

A análise de comportamentos assintóticos auto-similares
envolve um problema em escalas múl\-ti\-plas. A palavra
\emph{escala} se refere à escolha de uma unidade de medida, tanto
para a variável espacial $x$ quanto para a variável temporal
$t$. Estas escalas são interdependentes e em geral elas surgem de
uma determinada simetria da EDP em questão -- ou parte dela.
Essa simetria consiste no fato de que existem
números reais $\alpha$ e $\beta$ e uma função $\phi(\cdot)$ tais
que a equação
admita uma solução da forma
\be
  \label{eq:sssol}
  u(x,t)= \frac{1}{t^{\alpha}}\ \phi\left(\frac{x}{t^{\beta}}\right).
\ee
Essas soluções são ditas \emph{auto-similares} ou
\emph{invariantes por mudança de escala}.
Uma função homogênea de grau $-\alpha$ é o exemplo mais
simples de uma função auto-similar. De fato, considere a mudança de
escala $t\to L\, t$, onde $L$ é um fator positivo, e considere
o conjunto ${\cal C}$ das funções $f(\cdot):\R_+^*\to\R$. Dado
$\alpha\in\R$, defina, sobre ${\cal C}$, a operação de mudança de escalas
$\tilde f(t) \equiv L^\alpha f(Lt)$. Se
$y(t)\equiv At^{-\alpha}$, $t>0$, então $y(t)\in {\cal C}$ e
$\tilde y(t) = y(t)$, isto é, $y(t)$ é invariante pela mudança de
escala em  ${\cal C}$.

De forma análoga, assumindo a existência dos expoentes
$\alpha$ e $\beta$ em~(\ref{eq:sssol}) e dados $L>0$  e
$u:\R\times\R_+^*\to\R$, definimos a mudança de escalas
$u\mapsto v$ por
\be
  \label{eq:rescale}
  v(x,t) = L^\alpha u(L^\beta x,Lt).
\ee
Claramente, soluções da forma~(\ref{eq:sssol}) são invariantes
por esta mudança. Os expoentes $\alpha$ e $\beta$ são ditos
\emph{expoentes críticos} e a função $\phi$ é chamada \emph{função
perfil}. Em alguns casos, a equação completa pode não ter tal
simetria mas, sob certas condições, ela ainda é
``assintoticamente'' simétrica e a solução do problema de valor
inicial se comporta, para tempos longos, como
\be
  \label{eq:compass}
  u(x,t)\approx \frac{A}{t^{\alpha}}\ \phi\left(\frac{x}{t^{\beta}}\right).
\ee
O prefator $A$ é supostamente não-nulo e pode depender da equação
diferencial em questão (dos parâmetros e dos termos não-lineares
interpretados como {\it perturbações}) e da condição inicial.
Contudo, os expoentes críticos e a função perfil são, muitas
vezes, \emph{universais},
isto é, eles não dependem nem do dado
inicial e nem das perturbações adicionadas à equação, desde que
escolhidas dentro de uma classe adequada.

Definimos \emph{classes de universalidade\label{ind:univ}} como conjuntos
de problemas que apresentam um mesmo comportamento assintótico.
Comportamentos da forma~(\ref{eq:compass}) são completamente determinados
por: expoente $\alpha$, expoente $\beta$, função de perfil $\phi$ e
prefator $A$. A classe de universalidade ser\'a definida em termos de
$\alpha$, $\beta$ e $\phi$.

Por exemplo, dizemos que duas condições iniciais estão na mesma classe de
universalidade se as soluções dos respectivos problemas de valor inicial
tiverem um comportamento assintóticoda forma~(\ref{eq:compass}) com os
mesmos expoentes $\alpha$ e $\beta$ e mesma fun\c c\~ao de perfil $\phi$.

Estendendo a noção de expoente crítico, diremos que $\alpha$ e
$\beta$ são \emph{expoentes críticos} se $u$ satisfizer \be
  \label{eq:long-time}
  u(x,t) \approx \frac{A}{g_1(t)}\phi\left( \frac{x}{g_2(t)}\right),
\ee
onde $g_1$ e $g_2$ são funções positivas
tais que ${1 \over \ln t} \ln g_1(t)
\to \alpha$ e ${1 \over \ln t} \ln g_2(t)
\to \beta$ quando $t\to\infty$.
A relação~(\ref{eq:long-time}) é mais precisamente
formulada através da existência do limite
\be
  \label{eq:long-time-lim}
  g_1(t)\ u( g_2(t)\ \cdot\ ,t) \tot A\ \phi(\cdot).
\ee
Essa definição engloba tanto o comportamento descrito
por~(\ref{eq:compass}) como
também outros casos de interesse, como, por exemplo,
\be
  \label{eq:compasslog}
  u(x,t) \approx \frac{1}{(t\ln t)^{\alpha}}\phi\left(\frac{x}{t^{\beta}}\right).
\ee

\section{Problemas estudados}

Neste trabalho estaremos interessados em estudar
o comportamento assintótico de
soluções de problemas de valor inicial difusivos.
Listamos abaixo os problemas de interesse.

\begin{prob}[Difusão Não-Linear com Coeficiente Periódico]
\label{prob:periodico}
$$
\left\{
\ba{l}
  u_t = (1 + \mu g(x))\,u_{xx} + \lambda F(u, u_x, u_{xx}) \\
  u(x,0) = f(x),
\ea
\right.
$$
onde $g$ é contínua e periódica\footnote{Sem perda de generalidade, vamos assumir
que $g$ tem período $2\pi$.}, $\mu\in\R$
é tal que $1+\mu g(x)>0 \ \forall\, x \in\R$ e $F$ é uma
perturbação. A condição inicial $f(x)$ é suave e de
suporte compacto.
\end{prob}

Se $\mu=0$ temos a chamada \emph{equação de difusão não-linear}.
Em particular, se $F(u, u_x, u_{xx})= -uu_x + \mbox{perturbações
de ordem superior}$, então a equação é conhecida como
\emph{equação de Burgers não-linear}. Ainda, se a perturbação
de ordem superior depender somente de $u$, isto é, se $F$ for da
forma $F(u, u_x, u_{xx})= -uu_x + g(u)$, então a equação é dita de
\emph{reação-difusão-convecção}.

No caso geral ($\mu\ne0$), temos um modelo
aproximado para o fenômeno da difusão
em
um meio composto e condutor, caracterizado por uma estrutura
microscópica
periódica, além das não-linearidades mencionadas acima.
Esse meio é \emph{não-homog\^eneo}, pois o coeficiente de difusão
$K(x)\equiv (1+ \mu g(x))$ não é constante com respeito \`a
variável $x$.

\begin{prob}[Difusão Não-Linear com Coeficiente Dependente do Tempo]
\label{prob:deptempo}
$$
\left\{
\ba{l}
  u_t = K (t)\, u_{xx} + \lambda F(u, u_x, u_{xx}) \\
  u(x,0) = f(x)
\ea
\right. ,
$$
onde
\be
\label{eq:c_t}
K(t)=t^p + o(t^p),
\ee
$p>0$, $o(s)$ sendo uma ordem pequena
de $s$, ou seja,
$\dlim_{t\to\oo}\frac{K(t)}{t^p} = 1$.
A condição inicial $f(x)$ é suave e de suporte compacto.
\end{prob}

Equações desse tipo, em que a difusividade cresce com o tempo, aparecem
na descrição de certos escoamentos multifásicos em meios porosos.  Entre
estes se incluem escoamentos associados a processos de recuperação de
hidrocarbonetos em reservatórios de petróleo, e de remediação e controle
de poluição em aquíferos \cite{peaceman77}. Nestes escoamentos, o
campo de velocidades é proporcional ao gradiente de pressão, segundo a
lei de Darcy, e o fator de proporcionalidade é uma propriedade geológica
do meio poroso, denominada ``permeabilidade'' (recíproca da resistividade
ao escoamento).  Devido à sua grande variabilidade como função da posição
no meio poroso, e à impossibilidade de uma caracterização precisa desta
variabilidade, a permeabilidade é comumente modelada como uma função
aleatória da posição (campo estocástico) \cite{dagan89}. Em conseqüência, as
equações que descrevem o escoamento possuem coeficientes estocásticos e,
portanto, suas soluções sao estocásticas. Ou seja, previsões baseadas
nestas equações são inerentemente estocásticas.

Nesse contexto, a equação acima descreve o escoamento esperado (médio,
no sentido estatístico) na presença de heterogeneidades geológicas
estocásticas \cite{furt-glim-lind-pere, glim-lind-pere-zhan}.
A difusividade crescente como função do tempo resulta da ocorrência
dessas heterogeneidades em todas as escalas de comprimento.

Em geral, o comportamento assintótico das soluções dos problemas
descritos acima é auto-similar -- isto é, da forma~(\ref{eq:compass}) --
ou quase auto-similar, com uma correção logarítmica -- isto é,
da forma~(\ref{eq:compasslog}) --  e possui
características bastante interessantes.
Enfatizamos algumas abaixo.

Os expoentes $\alpha$ e $\beta$ são universais: independem
da condição inicial $f$, da perturbação periódica $\mu g$,
do coeficiente $K(t)$ e
da perturbação não-linear $F$, desde que escolhidos dentre
uma classe adequada (nesta classe, $\alpha$ e $\beta$
dependem apenas do expoente $p$ em~(\ref{eq:c_t})).
A função de perfil
$\phi$ também é universal dentro desta classe:
no Problema~\ref{prob:periodico}
depende apenas da média harmônica do coeficiente de difusão,
enquanto no Problema~\ref{prob:deptempo} depende apenas
do expoente $p$.

No caso linear e homogêneo ($\mu=0=\lambda$ no
Problema~\ref{prob:periodico}; $\lambda=0$ no
Problema~\ref{prob:deptempo}), o prefator $A$ é universal,
dependendo apenas do que chamamos de \emph{massa} do dado inicial
$f$. Mais especificamente, temos $A=\int_\R f(y)dy$, ou seja, $A$
não depende de nenhum dos parâmetros da equação nem de detalhes em
pequena escala do dado inicial $f$. No caso não-linear,
este prefator não é universal: o que os teoremas ou os estudos
numéricos fazem é mostrar a existência de algum $A$ estritamente
positivo para o qual vale o limite~(\ref{eq:long-time-lim}) ou,
no máximo, determiná-lo.

\section{A técnica do grupo de renormalização}

A técnica que nos permite estudar os problemas aqui apresentados e
vários outros semelhantes chama-se \emph{a técnica do grupo de
renormalização}, que daqui em diante chamaremos de \emph{RG}\footnote{A
sigla RG vem da expressão \textsl{Renormalization Group}, do
inglês.}. Do ponto de vista de Equações Diferenciais, a técnica
foi desenvolvida no final dos anos $80$ e início dos anos $90$ por
Goldenfeld e colaboradores (veja \cite{gold92} e referências
lá citadas), em sua versão formal, e por
Bricmont~et~al.~\cite{bric-kupa-lin,bric-kupa}, em sua
versão rigorosa. A versão numérica foi introduzida por Chen e
Goldenfeld~\cite{chen}.

A aplicação formal desta técnica constitui-se, grosso modo, dos seguintes
passos. Primeiro, procuramos por uma solução da EDP que seja da
forma~(\ref{eq:sssol}). Se a EDP contém termos não-lineares, nós os
ignoramos ou ignoramos parte deles e consideramos apenas a parte da
equação que nos interessa. A suposição de que a equação possui solução da
forma~(\ref{eq:sssol}) é muito forte e leva à inevitável redução do grau
de liberdade do problema. Isso significa que reduzimos a EDP a uma EDO,
cuja incógnita é $\phi:\R\to\R$. Se formos capazes de encontrar tal
solução teremos a solução auto-similar desejada. O próximo passo é
descrever a \emph{classe de universalidade} desta solução, ou seja, o
conjunto de dados iniciais e perturbações para os quais o
comportamento assintótico é ditado por esta solução auto-similar.
Consideramos a mudança de escalas~(\ref{eq:rescale}) e
\emph{renormalizamos} a equação diferencial em questão, bastando para isto
obter a equação satisfeita por $v$. A parte do operador diferencial que
consideramos ao encontrar a solução auto-similar deve ser invariante por
esta renormalização, ao passo que as partes restantes, às quais chamamos
de \emph{perturbações}, em geral sofrem alterações. Classificam-se então,
formalmente, as perturbações como \emph{irrelevantes} no caso em que a sua
contribuição no operador diferencial tende a se anular com os sucessivos
processos de renormalização; como \emph{marginais} as perturbações que
permanecem invariantes ou assintoticamente invariantes pela
renormalização; como \emph{relevantes} as que tendem a crescer neste
procedimento. Repetimos esta mudança de escalas acompanhada da
renormalização da equação por várias vezes e esperamos verificar que, sob
certas condições, essas iterações resultem em alguma convergência. Por
isso, a técnica do Grupo de Renormalização é considerada uma \emph{técnica
em escalas múltiplas}.

A aplicação rigorosa do método segue uma idéia simples e elegante:
ao invés de tomarmos o limite $t\to\infty$ de uma só vez,
estabelecemos um procedimento iterativo em que o limite
$t\to\infty$ é visto como a composição de um operador consigo
mesmo um número infinito de vezes. Esse operador chama-se
\emph{transformação do grupo de renormalização} e a sua definição
depende da equação em questão. Os pontos fixos da parte linear do
operador são as soluções invariantes por escala. Determina-se
então a \emph{bacia de atração} de cada ponto fixo. Para condições
iniciais dentro da bacia de atração, o fluxo do RG fornecerá então
o comportamento para tempos longos.
Resumindo, o procedimento se reduz aos seguintes passos: dada uma
escala $L>1$, integram-se as equações, evoluindo o dado inicial do
tempo 1 ao tempo $L^2$; reescalonam-se a variável espacial $x$ e a
solução $u$ por fatores de $L$ (esses fatores dependem dos
expoentes críticos); subtrai-se da solução reescalonada a
contribuição vinda do ponto fixo, o resto supostamente é
``irrelevante" e espera-se que seja contraído à medida que
iterarmos esse procedimento; a solução reescalonada é então
utilizada como condição inicial para uma nova equação, similar à
original, dita \emph{equação renormalizada}; o procedimento é
iterado até que o regime assintótico seja atingido.

Um exemplo da aplicação do RG analítico é o estudo para a versão
homogênea do Problema~\ref{prob:periodico}, isto é, com $\mu=0$.
Bricmont~et~al.~\cite{bric-kupa-lin} mostraram que, para
$F$ da forma $F(u,u_x,u_{xx}) = u^au_x^bu_{xx}^c$ e para expoentes
$a$, $b$ e $c$ tais que $a+2b+3c-3>0$, se
a condição inicial $f$ for suficientemente pequena então vale
$$
  u(x,t) \approx \frac A {\sqrt t} \phi_* \left(
  \frac x {\sqrt t}\right)
  \quando
  ,
$$
onde $A$ é uma constante positiva e
\be
  \label{eq:gauss}
  \phi_*(x) = \frac 1 {\sqrt {4\pi}}\,\,\exp\left(-{ {x^2} \over 4}\right).
\ee O sentido do símbolo $\approx$ acima é o de um limite
da forma (\ref{eq:long-time-lim}), mas o resultado acima é
provado na topologia de um certo espaço de Banach cuja formulação
exata não nos interessa neste momento. O que chama a atenção é o
fato de que, dentro de um conjunto consideravelmente amplo, o
teorema fornece um critério extremamente simples que classifica
aquelas perturbações $F$ que não modificam o comportamento
assintótico da solução, em relação àquele da equação do calor
(caso em que, além de $\mu=0$, temos $\lambda=0$). Fazendo
$$
  d_F\equiv a +2b + 3c-3,
$$
Bricmont et~al.~\cite{bric-kupa-lin} classificam as perturbações
como \emph{irrelevantes} se $d_F>0$, \emph{marginais} se $d_F=0$ e
\emph{relevantes} se $d_F<0$. O
resultado ora mencionado será apresentado de forma precisa mais
adiante (veja Teorema~\ref{teo:nldif} no
Capítulo~\ref{chap:rgana}). Observe também que o teorema
nos diz que os expoentes críticos
$\alpha$ e $\beta$ valem $1/2$ e que a função perfil é a distribuição
gaussiana, seja qual for a perturbação irrelevante adicionada à
equação linear.

No Problema~\ref{prob:deptempo}, a classificação formal das perturbações é
semelhante à do Problema~\ref{prob:periodico}, mas no lugar de $d_F$ temos
\be
  \label{eq:eta}
  \eta_F = a+2b+3c - \frac{p+3}{p+1}.
\ee Diremos que o termo não-linear é \emph{irrelevante} se $
\eta_F > 0 $; \emph{marginal} se $ \eta_F= 0$; \emph{relevante} se
$ \eta_F< 0$.

A primeira novidade deste problema é que qualquer $p>0$ faz com que as
perturbações da forma $u^3$, $uu_x$ e $u_{xx}$ deixem de ser marginais
para se tornarem irrelevantes (repare que $d_F=0$ para essas três
perturbações, ou seja, no caso do Problema~\ref{prob:periodico} essas
perturbações são marginais).

Usando as mesmas idéias introduzidas originalmente por Bricmont
et~al.~\cite{bric-kupa-lin}, Braga et~al.~\cite{braga04} provaram
que as soluções do Problema~\ref{prob:deptempo} linear, isto é, com
$\lambda=0$, se comportam como
\be
  \label{eq:assint-deptempo}
  u(x,t) \approx \frac {A_f} {\sqrt {t^{p+1}}} \,\, \phi_p \left( \frac x {\sqrt
  {t^{p+1}}}\right)\quando,
\ee onde
\be
  \label{eq:phi_p}
  \phi_p(x)=\sqrt{\frac{p+1}{4\pi}}\,\,\exp\left({-\frac{p+1}{4}\ x^2}
  \right)
\ee
e $A_f = \int_\R f(y)dy$ (veja o Teorema~\ref{teo:compu-0} no
Capítulo~\ref{chap:rgana}).

Braga et~al.~\cite{braga04} conjecturaram também que qualquer
perturbação não-linear irrelevante (isto é, com $\eta_F>0$) está
na mesma classe de universalidade do problema linear,
caracterizada pelos expoentes críticos $\alpha=\beta=(p +1)/ 2$ e
pela distribuição gaussiana $\phi_p$ dada em (\ref{eq:phi_p}) --
veja, no Capítulo~\ref{chap:rgana}, a Conjectura~\ref{conj:compassgeral}.

\section{A técnica de homogeneização}

O Problema~\ref{prob:periodico} foi apresentado como um modelo
aproximado de
um meio composto e condutor, caracterizado por uma estrutura
microscópica periódica, ou seja,
um meio \emph{não-homog\^eneo} cujo coeficiente de difusão
$K(x)=1+\mu g(x)$ é periódico de período $2\pi$.
É razoável supor que as
variações locais de $K(x)$ não são perceptíveis se observarmos esse
meio a uma escala $L\gg 2\pi$, de forma que, em grandes escalas, a
heterogeneidade do meio seja cada vez menos perceptível, se
aproximando de um meio homogêneo (isto é, com $K$ constante).
Por isso, quando nos concentramos apenas no termo $K(x)$ do
Problema~\ref{prob:periodico}, temos um \emph{problema em duas escalas}.
A teoria de homogeneização \cite{papa78,
jikov94,tartar77} para o
Problema~\ref{prob:periodico} (com
$\lambda=0$, isto é, sem a parte não-linear) dá uma resposta rigorosa
ao argumento heurístico aqui apresentado. Mais que isso, nos diz que a
solução $u$ se comporta como a solução de
\begin{equation}
\label{eq:efet}
  \left\{
    \ba{l} u_t = \sigma u_{xx} \\ u(x,1) = f(x), \ea
  \right.
\end{equation}
onde $\sigma = \langle K^{-1} \rangle ^{-1}$, sendo que $\langle
\cdot \rangle$ denota a média sobre o período de uma função periódica
qualquer (essa notação será mantida em todo este texto).

A substituição de $K(x)$ pela sua média harmônica
no caso unidimensional possui uma explicação heurística simples.
Primeiro note que o coeficiente
de difusão representa a facilidade de propagação (de calor,
do local mais quente para o mais frio; de partículas, do local
de maior para o de menor concentração; ou de cargas elétricas,
do local de maior para o de menor potencial elétrico) enquanto sua
recíproca representa dificuldade.
Agora considere, por exemplo, um sistema elétrico constituído de
diferentes condutores em série.
Sabemos, da teoria elementar de circuitos elétricos~\cite{halliday},
que se quisermos substituir esses condutores por outros
idênticos entre si (isto é, com uma condutância \emph{homogênea}),
então o valor dessa condutância deve ser dado pela média harmônica das
condutâncias originais. O mesmo fenômeno ocorre na difusão de calor:
podemos imaginar o meio condutor não-homogêneo como sendo uma superposição em
série de inúmeros condutores de espessura $\Delta x$ e, se quisermos
substituí-los por um condutor somente (ou vários com
mesma espessura $\Delta x$ e idênticos entre si),  temos que considerar a
média harmônica desses condutores infinitesimais.

No Capítulo~\ref{chap:homog} desta dissertação vamos estudar
rigorosamente alguns problemas de ho\-mo\-ge\-nei\-za\-ção. Em particular,
o Corolário~\ref{cor:parab} justifica a substituição de $K(x)$
por $\sigma$ em (\ref{eq:efet}).

A técnica de homogeneização foi desenvolvida para estudar o
comportamento assintótico de soluções de equações diferenciais
parciais com coeficientes periódicos.  Problemas dessa natureza t\^em
duas escalas de comprimento típicas: a escala do meio difusivo,
associada ao período do coeficiente
de difusão, e a escala de observação,
associada ao ponto espacial e ao tempo no qual se fazem medidas
experimentais. Por outro lado, a técnica do grupo de renormalização
\cite{bric-kupa,bric-kupa-lin,gold92} foi
desenvolvida para
tratar problemas com um número infinito de escalas, como por exemplo a
equação de difusão não-linear (Problema~\ref{prob:periodico} com $\mu=0$).
Se substituirmos $K(x)$ por $\sigma$ no Problema~\ref{prob:periodico},
podemos aplicar o teorema de Bricmont et~al.~\cite{bric-kupa-lin}
mencionado acima. Com esse argumento Braga et~al.~\cite{braga03a}
conjecturaram que o mesmo resultado continua válido sem a hipótese de
homogeneidade para o coeficiente de difusão ($\mu = 0$), exceto
porque a função de perfil é \emph{renormalizada}. Em outras
palavras, se $F$ pode ser escrita na forma $F(u,u_x,u_{xx})=
u^au_x^bu_{xx}^c$, se $d_F>0$ e se $\lambda$ é suficientemente pequeno,
então
\be
  \label{eq:assint-periodico}
  u(x,t) \approx \frac A {\sqrt t}\ \phi_\sigma \left(
  \frac x {\sqrt t}\right) \quando ,
\ee
onde
\be
  \label{eq:phi_sigma}
  \phi_\sigma(x) =
\frac 1 {\sqrt {4\pi\sigma}}\,\,\exp\left(-{{\ x^2} \over 4\sigma}\right).
\ee

Para uma formulação precisa de~(\ref{eq:assint-periodico}),
veja a Conjectura~\ref{conj:oscilante} no
Capítulo~\ref{chap:rgana}.

O termo $\mu g(x)$ tem representado um obstáculo para o estudo rigoroso do
Problema~\ref{prob:periodico} usando-se o grupo de renormalização
analítico. Até onde sabemos, esta é a primeira vez em que se estabelece
uma conexão entre a técnica de homogeneização e o grupo de renormalização
para EDP's (no contexto de física estatística, veja o trabalho de
Naddaf e Spencer~\cite{naddaf97}).

Duro e Zuazua~\cite{zuazua2000} consideraram a equação
$$
  u_t -
  \mbox{div}(a(\vec{x})\nabla u) = \vec{d}\cdot\nabla(|u|^{q-1}u)
  \,\,\,\,\,
  \mbox{em} \,\,\,\,\, {\R}^N\times(0,\infty) ,
$$
onde $\vec{d}\in{\R}^N$ e $a(\vec{x})$ é uma matriz suave, simétrica,
periódica e uniformemente elíptica. Usando a teoria de homogeneização eles
provaram, dentre outros resultados, que a classe de universalidade do
comportamento auto-similar da equação linear ($\vec d = 0$) inclui
perturbações com $q>1 + 1/N$. Apesar de não provarmos rigorosamente,
nossos resultados numéricos indicam que o
resultado de Duro e Zuazua deve continuar válido
para outras perturbações não-lineares\footnote {Chamamos atenção para uma
sutil diferença entre a natureza das equações. O problema estudado por
Duro e Zuazua, quando restrito ao caso unidimensional, inclui um termo da
forma $(K(x)u_x)_x$ enquanto o problema que estudamos computacionalmente
neste trabalho tem o termo da forma $(K(x)u_{xx})$.}.

\section{O grupo de renormalização numérico}

Partindo exatamente das mesmas idéias em que se baseiam o RG formal
e o RG analítico, surge naturalmente uma versão numérica
do mesmo procedimento, o chamado \emph{grupo de renormalização
numérico}, ou simplesmente \emph{RG numérico}. Até onde sabemos,
a primeira versão deste método
foi proposta por Chen e Goldenfeld em 1995~\cite{chen}. O método
numérico consiste, assim como no caso analítico, em evoluir a
condição inicial segundo o fluxo dado pela EDP e depois fazer uma
mudança de escalas. A primeira vantagem do RG numérico,
obviamente, é que não
precisamos saber resolver a equação analiticamente ou restringir nossas
hipóteses aos dados e parâmetros ``suficientemente pequenos'',
uma vez que a mesma será resolvida numericamente. Outra diferença
está no cálculo dos expoentes $\alpha$ e $\beta$: no algoritmo
numérico que vamos utilizar, os expoentes críticos são calculados
dinamicamente, ao fim de cada iteração do algoritmo, mesmo que
não sejam conhecidos \textsl{a priori}\footnote
{Isso é particularmente útil no caso de equações cujos expoentes
críticos não podem ser obtidos a partir de simetrias da
EDP, como no caso da equação de Barenblatt. Neste caso dizemos
que o problema possui \emph{expoentes críticos anômalos}.}.
Por último, o RG numérico permite calcular prefatores nos
casos em que os teoremas se limitam a garantir sua existência.

Nesta dissertação vamos apresentar a versão do RG numérico que foi
implementada e utilizada nos trabalhos feitos em colaboração com G.~Braga,
F.~Furtado e J.~Moreira. Nesses trabalhos estudamos numericamente o
comportamento assintótico de soluções dos problemas já apresentados e
determinamos os expoentes críticos, pré-fatores e funções de perfil, assim
como tentamos caracterizar suas classes de universalidade. Essa versão do
RG numérico, que considera a renormalização das equações e permite o
cálculo dos prefatores, além de uma malha discreta fixa, foi introduzida
originalmente por G.~Braga, F.~Furtado e V.~Isaia -- veja a tese de
doutorado de Isaia~\cite{isaia02}.

Os problemas estudados nesta disserta\c c\~ao s\~ao todos em apenas uma
dimens\~ao espacial. Ressaltamos, entretanto, que tal restri\c c\~ao n\~ao
\'e, de forma alguma, uma restri\c c\~ao intr\'\i nseca ao m\'etodo do
grupo de renormaliza\c c\~ao num\'erico. Veja, novamente, a tese de
doutorado de Isaia~\cite{isaia02} para o estudo de problemas com mais
dimens\c c\~oes espaciais.

Para o Problema~\ref{prob:periodico}, a verificação
do comportamento assintótico (\ref{eq:assint-periodico}),
incluindo o estudo de um caso particular de perturbações
relevantes,
bem como a descrição e a validação do método numérico utilizado,
foram publicados
na revista \textsl{Multiscale Modeling and Simulation}~\cite{
braga03a}, da SIAM.

Para o Problema~\ref{prob:deptempo}, além de verificarmos
o comportamento (\ref{eq:assint-deptempo}), fizemos o estudo
do comportamento próximo da linha que separa perturbações
relevantes e irrelevantes.
Neste caso,
nos restringimos ao caso particular da perturbação
não-linear
$$
  \lambda F = -u^a,
$$
onde $a\in(1,\oo)$, e trabalhamos no espaço de parâmetros $(p,a)$.
A condição para que o termo $-u^a$ seja marginal ($\eta_F = 0$)
se reduz a
\be
  \label{eq:a-critico}
  a = \frac{p+3}{p+1},
\ee
 cujo gráfico no plano $(p,a)$ é o que chamamos de \emph{curva crítica}
(veja a Figura~\ref{fig:pht}). Estudando numericamente o comportamento de
soluções próximas a essa curva e analisando o expoente crítico $\alpha$,
verificamos que a mesma é de fato a interface entre as regiões relevante e
irrelevante e dizemos que o sistema sofre uma \emph{transição de fase}
sobre esta curva. Estes resultados estão sendo redigidos na forma
de artigo~\cite{braga04} a ser submetido para publicação.

Partes do conte\'udo desta disserta\c c\~ao e da disserta\c c\~ao de
mestrado de Moreira~\cite{moreira} est\~ao tamb\'em publicadas nas
notas~\cite{braga03b} do mini-curso oferecido no $57^{\rm o}$ Semin\'ario
Brasileiro de An\'alise.

Todos os casos aqui mencionados são discutidos no Capítulo~\ref
{chap:rgana} e os resultados numéricos são apresentados no Capítulo~\ref
{chap:resul}.

\section{Considerações finais}

Este trabalho está dividido da seguinte maneira: no
Capítulo~\ref{chap:rgana} estudamos formalmente os
Pro\-ble\-mas~\ref{prob:periodico}~e~\ref{prob:deptempo} através do grupo
de re\-nor\-ma\-li\-za\-ção. Antes disso, as idéias e a motivação do
método serão introduzidas para a equação do calor. No
Capítulo~\ref{chap:rgnum} descrevemos uma versão numérica do grupo de
renormalização. O método numérico consiste na integração da equação por um
curto intervalo de tempo, seguido de um reescalonamento. O mesmo fornece
um quadro detalhado do regime assintótico: função perfil, expoentes
críticos, pré-fatores e parâmetros efetivos. O método é aplicável
a problemas com soluções assintoticamente auto-similares, inclusive no
caso crítico em que aparece correção logarítmica (decaimento com $(t\ln
t)^\alpha$). O Capítulo~\ref{chap:resul} contém, principalmente, a
verificação numérica das conjecturas e dos teoremas apresentados no
Capítulo~\ref{chap:rgana}, com ênfase na verificação de
(\ref{eq:assint-periodico}) e (\ref{eq:assint-deptempo}), além da
transição de fase em (\ref{eq:a-critico}). No Capítulo~\ref{chap:homog}
vamos apresentar a teoria de homogeneização para equações
lineares. Um dos objetivos é obter resultados rigorosos que justifiquem a
substituição do coeficiente não-homogêneo por um homogêneo no estudo do
Problema~\ref{prob:periodico}. Pretendemos apresentar ao leitor
um texto auto-contido sobre a teoria de homogeneização, para equações
elípticas e parabólicas lineares, em espaços de Sobolev. O leitor deve
estar familiarizado com a Teoria de Medida e Integração de Lebesgue e,
caso não tenha feito um curso de Análise Funcional, deverá ter um livro
introdutório à mão (por exemplo, \cite{oliveira01}) para buscar
alguns conceitos e resultados básicos.

Listamos aqui alguns problemas que não são tratados neste texto para os
quais acreditamos ser possível -- e interessante -- fazer um estudo
numérico cuidadoso. Alguns deles chegaram a ser estudados pelo autor e
colaboradores, mas não foram concluídos. No Problema~\ref{prob:periodico}
pode-se substituir o termo difusivo por um que seja conservativo, a
saber, $u_t = (K(x)u_x)_x$.  O estudo da equação $u_t = [K(x)(u^m)_x]_x$,
motivada pela difusão em meios porosos não-homogêneos, é também de grande
interesse; o autor desconhece qualquer estudo formal, numérico ou
analítico acerca da homogeneização desta equação. Também o estudo da
transição de fase no espaço de parâmetros $(m,a)$ da equação $u_t =
(u^m)_{xx} - u^a$ é um representante de todo um conjunto de problemas (em
que se enquadra o Problema~\ref{prob:deptempo} desta dissertação) em que
um parâmetro do termo difusivo determina o expoente crítico de decaimento
e há uma competição entre o regime linear e não-linear, dando origem a
uma transição de fase. Em todos os problemas em que ocorre essa transição
de fase, mesmo que com apenas um parâmetro (como $u_t = u_{xx} - u^a$), o
estudo numérico fino das corre\c c\~oes logar\'\i tmicas que surgem nos
pontos críticos é também questão interessante na opinião do autor. Em
outra direção ortogonal, podemos pensar em adaptar as idéias do método
para o estudo de equações (não necessariamente hiperbólicas) para as
quais a solução se comporta assintoticamente como uma ``onda viajante''.
Para este último caso o autor chegou a implementar uma versão que
forneceu resultados preliminares bastante satisfatórios. Por último, a
formulação de um método mais geral, que seja apropriado para estudar
problemas que apresentam decaimento, espalhamento e transla\c c\~ao em
seu comportamento assintótico, é questão de grande importância e alguns
dos colaboradores do autor estão trabalhando nessa direção.

Por último, observamos que decidimos definir $\phi_*$, $\phi_\sigma$ e
$\phi_p$ por (\ref{eq:gauss}), (\ref{eq:phi_sigma}) e (\ref{eq:phi_p}),
respectivamente, para manter o texto o mais limpo possível. Com esta
notação, o sub-índice não é um parâmetro propriamente dito; por exemplo,
$\phi_* = \phi_\sigma = \phi_p$ se $\sigma=1$ e $p=0$; assim, não faz
sentido escrever $\phi_{3,14}$.

%% file: rgana.tex
\chapter{O Grupo de Renormalização Formal e Analítico}
\label{chap:rgana}

Neste capítulo vamos estudar,
usando o método do Grupo de Renormalização -- RG,
os problemas apresentados na Introdução. As idéias do método, assim como
sua motivação, são introduzidas de forma natural para a equação do calor,
devido às suas simetrias. Por essa razão a Seção~\ref{sec:rgcalor} é
dedicada ao estudo desta equação.
Nessa seção nós introduzimos a operação por
mudança de escalas e mostramos que a equação do calor tem uma
solução simétrica com respeito a esta operação. Essa solução
é utilizada para provarmos a Proposição \ref{prop:compasseqcalor},
sobre o comportamento assintótico de soluções do
problema de valor inicial associado.

Motivados pela operação de mudança de escalas, na Seção \ref{sec:tran-grup-reno}
nós definimos o operador RG. Vamos utilizar alguns resultados
básicos sobre a transformada de Fourier aplicada à equação do calor. Nossa preocupação
principal não está no preenchimento detalhado das definições e demonstrações, mas na
idéia do RG como uma ferramenta para se provar teoremas a respeito do comportamento
assintótico de soluções de equações difusivas. Por essa mesma razão, o leitor que não
esteja familiarizado com esses resultados pode ler o texto sem se preocupar com tais
detalhes.
O nosso objetivo é apresentar uma outra
prova da Proposição \ref{prop:compasseqcalor}, utilizando, para isto,
o reescalonamento da solução do pvi
iterativamente.  Os ingredientes básicos
para  a realização desse procedimento são:
a Propriedade de Semigrupo~\ref{teo:semi-grupo}
e o Lema da Contração~\ref{teo:contr}.

Na Seção~\ref{sec:rggeral}, os Problemas \ref{prob:periodico} e
\ref{prob:deptempo} da Introdução
serão considerados  sob o ponto de vista formal e
heurístico. Observamos que a aplicação do método
descrito acima depende fortemente
da invariância da equação do calor por uma mudança de escalas. Contudo,
em geral as equações não são invariantes sob este reescalonamento.
Nesse caso, para analisar os Problemas~\ref{prob:periodico} e
\ref{prob:deptempo} sob o ponto de vista do RG, nós consideraremos
uma equação diferente em cada escala. Dizemos que a equação foi
{\em renormalizada}.
No processo de renormalização das equações, esperamos que fique latente a
simplicidade e o poder da mudança de escalas: é possível determinar
toda uma classe de perturbações não-lineares $F$, para as quais o
comportamento assintótico das soluções continua sendo o mesmo da equação
linear.

\section{A equação do calor e mudança de escalas}
\label{sec:rgcalor}

{Nesta seção vamos considerar o problema de valor inicial
para a equação do calor:
\be
\label{eq:calor}
  \left\{ \ba{cl}
  u_t = u_{xx}, & \; x \in \R,\ t>0 \\
  u(\cdot,0) = f, & \; f \in C_0(\R),
  \ea \right.
\ee
onde $C_0(\R)$ é o conjunto das funções contínuas
de suporte compacto em $\R$.

Enunciamos, a seguir, o teorema que nos garante a existência (global) e unicidade do pvi
(\ref{eq:calor}), como também nos fornece uma representação integral para a solução. Esse
teorema pode ser formulado sob condições mais fracas,
mas para os nossos propósitos é suficiente o caso em que o dado inicial seja
contínuo de suporte compacto. A prova da proposição abaixo pode ser encontrada em
\cite{iorio}.

\begin{prop}
\label{exis-unic}
O problema de valor inicial \ref{eq:calor},
com dado inicial $f\in C_0(\mathbb{R})$,
possui uma \'unica solução $u(\cdot, t)\in C^\infty(\R)$, $t>0$.
Essa solução \'e dada pela f\'ormula
\be
  \label{eq:sol}
  u(x,t) = \frac 1 {\sqrt{ 4 \pi t}} \int_\R {e^{-(x-y)^2 \over 4t}
  } f(y) dy
\ee
para qualquer $t>0$ e $u(x,t)\to f(x)$ quando $t\downarrow0$.
\end{prop}
}

Seja $L$ um n\'umero real estritamente positivo e
façamos a seguinte mudança de escalas:
\be
  \label{eq:mudanca}
   x \mapsto Lx \quad t \mapsto L^2t \quad u \mapsto Lu.
\ee
Chamemos de $v(x,t)$ \`a função $u(x,t)$ devidamente
reescalonada, isto \'e,
\be
  \label{eq:mudanca.1}
  v(x,t) = L u (L x , L^2 t).
\ee
Temos que $v_t = L^3 u_s $ e $ v_{xx} = L^3 u_{yy}$, onde $s=L^2t$ e
$y=Lx$, ou seja,  $u_t=u_{xx}$ se e somente se $v_t=v_{xx}$.
Dizemos que a equação do
calor \'e \emph{invariante pela mudança de escalas} (\ref{eq:mudanca.1}), para
qualquer $L>0$.

Vamos procurar soluções que satisfaçam $u(x,t) = Lu(Lx,L^2t)$.
Substituindo formalmente $L$ por $1 / {\sqrt t}$,
temos
$
   u(x,t) = \frac 1 {\sqrt t}\ u(x / {\sqrt t},1),
$
o que nos motiva a procurar soluções da forma
\be
\label{eq:invar.1}
   u(x,t) = {1 \over \sqrt t}\ \phi\left({x \over \sqrt
t}\right).
\ee

Qualquer função $u$ que se escreve na forma~(\ref{eq:invar.1}) é
obviamente invariante pela mudança de escalas~(\ref{eq:mudanca.1}).
Para determinar a função $\phi(\cdot)$, n\'os substituímos
(\ref{eq:invar.1})
em (\ref{eq:calor}), obtendo uma
equação diferencial ordin\'aria para $\phi(\xi)$:
\be \label{eq:ord}
  0 = \phi'' + \frac 1 2 \xi \phi' + \frac 1 2 \phi
  = \left(\phi' + \frac { \xi \phi} 2\right)'.
\ee

A equação acima admite soluç\~oes da forma
$C\exp(-\xi^2/4)$ e, escolhendo a constante $C$ de modo que $\int_\R
\phi(\xi)d\xi=1$, obtemos a distribuição gaussiana $\phi_*$ dada
por~(\ref{eq:gauss}).
Substituindo $\phi$ por $\phi_*$
em~(\ref{eq:invar.1}), temos
uma \emph{solução invariante pela mudança de escalas},
que satisfaz a
seguinte relação: $\sqrt t\, u(\sqrt t\ x,t) = \phi_*(x)$.
A proposição abaixo mostra como
essa relação est\'a
relacionada com
o comportamento assintótico da solução
geral.

\begin{prop}
\label{prop:compasseqcalor}
  Seja $u(x,t)$ a solução da equação do calor com dado
inicial $f\in  C_0(\R)$. Então
  \be
     \label{eq:limitepontual}
     \sqrt t \ u (\sqrt t\ x, t) \tot A_f \phi_*(x),
     \quad \forall\ x\in\R,
  \ee
  onde $A_f = \int_\R f(y)dy$.
\end{prop}

\prova
  Substituindo $x$ por $x \sqrt t $ e multiplicando
  (\ref{eq:sol}) por $\sqrt t$, obtemos
  \begin{eqnarray}
     \sqrt t\ u (\sqrt t\ x, t) - A_f\phi_*(x)
     &=& \frac 1 {\sqrt{4 \pi}} \int_\R
     e^{-(x\sqrt t - y)^2 \over 4 t}       \nonumber
     f(y) dy - \phi_*(x) \int_\R f(y)dy
     \\ \label{eq:dominada}
     &=& \int_\R
     \left(\frac 1 {\sqrt{4 \pi}}\ e^{-(x\sqrt t - y)^2 \over 4 t}
     - \phi_*(x) \right) f(y)dy
     \\ \label{eq:convintegr}
     &=& \int_{\R} \phi_*(x)\left(e^{2 x y
     \sqrt t -y^2 \over 4 t}-1\right) f(y) dy.
  \end{eqnarray}

  Primeiro note que o integrando de (\ref{eq:dominada}) é
  dominado em módulo por $|f(y)|/\sqrt \pi$, que é integrável.
  Note também que o integrando de (\ref{eq:convintegr})
  converge a zero quando $t\to\oo$, isso para todo
  $y\in\R$ (podemos considerar $x$
  como um parâmetro fixo). Aplica-se o Teorema da Convergência
  Dominada de Lebesgue para se obter (\ref{eq:limitepontual}).
\eop

\section{A transformação do grupo de renormalização}
\label{sec:tran-grup-reno}

O limite na Proposição~\ref{prop:compasseqcalor} diz
que soluções da equação do calor t\^em o
comportamento assint\'otico ditado pela solução invariante por mudança
de escalas. Vimos anteriormente que a equação do calor \'e invariante por
essa
mudança. Isso nos motiva a definir o operador do grupo de
renormalização (veja
\cite{bric-kupa,bric-kupa-lin}). Esse operador age num
espaço de funç\~oes, mais especificamente no espaço dos dados
iniciais, e o mesmo consiste em evoluir, no tempo, o dado inicial $f(x)$,
seguindo o fluxo dado pela equação do calor (essa evolução \'e
realizada do tempo $t=1$ ao tempo $t=L^2$) e, em seguida,
mudamos a escala do domínio e do contra-domínio do dado
evoluído. Denotando esse operador por $R_L$, onde $L>1$ \'e o fator de
escala\footnote{Daqui em diante vamos sempre assumir que $L>1$.},
podemos descrever essa seq\"u\^encia de operações pelo diagrama:
\be
    f = u_f(\cdot,1)
  \mapsto
    u_f(\cdot,L^2)
  \mapsto
    L u_f(L\,\cdot,L^2) = R_Lf,
\ee
onde $u_f$ denota a solução da equação do calor com condição inicial
$u(x,1) = f(x)$. De outra forma:
\be
  \label{eq:rg}
  (R_Lf)(x)\equiv L u_f(Lx,L^2).
\ee

É claro que $R_L$ \'e um operador linear, uma vez que a equação do
calor \'e linear e homogênea.

Vamos estudar como $R_L$ atua conjugado com a transformada de Fourier
(veja a definição e principais propriedades dessa transformada
no Ap\^endice~\ref{sec:tf}).

No caso da equação do calor, para a qual já sabemos da existência e unicidade
de uma solução global, a Propriedade de Semigrupo,
que enunciamos e provamos abaixo,
segue trivialmente da definição de $R_L$.
Entretanto, vamos usar a transformada de Fourier para
ilustrar um outro aspecto do RG analítico: podemos provar a existência da
solução global do problema a partir da existência local.

Primeiramente, se $u$ denota a única solução da equação do calor
com condição inicial $u(\cdot,0)=f$, e $u_f$ denota a única solução
com condição inicial $u_f(\cdot,1)=f$,
então temos
\be
  \label{eq:soltransf.1}
  \label{eq:soltransf}
  \hat u(k,t) =  \ e^{-k^2t} \hat f(k) \quad \mbox{e} \quad
  \hat u_f(k,t) =  \ e^{-k^2(t-1)} \hat f(k).
\ee
Para uma prova deste fato, veja \cite{iorio}.
Aplicando (\ref{eq:proptf}) a (\ref{eq:rg}), e substituindo
(\ref{eq:soltransf}), temos
\be \label{eq:rgfour}
  \widehat{(R_L f)} (k) = \widehat{(Lu_f(L\cdot,L^2))}(k) =
  \hat u_f(k/ L,L^2) =  e^{-(\frac{k}{L})^2(L^2 - 1)}
  \hat f(k / L).
\ee
Uma propriedade fundamental do operador $R_L$ \'e a seguinte:
\begin{teo}[Propriedade de Semigrupo]
\label{teo:semi-grupo}
  Sejam $L_1>1$ e $L_2>1$. Então
  \be \label{eq:semi-grupo}
    R_{L_1} \circ R_{L_2} = R_{L_1 L_2}.
  \ee
\end{teo}
\prova
  Basta desenvolver (\ref{eq:rgfour})
  para mostrar que $\widehat{(R_{L_1} R_{L_2} f)}
  = \widehat{(R_{L_1L_2} f)}$ e o resultado segue, como vemos:
\begin{eqnarray*}
\widehat{(R_{L_1} R_{L_2} f)\ }\ (k) &=&
  e^{-( k / {L_1})^2(L_1^2-1)}
  \widehat{(R_{L_2}f)}(k / {L_1}) \\
&=&
  e^{-(\frac k {L_1})^2(L_1^2-1)}
  e^{-(\frac k {L_1 L_2})^2(L_2^2-1)}
  \hat f( k / {L_1 L_2}) \\
&=&
  e^{ -k^2(1 - \frac 1 {L_1^2}) }
  e^{ -k^2(\frac 1 {L_1^2} - \frac 1 {(L_1L_2)^2}) }
  \hat f( k / {L_1 L_2}) \\
&=&
  \widehat{(R_{L_1L_2}f) \ } \ (k).
  \end{eqnarray*}
A prova está concluída.\eop

Também podemos ``diagonalizar'' o operador $R_L$, que
tem os polinômios de Hermite com peso gaussiano como auto-vetores no
espaço de Fourier.

A transformada de Fourier da gaussiana $\phi_*(x)$
dada por (\ref{eq:gauss}) é $\widehat {\phi_*}(k) = e^{-k^2}$.
De (\ref{eq:soltransf}) temos
\be
  (\widehat{R_L \phi_*})  (k) = e^{-k^2(1 - \frac 1 {L^2})}
  e^{-\frac {k^2}{L^2}} = \widehat {\phi_*}(k),
\ee
ou seja, a distribuição gaussiana $\phi_*$ \'e um ponto
fixo do RG.
Da mesma forma, se
$f_n(x)$ \'e a função cuja transformada de Fourier \'e $\widehat {f_n}(k) =
k^ne^{-k^2}$ (esses são, essencialmente, os
polin\^omios de Hermite com peso gaussiano), então $f_n$ \'e um
autovetor de $R_L$ com autovalor $(1/L)^{n}$.
Dizemos que a gaussiana
$\phi_*(x)$ est\'a na direção \emph{marginal} e que os outros autovetores
estão na direção \emph{irrelevante}. O subespaço gerado pelos
autovetores irrelevantes tamb\'em \'e chamado de \emph{variedade
estável}.

Definimos o espaço $\B_4$ por
$$
  \B_4 = \left\{ f(x) \in L^2(\R): \hat f \in C^1(\R) \mbox{ e }
  \|f\|_{\B_4}<+\oo \right\},\ \ {\rm onde}
$$
\[
  \|f\|_{\B_4} = \sup_{k\in\R} \left[
  (1 + k^4) ( |\hat f(k)| + |\hat f'(k)| )
  \right].
\]

Pode-se mostrar que $\B_4$ é um espaço de Banach,
que os elementos de $\B_4$ são funções contínuas e
que a norma $\|\cdot\|_{\B_4}$ domina as normas $L^1(\R)$,
$L^2(\R)$ e $L^\oo(\R)$ (veja, por exemplo, a tese de mestrado de
Moreira~\cite{moreira} ou as notas de mini-curso~\cite{braga03b}).

\begin{teo}[Lema de Contração]
  \label{teo:contr}
Se $f\in\B_4$ e $u(x,t)$ denota a solução da
equação do calor tal que $u(\cdot,0)=f$,
então $u(\cdot,t)\in\B_4\ \forall \ t>0$.
Em particular,
$R_L \B_4 \cont \B_4$.

Além disso, se $\hat{f}(0)=0$,
vale
  $$
    \|R_L f\|_{\B_4} \leq \frac{C}{L}\|f\|_{\B_4},
  $$
onde $C=C(L)>0$ é decrescente em relação a $L$ e não depende de $f$.
\end{teo}

\prova Começaremos pela segunda parte. A equação do calor pode
ser vista como um caso particular do Problema~\ref{prob:deptempo}
com $p=\lambda=0$ e $K(t)=1$. A prova que apresentamos aqui é uma
adaptação para este caso particular de um prova mais geral -- veja
o Lema de Contração em \cite{braga04}. A partir de
(\ref{eq:rgfour}) temos $
  \widehat{R_L f}(k) =
  e^{-(\frac{k}{L})^2(L^2 - 1)} \hat f(k / L)
$
e
$
  \widehat{R_L f}'(k) =
  -(L^2-1)(2k/L^2) e^{-(\frac{k}{L})^2(L^2 - 1)} \hat f(k / L) +
      L^{-1}{e^{-(\frac{k}{L})^2(L^2 - 1)}} \ \hat f'(k / L) .
$
Logo, vale que
$$
  { |\widehat {R_Lf}(k)| + |\widehat {R_Lf}'(k)|
  \over e^{-(\frac{k}{L})^2(L^2 - 1)} }
  \leq
   \left|\hat f(k / L)\right| + \left| 2k \hat f(k / L) \right|
   + \left| \frac 1 L \hat f'(k/L) \right|.
$$
Valem as estimativas $\hat f'(k/L) \leq \|f\|_{\B_4}$ e
\be
  \label{eq:estmfkl}
  |\hat f(k/L)| = \left| \hat f(0) + \int_0^{k/L} \hat f'(s)ds \right|
  \leq \int_0^{|k/L|}   \frac{\|f\|_{\B_4}}{1 + s^4}\ ds \leq
  \left|  k / L \right| \|f\|_{\B_4},
\ee
uma vez que $\hat f(0)=0$. 
Assim, temos
$$
  (1+k^4)(|\widehat {R_Lf}(k)| + |\widehat {R_Lf}'(k)|) \leq
  { (1+k^4)\left(1 + |k| + 2 |k^2| \right)
  e^{-(\frac{k}{L})^2(L^2 - 1)}
  \over L }\ \|f\|_{\B_4}.
$$
Definindo $C$ como
$$
  C(L) = 
  \sup_{k\in\R} \left\{
  (1+k^4)\left(1 + |k| + 2|k^2| \right)
  e^{-(\frac{k}{L})^2(L^2 - 1)} \right\} < +\oo, \quad L>1,
$$
temos que $\|R_L f\|_{\B_4} \leq \frac{C}{L}\|f\|_{\B_4}$.
Além disso, $C(L)$ é positivo e decrescente, logo $C/L$ é estritamente
decrescente e $\lim_{L\to\oo}C/L = 0$.

A primeira parte pode ser mostrada de forma idêntica,
exceto porque (\ref{eq:estmfkl}) não vale em geral, o que
não é problema, pois a exponencial de (\ref{eq:soltransf.1})
domina qualquer polinômio e com isso podemos mostrar que
$u(\cdot,t)\in\B_4\ \forall\ t>0$.
\eop

\obs Segue do Teorema~\ref{teo:contr} que existe $L_0>1$ tal que $R_L$ é
uma contração em $\{ f \in {\cal B}_4:\hat f(0)=0\}$ para qualquer
$L>L_0$.

\begin{cor}[Convergência]
\label{cor:conv}
Se $L>L_0$ e $f\in {\cal B}_4$ então
\be
   \label{eq:converg}
   R_L^n f \ton^{\B_4} A_f\phi_*,
\ee
onde $A_f=\int_\R f(y)dy$.
\end{cor}
\prova Defina $f_n = R_L^nf$ e $g_n = f_n - A_f \phi_*$, $n=0,
1,2, \dots$. Pela linearidade de $R_L$ e como $\phi_*=R_L\phi_*$,
temos que $g_{n+1} = f_{n+1} - A_f\phi_* = R_L f_n - R_L A_f \phi_*
= R_L g_n$. \'E suficiente ver que
$\widehat {g_n}(0) = 0 \ \forall n \in \N$, pois, pelo
Teorema~\ref{teo:contr} teremos $g_n\ton^{\B_4}0$. Mas $\widehat
{g_0}(0) = \int_\R g_0(y)dy = \int_\R (f-A_f\phi_*)(y)dy = \int_\R
f(y)dy - A_f\int_\R \phi_*(y)dy = 0$. Seja $n\ge 0$ e suponha que
$\widehat{g_n}(0)=0$, temos, por (\ref{eq:rgfour}), que
$\widehat{g_{n+1}}(0)=\widehat{R_Lg_n}(0)=0$, o que, por indução,
conclui esta prova. \eop

O corolário abaixo traz um resultado similar à
Proposição~\ref{prop:compasseqcalor}:
\begin{cor}
\label{cor:compasseqcalor2}
  Seja $u(x,t)$ a solução da equação do calor com dado
inicial $u(\cdot,0) = f\in \B_4$. Então vale
  \be
    \label{eq:resclim}
     \sqrt {t_n} \ u (\sqrt {t_n}\ x, t_n) \ton^{\B_4} A_f \phi_*(x)
  \ee
  para qualquer seqüência $t_n = L^{2n}$ com $L>L_0$,
  onde $A_f = \int_\R f(y)dy$.
  Em particular, temos convergência uniforme e convergência nas
  topologias de $L^1(\R)$ e $L^2(\R)$.
\end{cor}
\prova
Como a equação que estamos considerando tem
sua condição inicial no tempo $t=0$ e o grupo de renormalização
foi definido para soluções com condição inicial em
$t=1$, definimos a função $g$ como $g = u(\cdot,1)$.
Segue de (\ref{eq:soltransf.1}) que
$\hat g(0) = \hat f(0)$ e, portanto,
$\int_\R g(y)dy = \int_\R f(y)dy = A_f$.

Agora, veja que da \emph{Propriedade de Semigrupo} obtém-se
$R_L^n = R_{L^n}$ e,
pela definição de $R_L$~(\ref{eq:rg}), temos
$$
   R_L^ng(x)
   =
   R_{L^n}g(x)
   =
   {L^{n}} u (L^n x, L^{2n})
$$
Segue do Corolário~\ref{cor:conv} que
$$
  {L^{n}} u (L^n x, L^{2n}) \ton^{\B_4} A_f \phi_*(x),
$$
concluindo a demonstração.
\eop

\obs O limite pode ser mostrado para $t\to\oo$ em geral,
preenchendo-se os intervalos $[L^{2n},L^{2n+2}]$.
Entretanto, esse último passo foge ao objetivo desta dissertação, pois é
uma parte mais técnica e que ilustra pouco da utilidade do RG. Além
disso, estamos considerando sistemas cujas soluções são ``bem comportadas''
e é pouco razoável a possibilidade de tal limite valer para toda subseqüência
da forma acima mas não valer para $t\to\oo$ em geral. Por último,
no algoritmo numérico que vamos utilizar neste trabalho, somente
essas subseqüências serão consideradas. O leitor interessado pode consultar
\cite{braga04,braga03b,bric-kupa,bric-kupa-lin,moreira}.

\section{Renormalização das equações}
\label{sec:rggeral}

Na seção anterior definimos o grupo de renormalização motivados pelo fato
de a equação do calor ser invariante pela mudança de escalas
(\ref{eq:mudanca.1}). O que est\'a por tr\'as do fato de o m\'etodo
determinar o comportamento assint\'otico do pvi
\'e a Propriedade de
Semigrupo (para a equação do calor,
veja como a demonstração do Corolário~\ref{cor:compasseqcalor2}
utiliza o Teorema~\ref{teo:semi-grupo}).
Na verdade, a transformação do
RG aplicada \`a equação do calor \'e um caso particular da aplicação
desse m\'etodo. No caso mais geral (veja \cite{bric-kupa,bric-kupa-lin}),
a equação dada não
\'e necessariamente invariante pela ação das operaç\~oes
(\ref{eq:mudanca}). Por isso temos que ampliar o conceito do RG, fazendo o
mesmo agir tamb\'em sobre a equação diferencial dada. O que se faz \'e
definir a solução reescalonada $v$ como em (\ref{eq:mudanca.1}) e
considerar a equação satisfeita por $v$. Assim, o esquema para o caso
geral \'e:
\be
\ba{c}
  {\rm eq}_0 \\
  f_0
\ea
\mapsto
\ba{c}
  {\rm eq}_1 \\
  f_1
\ea
\mapsto
\ba{c}
  {\rm eq}_2 \\
  f_2
\ea
\mapsto \cdots
\mapsto
\ba{c}
  {\rm eq}_n \\
  f_n
\ea
\mapsto \cdots,
\ee
onde o $n$-\'esimo mapeamento evolui o $n$-\'esimo dado inicial $f_n$ pela
$n$-\'esima equação ${\rm eq}_n$, com $n=0,1, 2, \dots$. Portanto a
transformação do grupo de renormalização ter\'a, no caso geral,
duas componentes: uma atuando sobre o espaço dos dados iniciais e a
outra atuando sobre a equação diferencial em questão. Com isso mantemos a
propriedade de semigrupo, fato essencial para que possamos fazer a
evolução temporal do dado inicial usando o RG.

Seja $u(x,t),\ x\in\R,t\geq 1$, solução global\footnote
{Não estamos preocupados neste momento com a questão da existência
global de tal solução, mas simplesmente assumindo este fato.}
de um dos dois problemas
apresentados no Capítulo~\ref{chap:intro}. Ao invés de tomarmos o
limite~(\ref{eq:long-time-lim}) de uma só vez, o mesmo será fruto de um
processo iterativo. Para um $L>1$ fixo e expoentes $\alpha>0$ e $\beta>0$,
definimos\footnote{
A forma mais ``correta'' seria definir
$u_n(x,t) = L^{\alpha} u_{n-1}(L^{\beta} x,Lt)$.
A substituição de $L$ por $L^2$ neste capítulo
pode ser vista como uma questão simplesmente estética.
O RG numérico, que veremos no Capítulo~\ref{chap:rgnum},
utiliza $L$ como fator de escala ao invés de $L^2$.}
uma seqüência $(u_n)_{n\in\N}$ de funções, por $u_0 = u$ e, para
$n=1,2,3,\dots$
\be
  \label{eq:fixscale}
  u_n(x,t) = L^{2\alpha} u_{n-1}(L^{2\beta} x,L^2t).
\ee
Diremos que $(u_n)$ é uma \emph{seqüência de soluções reescalonadas}. A
partir da definição~(\ref{eq:fixscale}) e sabendo-se qual equação
diferencial é satisfeita por $u_{n-1}$, é um exercício simples de cálculo
obter a equação diferencial satisfeita por $u_n$. A estas equações
chamaremos \emph{seqüência de equações renormalizadas}.

Vamos definir $f_n$ por
\be
  \label{eq:deffn}
  f_n(x) = u_n(x,1).
\ee
Segue de (\ref{eq:deffn}) e (\ref{eq:fixscale}) que
\be
  \label{eq:fnunmu}
  f_n(x) = L^{2\alpha} u_{n-1}(L^{2\beta} x, L^2).
\ee
Observe que $f_{n+1}$ pode ser obtida evoluindo-se $f_{n}$ de $t=1$ até $t=L^2$,
de acordo com o
fluxo que determina $u_{n}(\cdot,t)$, e aplicando-se as mudanças de
escalas $x\mapsto L^{2\beta} x,\ u\mapsto L^{2\alpha} u$.
Observe também que uma conseqüência imediata de~(\ref{eq:fixscale}) é
\be
  \label{eq:fnu}
  f_n(x) = u_n(x,1) = L^{2n\alpha} u \left( L^{2n\beta} x, L^{2n}  \right).
\ee
Ressaltamos a importância da relação acima. Suponha que, para o
par de expoentes $\alpha$ e $\beta$ dados, sejamos capazes de
mostrar que $f_n\to\phi$, onde $\phi:\R\to\R$ não é identicamente nula.
Neste caso, teremos provado o limite~(\ref{eq:long-time-lim}) para a
subseqüência $t_n=L^{2n}$. Para provar o mesmo limite para qualquer seqüencia
$t\to\oo$ de forma rigorosa, ficaria
faltando então preencher os intervalos $[L^{2n},L^{2n+2}]$
-- veja a observação após
o Corolário~\ref{cor:compasseqcalor2}. Essa é, essencialmente, a forma
com que a técnica do RG analítico é aplicada.
A seguir, aplicaremos
as idéias explicadas neste parágrafo para o problema de difusão não-linear.
Os Problemas \ref{prob:periodico} e \ref{prob:deptempo} serão considerados
nas subseções seguintes.

Consideremos o problema de difusão não-linear
(Problema~\ref{prob:periodico} com $\mu=0$).
Vamos assumir que a perturbação é da forma $F = u^au_x^bu_{xx}^c$.
Se $u$ é solução global do problema, então a \emph{$n$-ésima equação
renormalizada} será aquela satisfeita por $u_n$ definida em
(\ref{eq:fixscale}), dada por
\be
  \label{eq:nldifrenorm}
    u_t = \chi_nu_{xx} + \lambda_n u^au_x^bu_{xx}^c,
\ee
onde
$$
  \chi_n = (L^{1-2\beta})^{2n} \quad \mbox{e} \quad \lambda_n =
  (L^{1-(b+2c)\beta-(a+b+c-1)\alpha})^{2n}\ \lambda.
$$

Se estamos interessados em caracterizar os problemas cujas soluções se
comportem como no caso da equação do calor, esperamos algum tipo de
convergência e por isso não é interessante que
o termo $\chi_n$ tenda a $0$ ou $\oo$ quando $n$ cresce. Para isso,
devemos ter $\beta=1/2$. Como a equação do calor é conservativa\footnote
{Dizer que a equação é \emph{conservativa} significa que suas soluções
satisfazem
$\int_\R u(x,t)dx = \int_\R f(y)dy\ \forall\ t>1$}, devemos ter $\alpha =
\beta$. Substituindo $\alpha = 1/2 = \beta$ temos que
$$
  \chi_n = 1 \quad \mbox{e} \quad \lambda_n =  (L^{-d_F})^n\ \lambda,
$$
de forma que a parte não-linear vai tender a se anular,
permanecer constante, ou tender a infinito, se $d_F = a+2b+
3c- 3$ for positivo, nulo ou negativo, respectivamente.

Ainda
que não constitua uma prova rigorosa,
o argumento
acima permite classificar, mediante um exercício simples de derivação,
uma gama de perturbações não-lineares da forma $F=u^au_x^bu_{xx}^c$
que estão na mesma classe de universalidade\footnote
{Veja a discussão sobre classes de universalidade na Introdução,
página \pageref{ind:univ}.}
da equação linear, caracterizada pelos expoentes críticos $\alpha$ e $\beta$.

Com essa abordagem, Bricmont
et~al.~\cite{bric-kupa-lin} provaram o teorema
abaixo\footnote{Um
material acessível e minuncioso em que estudam-se perturbações da forma
$F=u^a,\ a=4,5,6,\dots,$ é a tese de mestrado de
Moreira~\cite{moreira}.
O resultado de Bricmont et~al. é mais geral do que o Teorema~\ref{teo:nldif}
que enunciamos aqui.}.

\begin{teo}[Bricmont et al., 1994]
  \label{teo:nldif}
Considere o Problema~\ref{prob:periodico} com $\mu=0$
e condição inicial $u(\cdot,0)=f\in \B_4$.
Suponha que a perturbação $F$ possa ser escrita
na forma $F(u,u_x,u_{xx})= u^au_x^bu_{xx}^c$.
Nessas condições,
se $a+2b+3c-3>0$ e $\lambda$ é suficientemente pequeno,
então  existe $A\in\R$ tal que
$$
  \sqrt t\ u(\sqrt t\ x, t) \ton^{\B_4} A\phi_*(x)
  ,
$$
onde $u$ é a solução do problema de valor inicial
e $\phi_*$ é definido em (\ref{eq:gauss}).
\end{teo}

Em outras palavras, Bricmont et~al. mostraram que $u$ é assintoticamente
auto-similar da forma
$$
  u(x,t) \approx \frac A {\sqrt t}\ \phi_* \left(
  \frac x {\sqrt t}\right).
$$

\begin{defi}[A transformação do grupo de renormalização]
\label{def:rg}
Para $n=0,1,2,\dots$, definimos o operador $R_{L,n}$ como
$$
  R_{L,n}f = L u_{f}( L \cdot, L^2 ),
$$
onde $u_{f}(x,t)$ denota a solução da $n$-ésima equação renormalizada
(\ref{eq:nldifrenorm}) com condição inicial $u_{f}(x,1)=f(x)$.
\end{defi}

De acordo com a
Definição~\ref{def:rg} e de acordo com
as equaç\~oes (\ref{eq:fixscale}), (\ref{eq:deffn}) e (\ref{eq:fnunmu}),
temos que
$$
  f_{n+1}(x) = L u_n ( Lx, L^2 ) = R_{L,n}f_n.
$$

Usando a identidade acima e denotando por $R_L^n$ a composição
$R_{L,n-1}\circ R_{L,n-2}\circ \cdots \circ
R_{L,1}\circ R_{L,0}$, segue de (\ref{eq:fnu})
que
\be
\label{eq:semi-prop}
  f_n(x) = L^n u (L^nx, L^{2n}) =  R_{L^n,0}f(x) = R^n_Lf(x), \quad
  n= 1, 2, \dots,
\ee
ao que chamaremos de ``propriedade de semigrupo''.

\subsection{Difusão com coeficiente periódico}

Considere o Problema~\ref{prob:periodico} com $F=u^au_x^bu_{xx}^c$.
Seguindo o mesmo procedimento feito para a equação de difusão
não-linear acima, temos que a \emph{$n$-ésima equação renormalizada}
será
\be
  \label{eq:rescaledeqb}
  u_t = \chi_n \,
  [1 + \mu g(\omega_nx) ] u_{xx}
  + \lambda_n \, u^au_x^b  u_{xx}^c,
\ee
onde
$$
\chi_n = (L^{1 - 2 \beta})^{2n}, \quad \omega_n = (L^\beta)^{2n}, \quad
\lambda_n = \lambda (L^{1 - (b + 2c) \beta -(a+b+c-1)\alpha})^{2n}.
$$

Observe que estamos no mesmo caso da difusão não-linear exceto pela
adição do termo $\mu g(x)$. No caso de $\mu = 0 \ne \lambda$,
o Teorema~\ref{teo:nldif} pode ser aplicado para garantir que
\be
  \label{eq:compmuz}
  u(x,t) \approx \frac A {\sqrt t}\ \phi_* \left( \frac x {\sqrt t} \right).
\ee
Gostaríamos de entender como o termo $\mu g(x)$ afeta a classe
de universalidade do comportamento assintótico auto-similar dado
por (\ref{eq:compmuz}). Para isso, é natural escolher novamente
$\alpha=\beta = 1/2$. Nesse caso temos
$$
  \chi_n = 1, \quad \omega_n = L^{n}, \quad
  \lambda_n = \lambda (L^{-d_F})^{n}.
$$

Consideremos,  inicialmente, o Problema 1 supondo que
$\mu \not= 0 = \lambda$  (equação linear com coeficiente
de difusão peri\'odico).
Repare que, a cada renormalização, temos uma equação similar, porém
com o termo $\mu g(x)$ oscilando cada vez mais rapidamente.
A teoria de homogeneização,
que estudaremos no Capítulo~\ref{chap:homog}, nos diz que,
para grandes escalas,
as soluções desse problema -- dito \emph{não-homogêneo}, pois
reflete as propriedades de um meio não-homogêneo -- se aproximam da solução
de uma equação com coeficiente cons\-tan\-te, chamada \emph{equação
homogeneizada}. Nesse caso, conforme a discussão feita
no Capítulo~\ref{chap:intro} (e que será vista com mais detalhes no
Capítulo~\ref{chap:homog}), este coeficiente constante é dado pela média
harmônica do coeficiente $1+\mu g(x)$.
Tome $\epsilon_n^{-1} = L^{n}$. Observando que $\epsilon_n \ton 0$,
vemos que o estudo do comportamento assintótico da solução desse problema
é equivalente ao estudo, quando $\epsilon \to 0$, da solução da equação
$$
  u_t = [1+\mu g(\epsilon^{-1}x)]u_{xx}.
$$

De acordo com  a teoria de homogeneização, a solução da equação acima
se aproxima, quando $\epsilon\to 0$,
de um
múltiplo da solução de
\be
  \label{eq:line-homo-eqb}
  u_t = \sigma u_{xx},
\ee
onde $\sigma = \langle(1+\mu g)^{-1}\rangle^{-1}$ é a média harmônica do
coeficiente de difusão.

Consideremos,  agora, o caso geral para o Problema 1, isto \'e,
suponha agora  que
$\lambda \not= 0 \not=\mu $  (equação não-linear com coeficiente
de difusão peri\'odico).
A equação renormalizada é
dada por~(\ref{eq:rescaledeqb}).
Quando $n\to\oo$,  o termo não-linear
$\lambda_n F$ é supostamente levado a zero se $d_F>0$, enquanto que o
coeficiente periódico $1+\mu g(\omega_n x)$ oscila cada vez mais
rapidamente, produzindo um resultado
cada vez mais parecido com o de um coeficiente constante
$\sigma = \langle(1+\mu g)^{-1}\rangle^{-1}$.

Somos então levados a supor que
as soluções do Problema 1, dentro de uma classe de condições iniciais
adequada e com $F$ satisfazendo $d_F>0$, se comportam assintoticamente
como as soluções de (\ref{eq:line-homo-eqb}). Fazendo a mudança de variáveis
$x \mapsto x / \sqrt \sigma $, obtemos a equação do calor. Concluímos
então que as soluç\~oes do Problema 1 se comportam da mesma forma que
as soluç\~oes da equação de difusão não-linear, exceto pelo fator
$\sqrt\sigma$ no argumento da função de perfil $\phi_*$.
Com esses argumentos, Braga et~al.~\cite{braga03a} estabeleceram e verificaram
numericamente a conjectura abaixo.

\begin{conj}[Braga et al., 2003]
\label{conj:oscilante}
Considere o Problema~\ref{prob:periodico}. Suponha que a perturbação $F$
possa ser escrita na forma $F(u,u_x,u_{xx})= u^au_x^bu_{xx}^c$.
Se a condição inicial é dada por $u(\cdot,0)=f\in
C_0^\oo(\R)$, $a+2b+3c-3>0$ e $\lambda$ é suficientemente pequeno, então
existe uma única solução $u(x,t),\ x\in\R,t\geq 0$. Além disso, existe
$A\in\R$, tal que
$$
  \sqrt t u(\sqrt t x, t) \to A\phi_\sigma(x)
  \quad \mbox{quando} \quad t \to \oo,
$$
onde
$\phi_\sigma$ é definido
em~(\ref{eq:phi_sigma}), o limite valendo em $L^1(\R)$ e
$L^\oo(\R)$.
\end{conj}

Em outras palavras, temos que $u$ é assintoticamente
auto-similar, da forma
$$
  u(x,t) \approx \frac A {\sqrt t}\ \phi_\sigma \left(
  \frac x {\sqrt t}\right),
$$
com expoentes críticos $\alpha=\beta=\frac 12$ e
função de perfil $\phi_\sigma$.

\subsubsection{Um caso de perturbação relevante}

Para um caso particular de perturbações não-lineares
relevantes, é possível ainda fazer o estudo do
comportamento assintótico das soluções a partir da renormalização
da equação.

Considere o Problema~\ref{prob:periodico} com $\lambda F = -u^a,\
a\in(1,+\oo)$.
O mesmo argumento que justifica a Conjectura~\ref{conj:oscilante}
permite classificar as perturbações da forma acima como
irrelevantes se $a>3$. Pensemos agora no caso $a<3$ e
fixemos nossa atenção na parte não-linear do problema.
A renormalização com $\alpha=1/2$ faz esse termo crescer.
Podemos pensar que isso se deve ao fato de que
o termo não-linear produz um decaimento mais rápido do que
$t^{-1/2}$ -- as soluções de $u_t = -u^a$
são da forma $u = A\, \phi(x)\, t^{-1/(a-1)}$ com qualquer
$\phi$. Por outro lado, a parte linear não pode ser descartada
pois ela é responsável pela difusão. Nesse
caso, a renormalização correta é dada pelos expoentes
$\beta=1/2,\ \alpha=1/(a-1)$. Essa escolha de $\beta$ garante
que a parte linear seja invariante pela renormalização.
A escolha de $\alpha$ garante o mesmo para  a parte não-linear. Essa nova
escolha de expoentes determina um novo operador $R_L$ e uma
nova renormalização da equação. Em relação a essa nova renormalização,
tanto $u_{xx}$ quanto $u^a$ são termos \emph{marginais}.

Para o caso homogêneo ($\mu=0$), pode-se
seguir o mesmo raciocínio que
nos levou à equação (\ref{eq:ord}) para obter a equação a ser satisfeita
pelo ponto fixo do novo RG.
Essa equação possui uma solução $\phi_a$
positiva, de quadrado integrável,
que é única a menos de uma constante multiplicativa  -
veja~\cite{bric-kupa-lin} e referências lá citadas.
Bricmont et~al.~\cite{bric-kupa-lin} provaram que
se a condição inicial $u(\cdot,0)=f$ satisfaz
$f\in\B_4$ e $f(x)\geq 0\ \forall\ x\in\R$, então
$$
  u(x,t) \approx \frac A {t^{\frac 1 {a-1}}}\ \phi_a
  \left(  \frac x {\sqrt t} \right).
$$

Para o caso em que $\mu\ne0$, Braga et~al.~\cite{braga03a}
verificaram numericamente a conjectura abaixo.

\begin{conj}[Braga et al., 2003]
\label{conj:periodicorelevante}
Considere o Problema~\ref{prob:periodico} com
condição inicial $u(\cdot,0)=f\in C_0^\oo(\R),\,f\geq0$.
Suponha que  perturbação $F$ possa ser escrita como $F(u)=u^a$
para $u\geq0$ e que $\lambda<0$. Se $a\in(1,3)$ e
$\lambda$ é suficientemente pequeno, então existe
uma única solução $u(x,t),\ x\in\R, t>1$. Além disso, 
\begin{equation}
\label{eq:long-time-beha-rele}
  t^{\frac 1 {a-1}}u(\sqrt{t} x, t)\rightarrow
  f_a\left(\frac{x}{\sqrt{\sigma}}\right) \,\,\,\,\,\,
  \mbox{quando} \,\,\,\,\,\, t\rightarrow\infty,
\end{equation}
onde $\sigma = \langle(1+\mu g(x))^{-1}\rangle^{-1}$ e
$f_a$ é uma função tal que
$$
 u(x,t) \equiv \frac{1}{t^{\frac 1 {a-1}}}f_a\left(\frac{x}{t^{1/2}}\right)
$$
é solução do problema com $\mu=0$.
\end{conj}

\subsection{Difusão com coeficiente dependente do tempo}
\label{sec:anadeptempo}

Considere o Problema~\ref{prob:periodico} com $F=u^au_x^bu_{xx}^c$.
Reescrevemos a equação como:
$$
  u_t = [t^p + o(t^p)]\, u_{xx}
+ \lambda_n \, u^a u_x^b u_{xx}^c.
$$
Neste caso, a $n$-ésima equação renormalizada é dada por
\be
  \label{eq:renormnltp}
  u_t = \chi_n \left[ t^p + \frac{o(L^{np}t^p)}{L^{np}} \right]\, u_{xx}
  + \lambda_n \, u^a u_x^b u_{xx}^c,
\ee
onde
$$
  \chi_n = (L^{p+1-2\beta})^{2n}, \quad \lambda_n =
  \lambda (L^{1 - (b + 2c)\beta -
  (a+b+c-1)\alpha})^{2n}.
$$

Fazendo $\beta=\alpha=\frac{p+1}2$, temos
$$
  \chi_n = 1, \quad \lambda_n =
  \lambda (L ^{-(p+1) \eta_F})^n,
$$
onde
\be
  \label{eq:eta2}
  \eta_F = a+2b+3c - \frac{p+3}{p+1}.
\ee

Observe que $\lim_{s\to\oo}o(s)/s = 0$ e, portanto,
$$
  \sup_{t\in[1,L^2]} \frac{o(L^{np}t^p)}{L^{np}}
  \ton 0,
$$
de forma que, para grandes valores de $n$, este termo se torna
desprezível
na equação (\ref{eq:renormnltp}).
Por outro lado, se $\eta_F>0$, temos que $\lambda_n$
decresce rapidamente
a 0 quanto $n\to\oo$. Assim, somos levados a supor que o regime assintótico
das soluções do Problema~\ref{prob:deptempo} com $\eta_F>0$, e cuja
condição inicial $u(\cdot,t)=f$ é suave e decai suficientemente rápido,
está na mesma classe de universalidade da solução auto-similar de
\begin{equation}
\label{eq:difu.tp}
u_t = t^p\, u_{xx}.
\end{equation}
De fato,
$$
  u(x,t) = \frac A {t^\alpha}\ \phi_p \left(\frac x {t^\beta}\right),
$$
onde $\alpha = \beta= \frac {p+1}2$ e $\phi_p$ é dada por
(\ref{eq:phi_p}) é uma solução auto-similar de (\ref{eq:difu.tp}).

Para o caso linear $\mu=0$, Braga et.~al~\cite{braga04} obtiveram
o seguinte resultado.

\begin{teo}[Braga et al., 2004]
\label{teo:compu-0}
Considere o Problema~\ref{prob:deptempo} linear, isto
é, com $\lambda=0$. Se a condição inicial $u(\cdot,0)=f\in\B_4$,
então a solução $u(x,t)$ satisfaz
$$
  \label{eq:compu-0}
  \sqrt{t^{p+1}} u(\sqrt{t^{p+1}} \cdot,t) \tot^{\B_4}
  A_f\, \phi_p(\cdot),
$$
onde
$\phi_p(x)$ é dada por (\ref{eq:phi_p})
e $A_f = \int_\R f(y)dy$.
\end{teo}

Em outras palavras, $u$ é assintoticamente auto-similar com
expoentes $\alpha=\beta=\frac {p+1}2$, e o prefator
$A_f$ depende apenas da integral do dado inicial $f$.

Para o caso não-linear, Braga et~al. propuseram a seguinte conjectura.

\begin{conj}[Braga et al., 2004]
\label{conj:compassgeral}
Considere o Problema~\ref{prob:deptempo} com $F=u^au_x^bu_{xx}^c$ e
defina $\eta_F$ por (\ref{eq:eta2}). Então, se $\eta_F>0$ e
$\lambda$ é suficientemente pequeno, continua
valendo a conclusão do Teorema~\ref{teo:compu-0}, substituindo-se
$A_f$ por uma constante   $A$, que não depende apenas de $f$.
\end{conj}

\subsubsection{Transição de fase}

Considere a equação
\be
  \label{eq:tpua}
  u_t = t^p u_{xx} - u^a,
\ee
com espaço de parâmetros $(p,a)$.

Observe que a condição para que o termo $-u^a$ seja marginal,
isto é, para que $\eta_F$ dada por (\ref{eq:eta2}) se anule,
se reduz a
\be
  \label{eq:ap}
  a = a_c(p) = \frac {p+3}{p+1}.
\ee

Vamos analisar como o par de parâmetros $(p,a)$ afeta o expoente crítico
$\alpha$.

Para $a>a_c(p)$ temos que o termo $-u^a$ é supostamente irrelevante
pois tende a se anular sob a renormalização dada pelos expoentes $\alpha=
\beta = (p+1)/2$. Nesse caso, o comportamento assintótico das soluções
está na mesma classe de universalidade da equação linear e dizemos que
estamos no \emph{regime linear}.

Para $a<a_c(p)$, temos que o termo $-u^a$ tende a crescer
com as sucessivas renormalizações se os expoentes forem escolhidos como
no parágrafo acima. Nesse caso, a escolha correta dos expoentes
é $\alpha = 1/(a-1)$ e $\beta = (p+1)/2$. O espalhamento com
$t^\beta$ é causado pelo termo $t^pu_{xx}$ e o decaimento com
$t^\alpha$ é determinado por $-u^a$. O termo não-linear, que é considerado
relevante, passa a ser marginal em relação a esses novos expoentes críticos.
Nesse caso, dizemos que estamos no \emph{regime não-linear}.

O caso $a=a_c(p)$ é o mais delicado. Ambos os termos linear e não-linear
produzem o mesmo expoente crítico de decaimento $\alpha = \frac{p+1}2 = \frac 1{a-1}$.
Nesse caso, dizemos que estamos no \emph{regime crítico}.
No regime crítico aparece uma correção logarítmica ao decaimento
assintótico das soluções.

Podemos pensar que o termo determinante do expoente $\alpha$ será aquele
que, sozinho, causaria o maior decaimento. Nesse caso, temos
\be
  \label{eq:alphapa.2}
  \alpha(p,a) = \max\left\{\frac {p+1}2,\frac 1{a-1}\right\}.
\ee
Observe que a superfície dada pelo gráfico de $\alpha(p,a)$
tem suas singularidades justamente
sobre a curva $a = a_c(p)$. A Figura~\ref{fig:pht} mostra o gráfico
de $a_c(p)$ separando as regiões dos regimes linear e não-linear.
A Figura~\ref{fig:pht2} mostra o gráfico de $\alpha(p,a)$ com
singularidades em $a = a_c(p)$.

\begin{conj}[Braga et al., 2004]
\label{conj:phasetrans}
Considere a equação (\ref{eq:tpua})
com $p\in[0,\oo)$ e $a\in(1,\oo)$, com condição inicial $u(\cdot,0)=f$
satisfazendo
$f\in C_0^\oo(\R)$ e $f(x)\geq0\ \forall\ x\in\R$. Então existe uma única solução
$u(x,t),\,x\in\R,t\geq0$.

Se $a>a_c(p)$, então existe $A\geq0$, que depende de $p$, $a$ e $f$, tal que
$$
  \sqrt{t^{p+1}} u(\sqrt{t^{p+1}} \cdot,t) \tot^{\B_4}
  A_f\, \phi_p(\cdot),
$$
onde $\phi_p$ é dado por (\ref{eq:phi_p}).

Se $a = a_c(p)$, então existe $A\geq0$, que depende apenas de $p$, tal que
$$
  \sqrt{(t\ln t)^{p+1}} u(\sqrt{t^{p+1}} \cdot,t) \tot^{\B_4}
  A_f\, \phi_p(\cdot).
$$

Se $a<a_c(p)$, então existe $A>0$, que depende apenas de $a$ e $p$,
e $\phi_{p,a}\in L^1(\R)$ estritamente positiva tal que
$$
  t^{\frac{1}{a-1}} u(\sqrt{t^{p+1}} \cdot,t) \tot^{\B_4}
  A_f\, \phi_{p,a}(\cdot).
$$
\end{conj}

No Capítulo~\ref{chap:resul} vamos fazer um estudo em torno do ponto
$(p,a)=(1/2,7/3)$, que está sobre a curva crítica (\ref{eq:ap}).
Verificamos, entre outras coisas, que este é de fato um ponto de transição
de acordo com a Conjectura~\ref{conj:phasetrans} acima. Em particular,
vamos variar $p$ e $a$ separadamente e calcular o expoente $\alpha$
numericamente, verificando a fórmula (\ref{eq:alphapa.2}).

%% file: rgnum.tex
\chapter{O Grupo de Renormalização Numérico}
\label{chap:rgnum}

Neste capítulo vamos descrever uma
versão numérica do grupo de renormalização. O algoritmo, baseado nas
idéias do RG analítico, é dotado de elegância e simplicidade que merecem
destaque. Ademais, o método fornece um quadro detalhado do regime
assintótico: função perfil, expoentes críticos, prefatores e parâmetros
efetivos, sendo aplicável não só a problemas com soluções auto-similares
(que decaem como $\approx t^{-\alpha}$), mas também ao caso crítico
(ou marginal), em que aparecem
correções  logarítmicas ao decaimento $t^{-\alpha}$
(decaimento como $\approx (t\ln t)^{-\alpha}$).

Por tudo isso, o RG numérico parece ser a ferramenta ideal para estudar
quais aspectos do comportamento assintótico são universais e quais
dependem de parâmetros da EDP ou da condição inicial.
O método também dá suporte a
trabalhos rigorosos, pois os resultados obtidos podem encorajar tentativas
de se provar teoremas ou, eventualmente, apontar para a necessidade de
mudanças nas hipóteses.

Pelo que sabemos, a primeira versão numérica do RG para EDP's foi
proposta por Chen e Goldenfeld~\cite{chen}. Uma versão modificada
foi introduzida por Braga, Furtado e Isaia, sendo bem mais rápida
e eficiente que a primeira. Uma descrição detalhada desta versão
pode ser encontrada na tese de doutorado de Isaia~\cite{isaia02}.
Infelizmente, parte do ganho computacional desta última não se
aplica ao Problema~\ref{prob:periodico}, como veremos mais
adiante.

O método aqui apresentado aplica-se a equações da forma
\be
  \label{eq:genivp}
  \left\{ \ba{l}
  u_t = G(u,u_x,u_{xx},x,t) \\
  u(x,1) = f(x)
  \ea \right..
\ee
Observamos que os problemas apresentados no capítulo introdutório e que
serão estudados numericamente nesta dissertação podem ser escritos
na forma acima.

A idéia de uma seq\"u\^encia de mudanças de escalas e de
renormalização de equações será ge\-ne\-ra\-li\-za\-da para
seq\"u\^encias em que os expoentes de escala $\alpha$ e $\beta$
podem ser diferentes em cada passo.
Essa generalização será vista na Seção~\ref{sec:prelim},
juntamente outras noções preliminares
de que vamos precisar para explicar o funcionamento do método.
Cada iteração do método consiste basicamente nos seguintes
passos: \emph{evoluímos} uma dada função por um curto intervalo
de tempo, segundo o fluxo dado pela equação~(\ref{eq:genivp});
\emph{mudamos as escalas} do domínio e do contra-domínio da função obtida;
\emph{renormalizamos} a equação~(\ref{eq:genivp}). Esse
algoritmo numérico será descrito na Seção~\ref{sec:algor}.
Outras decisões de implementação, como a discretização
da equação, a escolha de uma escala $L$ e a mudança de escalas
no domínio da função serão discutidas na Seção~\ref{sec:imple}.
A validação do método numérico será feita no Capítulo~\ref{chap:resul}.

\section{Preliminares}
\label{sec:prelim}

Seja $u_0$ uma função real de $(x,t) \in {\mathbb{R}} \times
\mathbb{R}_+$. Para um $L > 1$ fixo e seq\"u\^encias de expoentes
positivos $(\alpha_n)_{n\in\N}$ e $(\beta_n)_{n\in\N}$, defina a
seq\"u\^encia de funções $(u_n)_{n\in\N}$, indutivamente, pela
seguinte \emph{mudança de escalas}:
\be
  \label{eq:scalingstep}
  u_{n}(x,t) = L^{\alpha_{n}}u_{n-1} \left( L^{\beta_{n}}x,Lt \right).
\ee

Repare que a definição~(\ref{eq:scalingstep}) acima implica em
\begin{equation}
  \label{eq:rescaling}
  u_{n}(x,t) =
  L^{\alpha_1 + \cdots + \alpha_n}
  u_0
  \left (L^{\beta_{1} + \cdots + \beta_n}x,L^nt\right ).
\end{equation}

Daqui em diante vamos supor que a função original $u_0$ é solução
global da Eq.~(\ref{eq:genivp}) com $f=f_0$ e $G=G_0$, de forma
que $(u_n)$ seja uma \emph{seq\"u\^encia de soluções
reescalonadas}. Neste caso, $u_n$ satisfaz ao seguinte
\emph{problema de valor inicial renormalizado}: \be
  \label{eq:renivp}
  \left\{ \ba{l}
  u_t = G_n(u,u_x,u_{xx},x,t) \\
  u(x,1) = f_n(x)
  \ea \right.,
\ee
onde $G_n$ e $f_n$ são obtidos recursivamente,  como explicamos a
seguir.

Primeiramente, o operador $G_n$ pode ser obtido a partir de $G_{n-1}$
aplicando-se substituição de variáveis e a regra da cadeia
\be
\label{eq:renop}
  G_n(u,u_x,u_{xx},x,t) = L^{\alpha_n+1} G_{n-1}
  ( L^{-\alpha_n}u, L^{-\alpha_n-\beta_n}u_x,
    L^{-\alpha_n-2\beta_n}u_{xx}, L^{\beta_n}x, Lt ).
\ee

Já a condição inicial $f_n$ pode ser obtida a partir de $f_{n-1}$
integrando-se a equação para $u_{n-1}$ no intervalo de tempo $[1,L]$ e
mudando-se as escalas do espaço e da própria função, pois,
pela definição~(\ref{eq:scalingstep}), temos
$$
  f_n(x) = L^{\alpha_n} u_{n-1} (L^{\beta_n}x, L).
$$
Temos, dessa forma, uma \emph{seq\"u\^encia de equações
renormalizadas}.

Queremos agora obter uma relação análoga à Propriedade de
Semigrupo. Note que por (\ref{eq:rescaling}) é possível relacionar
$f_n$ com a solução original:
\be
  \label{eq:rescaled-eqc}
  f_n(x) =
  L^{\alpha_n+\cdots +\alpha_1}u_0\left
  (L^{\beta_n+\cdots+\beta_1}x,L^n\right).
\ee

Uma simples análise da equação (\ref{eq:rescaled-eqc})  fornece a
relação entre a din\^amica da seq\"u\^encia $(f_n)$ quando
$n\to\infty$ e o comportamento assintótico de $u_0$. Para isso,
reescreva (\ref{eq:rescaled-eqc}) como
\be
  \label{eq:prefatores2}
  u_0(x,L^n) =
  \frac{A_n}{L^{n\alpha_n}}f_n\left(B_n\frac{x}{L^{n\beta_n}}\right),
\ee
onde $A_n$ e $B_n$ são chamados de \emph{prefatores} e podem ser
calculados por
\be
  \label{eq:prefatores}
  A_n = L^{n\alpha_n - (\alpha_n+\cdots +\alpha_1)}, \quad B_n =
  L^{n\beta_n - (\beta_n+\cdots +\beta_1)}.
\ee

Assim, se os limites $A_n\rightarrow A$, $B_n\rightarrow B$,
$\alpha_n\rightarrow \alpha$, $\beta_n\rightarrow \beta$ e
$f_n\rightarrow \phi$ são obtidos quando $n\rightarrow\infty$,
podemos esperar que
\[
L^{n\alpha}u_0(L^{n\beta}x, L^n)=
\frac{A_n}{L^{n(\alpha_n-\alpha)}}f_n\left(B_n\frac{x}{L^{n(\beta_n-\beta)}}\right)
\ton A \phi(Bx).
\]

O limite acima, por sua vez, estabelece a auto-similaridade
do comportamento assintótico de $u_0$.

\section{O algoritmo}
\label{sec:algor}

A partir das definições preliminares, podemos sistematizar um método
construtivo para a obtenção das seq\"u\^encias
$(\alpha_n)$, $(\beta_n)$, $(A_n)$, $(B_n)$ e $(f_n)$,
em que cada termo das seq\"u\^encias é calculado recursivamente.

Sejam $L$, $f_0$ e $G_0$ tais que o problema de valor inicial~(\ref{eq:genivp})
possui única solução $u_0(x,t)$ para todo $x\in\R$ e $t\geq 1$.
O algoritmo iterativo consiste dos seguintes passos, para $n=1,2,3,\dots$,
\ben
\item
  \label{item:integration}
  Obtenha $u_{n-1}(x,L)$ evoluindo do dado inicial $f_{n-1}$ do tempo
  $t=1$ até $t=L$, de acordo com a equação~(\ref{eq:renivp}) para
  $u_{n-1}$.
\item
  \label{item:choice}
  Calcule os expoentes $\alpha_{n}$ e $\beta_{n}$, conforme
  discutiremos mais adiante.
\item
  \label{item:scaling}
  Defina $f_{n}(x) = L^{\alpha_{n}}u_{n-1}(L^{\beta_{n}}x,L)$.
\item
  \label{item:prefactors}
  Calcule os prefatores $A_{n}$ e $B_{n}$ por~(\ref{eq:prefatores}).
\item
  Calcule os parâmetros do operador $G_n$ a ser usado na próxima
  iteração, de acordo com a relação~(\ref{eq:renop}).
\een

É o algoritmo numérico acima que usaremos para verificar as
conjecturas já apresentadas ou até mesmo para se determinar os
expoentes críticos sem que os mesmos sejam conhecidos \textsl{a
priori} (como, por exemplo, na equação de Barenblatt).

O expoente $\alpha$ pode ser sempre computado no passo~\ref{item:choice}
pela fórmula $\alpha_{n} = -\ln_Lu_{n-1}(0,L)$ de forma que
$f_{n}(0)=1$. A razão para isso é a suposta din\^amica auto-similar
que esperamos encontrar: no regime assintótico, a solução $u_0$ em
$x=0$ se comporta como
\be
  \label{eq:decay}
  u_0(0,t)\sim c\ t^{-\alpha}.
\ee
Como, por construção, $u_n(0,1)=1$ para $n=1,2,3,\dots$, e dada a relação
entre $u_n$ no intervalo de tempo $t\in[1,L]$ e $u_0$ em $t\in[L^n,
L^{n+1}]$ na equação~(\ref{eq:rescaling}), temos que
$$
  { u_0(0,L^{n+1}) \over u_0(0,L^n) } = {1 \over  L^{\alpha_{n+1}} }
$$
e, por~(\ref{eq:decay}), esperamos que, quando $n\rightarrow \infty$,
$\alpha_n$ se aproxime de $\alpha$.

Enquanto o cálculo dos expoentes $(\alpha_n)$ pode ser sempre
feito da forma descrita acima, a forma de se determinar os
expoentes $(\beta_n)$ depende do problema. Normalmente envolve uma
relação de escala entre $\alpha_n$ e $\beta_n$ de forma que uma
certa parte do operador diferencial (escolhida \textsl{a priori})
permaneça invariante sob a renormalização da equação.

Para ilustrar esta idéia, vamos rever brevemente o estudo da equação
$$
  u_t = u_{xx} + \lambda u^a u_x^b u_{xx}^c.
$$
A equação renormalizada é dada por
$$
  u_t = L^{1-2\beta_n}
  u_{xx} +
  L^{[1-(b+2c)\beta_n+(1-a-b-c)\alpha_n]} \lambda
   u^a u_x^b u_{xx}^c.
$$

Se escolhermos $\beta_n=\frac 1 2$, a parte linear da equação é mantida
invariante e a din\^amica associada \`a equação
\be
  \label{eq:difu}
  u_t = u_{xx}
\ee
pode ser explorada.

Por outro lado, se $\beta_n$ é escolhido de forma a
satisfazer a relação
\be
  \label{eq:scal}
  1 - (b+2c)\beta_n + (1 -a -b - c)\alpha_n = 0,
\ee
então é o operador não-linear
\be
  \label{eq:nonl}
  u_t = \lambda u^au_x^bu_{xx}^c
\ee
que permanece invariante. Neste caso, é a din\^amica
da equação~(\ref{eq:nonl}) que pode ser explorada.

Na linguagem do RG, a escolha de $\beta_n = \frac 1 2$
concentra a atenção na din\^amica de (\ref{eq:difu})
e, dessa forma, é própria para a investigação
de perturbações irrelevantes e marginais.
A escolha (\ref{eq:scal}) leva o foco para
(\ref{eq:nonl}) e permite a investigação dos casos
de perturbações relevantes.

Por último, a mudança de escala da variável espacial, no
passo~\ref{item:scaling}, pode ser feita de duas formas. A
primeira foi utilizada por Chen e Goldenfeld~\cite{chen} e
consiste em modificar o tamanho da malha $\Delta x$ sem mudar os
sítios discretos $j=0,1,2,\dots$, de forma que depois da $n$-ésima
iteração temos $(\Delta x)' = L^{-\beta} \Delta x$ e os pontos da
nova malha estarão localizados em $x = j(L^{-\beta}\Delta x)$. A
outra forma consiste em mudar os sítios discretos mantendo o
tamanho da malha fixo, de forma que os sítios discretos, após a
mudança de escalas, estarão em $x = (L^{-\beta} j)\Delta x$. Neste
caso, os valores da solução $u$ nos pontos $x = j \Delta x$ da
malha fixa t\^em que ser interpolados em cada iteração.

A última técnica para se mudar a escala espacial referida acima
foi proposta por Braga, Furtado e Isaia (veja a tese de doutorado
de Isaia~\cite{isaia02} e referências lá citadas). Utilizando este
ponto de vista conseguimos um custo constante para cada iteração
do algoritmo. Neste caso, temos um enorme ganho computacional, uma
vez que, para se calcular $u(\cdot,L^n)$ a partir de $u(\cdot,1)$
precisamos fazer $n$ integrações no intervalo $[1,L]$ a custo
constante.

\section{Implementação}
\label{sec:imple}

Esta seção descreve algumas características da versão do algoritmo
descrito acima que foi implementada pelo autor desta dissertação.
O algoritmo numérico foi implementado
em linguagem C e a geração dos arquivos de entrada, bem como
dos gráficos e tabelas a partir dos arquivos de saída, foi feita
com o auxílio do software Maple\,\copyright.
Essa versão foi utilizada na obtenção dos resultados apresentados
nos artigos \cite{braga03a,braga04}
e no Capítulo~\ref{chap:resul} desta dissertação.
O leitor que esteja interessado apenas em conhecer
o algoritmo em si, ou que queira implementar sua própria versão,
pode dispensar a leitura desta seção.

\subsubsection{Equação e parâmetros}

O programa considera a equação
$$
  u_t = G(u,u_x,u_{xx},x,t) = \chi (t^p + \delta t^r) [ 1+\epsilon H(-(u^m)_{xx}) ]
  [ 1 + \mu \cos (\omega x) ] (u^m)_{xx}  +
  \lambda u^au_x^bu_{xx}^c,
$$
onde $H$ é a função de Heaviside.
Os parâmetros do programa são:
$L$, $\chi$, $p$, $\delta$, $r$, $\epsilon$, $\mu$, $\omega$, $m$,
$\lambda$, $a$, $b$ e $c$, além da condição inicial, do número
de iterações desejadas e da opção para mudança de escalas
na variável espacial (veja este tópico mais adiante).
Normalmente, escolhemos $\chi=m=\omega=1$ e $\epsilon=0$ para o
estudo dos Problemas~\ref{prob:periodico}~e~\ref{prob:deptempo}. O
parâmetro $m$ é usado no estudo de equações que não consideramos
nessa dissertação -- portanto, em todas as simulações aqui mencionadas,
assumimos a escolha de $m=1$ -- e o parâmetro $\epsilon$ é usado para
estudarmos a equação de Barenblatt na validação do método (veja a
Seção~\ref{sec:validacao} no próximo capítulo).

\subsubsection{Discretização da equação}

Para a discretização da EDP, podemos utilizar qualquer esquema
apropriado. Escolhemos
um esquema simples de diferenças finitas que combina o método
de Euler para a discretização no tempo e, para a discretização no
espaço, a fórmula padrão de tr\^es pontos para o operador
Laplaceano e diferenças centradas para derivadas espaciais de primeira
ordem. Dentro das limitações de estabilidade, o esquema resultante é de
segunda ordem.

Para discretizar o espaço-tempo, fazemos $x=j\Delta x$ e $t = n \Delta
t$, onde $\Delta x$ é o tamanho da malha e $\Delta t$ é o passo no
tempo. Integramos a versão discreta da EDP para obtermos os dados
$(u(j\Delta x,L))_j$ que aproximam a função $u(x,L)$ nos sítios
discretos $x = j \Delta x$. Denotaremos por $u_j^n$ a aproximação para
$u(j\Delta x, n \Delta t)$.

As derivadas serão aproximadas da seguinte forma:
\be
   u_{xx}(j\Delta x,n \Delta t) \approx { u_{j+1}^n -2
   u_{j}^n + u_{j-1}^n \over (\Delta x)^2 },
\ee
\be
  u_x(j\Delta x,n \Delta t) \approx { u_{j+1}^n -
   u_{j-1}^n \over 2(\Delta x) }.
\ee

Assim, temos o esquema de Euler explícito
$$
  u_j^{n+1} = u_j^n + \Delta t \,  G\left(
    u_j^n,  { u_{j+1}^n - u_{j-1}^n \over 2(\Delta x) },
    {u_{j+1}^n -2 u_{j}^n + u_{j-1}^n \over (\Delta x)^2},
    j\Delta x, n\Delta t
  \right).
$$

\subsubsection{Condição de estabilidade}

Reescrevemos a equação como
$$
  u_t = G(u,u_x,u_{xx},x,t) = K(x,t,u_{xx}) u_{xx}  +
  \lambda u^au_x^bu_{xx}^c.
$$
Uma condição suficiente de estabilidade será\footnote
{Veja o livro de Ames~\cite{ames92}, Seção~2-18, p.~104.}:
$$
  \der G{u_{xx}} \pm \frac {\Delta x}2 \der G{u_x} \geq 0
  \qquad \mbox{e} \qquad
  1+\der G u \,\Delta t - 2 \der G{u_{xx}}
  \left(\frac {\Delta t}{{\Delta x}^2}\right) \geq 0.
$$

Tome $K',K''\in\R$ satisfazendo
$0<K'\leq K(x,t,u_{xx}) \leq K'' <\oo$.

Supondo que $\lambda=0$, temos
$$
  \frac {\Delta t}{\Delta x^2} \leq \frac 1 {2 {K''}}
$$
como condição suficiente de estabilidade.

Supondo que a solução $u$ da equação satisfaça\footnote{
A limitação imposta sobre a solução $u$ e suas derivadas
é satisfeita por funções da forma $\phi(x)=\exp(-kx^2/4),\,k\leq 3$,
por exemplo.
} $|u|\leq 1$,
$|u_x|\leq 1$ e $|u_{xx}| \leq 1$, temos as condições abaixo
como condições suficientes para a estabilidade:
\begin{eqnarray*}
  \Delta x &\leq& \frac 2 {b|\lambda|}({K'}-c|\lambda|) \\
  \frac {\Delta t}{\Delta x^2} &\leq&
  \frac 1 {\Delta x^2 a |\lambda| + 2 {K''} + 2c|\lambda|}.
\end{eqnarray*}

\subsubsection{Condição de fronteira}

A vantagem do método explícito é que ele torna a velocidade de
propagação da informação finita.
Isso possibilita o estudo de problemas de valor inicial
com dado inicial de suporte compacto adicionando-se mais um sítio
discreto na variável espacial a cada novo passo de tempo.
Dessa forma, a solução não é afetada por condições de contorno artificiais.

\subsubsection{Mudança de escala na variável espacial}

As duas formas de mudar a escala no passo~\ref{item:scaling} do algoritmo
foram implementadas. A interpolação implementada foi a interpolação linear
canônica. Uma desvantagem da interpolação é que ela não conserva a
integral discreta da função interpolada. No estudo das equações lineares
conservativas, em que o prefator depende apenas da integral do dado
inicial, utilizamos a mudança de escala sem interpolação para tornar o
resultado mais fino. No Problema~\ref{prob:periodico} com $\mu\ne 0$, as
oscilações do coeficiente de difusão têm seu período diminuído na mesma
proporção que os intervalos da malha discreta no passo de mudança de
escalas; neste caso, para enxergar o efeito da oscilação deste coeficiente
com certa fidelidade, o recurso da interpolação também foi deixado de
lado. Nas demais simulações utilizamos a interpolação.

Ressaltamos que \'e poss\'\i vel fazer uma interpola\c c\~ao que conserve
a integral discreta da fun\c c\~ao interpolada.  Para isso deve-se adotar
outro ponto de vista no momento da interpola\c c\~ao, que consiste em
associar os valores \`as c\'elulas e n\~ao aos pontos. Para o leitor
interessado em se aprofundar nesta dire\c c\~ao, indicamos o livro de
LeVeque~\cite{leveque92}.

\subsubsection{Fator de escala $L$}

A escolha do valor de $L$ é puramente empírica. Em princípio qualquer
fator de escala $L>1$ poderia ser usado. Para um valor de $L$ muito grande,
um elevado número de passos no tempo seria necessário para cada iteração
do RG, e a malha espacial seria rapidamente contraída em poucas iterações
devido à mudança de escala na variável espacial. Para um valor de
$L$ muito pequeno, seria necessário um número muito maior de iterações
para se obter a acurácia desejada.
Escolhemos um fator de escala $L=1,4$ nas simulações
em que utilizamos a interpolação no passo~\ref{item:scaling};
escolhemos $L$ entre $1,02$ e $1,021$ nos demais casos.

%% file: resul.tex
\chapter{Resultados Numéricos}
\label{chap:resul}

O objetivo deste capítulo é estudar os
Problemas~\ref{prob:periodico}~e~\ref{prob:deptempo} utilizando a
implementação do RG
numérico descrita no capítulo anterior.
Primeiramente, procuraremos validar o método utilizando o mesmo
no estudo de problemas de valor inicial em que soluções analíticas
são conhecidas ou para os quais existem teoremas sobre a forma
assintótica das soluções. Em seguida fazemos um estudo numérico
da convergência do método, através de um refinamento de malha.
Finalmente, usamos o algoritmo do RG para estudar os Problemas
\ref{prob:periodico} e \ref{prob:deptempo}
do Capítulo~\ref{chap:intro}. Cumpre notar que fazemos
um estudo padronizado da equação do calor, em que verificamos:
\begin{itemize}
\item a convergência do expoente $\alpha_n$ para o valor teórico
$\alpha=1/2$ com diferentes
condições iniciais, de forma a confirmar a universalidade
desse expoente crítico com respeito ao dado inicial;
\item a convergência
do prefator $A_n$ para algum $A>0$ em todas a simulações feitas;
\item   que o prefator $A$ é de fato dado por
$(4\pi)^{-1/2}\int_\R f(y)dy$;
\item a convergência das funções de perfil nas normas $L^1$ e $L^\oo$;
\item   que a
função de perfil $\phi$, obtida como $\lim_{n\to\oo}f_n$,
é de fato a gaussiana $\phi_*$.
\end{itemize}
Nos estudos subseqüentes, alguns dos gráficos referentes aos ítens
acima serão omitidos por não trazerem grande novidade.
As convergências de $(B_n)$ e $(\beta_n)$
são obtidas trivialmente uma vez que o procedimento escolhe
sempre o mesmo valor de $\beta_n$ para manter o termo linear
invariante sob a renormalização da equação.

A Seção~\ref{sec:validacao} é dedicada à validação do método.
Começamos pelo estudo padronizado da equação do calor, comprovando
a universalidade do expoente $\alpha$ e da função de perfil
$\phi$, assim como a expressão para o prefator em termos da integral
da condição inicial. Daí em diante, alguns gráficos
análogos, que foram gerados em outras simulações numéricas,
são omitidos
por se tratar de repetições monótonas que acrescentam poucas informações,
como a convergência da diferença relativa entre duas funções de perfil
para zero ou a verificação de que a função de perfil é dada pela gaussiana.
Estudamos então a equação de difusão não-linear e verificamos que
as perturbações não-lineares classificadas como \emph{irrelevantes}
de fato não alteram o comportamento assintótico exceto pelo prefator $A$.
A equação de Barenblatt é um exemplo bem ilustrativo: o expoente
$\alpha$ não pode ser determinado de forma exata analiticamente nem
formalmente, restando o recurso numérico como alternativa para se
calcular este expoente. A validação neste caso é feita comparando-se
o expoente calculado numericamente com aquele advindo de
uma aproximação de primeira ordem conhecida
na literatura. Por último, escolhemos uma das condições iniciais
e repetimos seu estudo com uma malha duas vezes mais fina. A semelhança
visual entre os gráficos obtidos, exceto por uma ligeira suavização
e maior precisão para a malha mais fina são os indicativos de convergência
do método com relação ao refinamento da malha.

Na Seção~\ref{sec:periodico}, estudamos o
problema de difusão com coeficiente periódico
(Problema~\ref{prob:periodico} da Introdução). Os resultados apresentados
nessa seção são essencialmente aqueles feitos no artigo \cite{braga03a}
e constituem-se principalmente da verificação da
Conjectura~\ref{conj:oscilante} do Capítulo~\ref{chap:rgana}.
O fato que traz mais novidade é a mudança da função de perfil
devido à homogeneização do coeficiente de difusão.

Na Seção~\ref{sec:deptempo} é feito o estudo dos teoremas e das
conjecturas da Seção~\ref{sec:anadeptempo} do Capítulo~\ref{chap:rgana},
que tratam do problema de difusão com coeficiente dependente do tempo
(Problema~\ref{prob:deptempo} da Introdução).
Esse estudo pode ser encontrado no artigo em elaboração
\cite{braga04}. As parte mais interessantes, na opinião do autor,
são a ampliação da classe de perturbações não-lineares
\emph{irrelevantes} para
$p>0$, que passa a incluir todas aquelas que eram marginais
nas outras equações, como $u^3$, $uu_x$ e $u_{xx}$, e a transição
de fase entre os regimes assintóticos linear e não-linear para
uma forma particular da equação com expoente
não-inteiro e solução positiva.

Como dito acima, os gráficos que trazem informação repetitiva
são omitidos nas Seções~\ref{sec:periodico}~e~\ref{sec:deptempo}.

Em todo este capítulo, o símbolo $\phi$ vai denotar a função $f_n$ obtida no
último passo da iteração do algoritmo, após terem sido obtidas as
convergências desejadas.

Chamamos a atenção para uma diferença
na definição das funções de perfil. A versão numérica que
implementamos (descrita no Capítulo~\ref{chap:rgnum}), calcula $\alpha_n$
de maneira a forçar que $f_n(0)=1$ para $n\geq 1$. Por este motivo,
será conveniente redefinir
$$
  \phi_*(x) = e^{-x^2/4}, \quad
  \phi_p(x) = e^{-(p+1)x^2/4} \quad \mbox{e} \quad
  \phi_\sigma(x) = e^{-x^2/4\sigma}.
$$
Isso implica na necessidade das constantes $\sqrt{1/ 4\pi}$,
$\sqrt{(p+1)/4\pi}$ e $\sqrt{1/4\pi\sigma}$, respectivamente, multiplicando
o prefator $A$ dos teoremas e conjecturas apresentados no
Capítulo~\ref{chap:rgana} e mencionados na Introdução. Em particular,
teremos as seguintes relações, que serão usadas mais adiante:
\begin{eqnarray}
  \label{eq:logphi*} -\ln \phi_* &=& x^2/4 \\
  \label{eq:logphip} -\ln \phi_p &=& (p+1)\cdot(x^2/4) \\
  \label{eq:logphisigma} -\ln \phi_\sigma &=& \frac 1 \sigma \cdot(x^2/4).
\end{eqnarray}

\section{Validação do método}
\label{sec:validacao}

\subsection{Equação do calor}

Vimos no Capítulo~\ref{chap:rgana} que o comportamento assintótico
das soluções da equação do calor com condição inicial
$u(\cdot,0)=f\in\B_4$ é dado por
\be
  \label{eq:ncomasscalor}
  u(x,t) \approx \frac A{t^\alpha}\ \phi\left(\frac x {
  t^\beta}\right),
\ee
com $\alpha=\beta=1/2$, $\phi=\phi_*$ e $A=(4\pi)^{-1/2}\int_\R f(y)dy$.
Segundo os resultados
teóricos já apresentados, $\alpha$, $\beta$ e $\phi$ são universais
com respeito à condição inicial, desde que a mesma esteja, por exemplo,
em $\B_4$. Por outro lado, o prefator
$A$ depende apenas do que chamamos de \emph{massa} do dado inicial
$f$, dada por $\int_\R f(y)dy$.

As três condições iniciais utilizadas, $f_1$, $f_2$ e $f_3$, são
exibidas nas Figuras~\ref{fig:f1},~\ref{fig:f2}~e~\ref{fig:f3},
respectivamente.
\insertfigure{fig:f1}{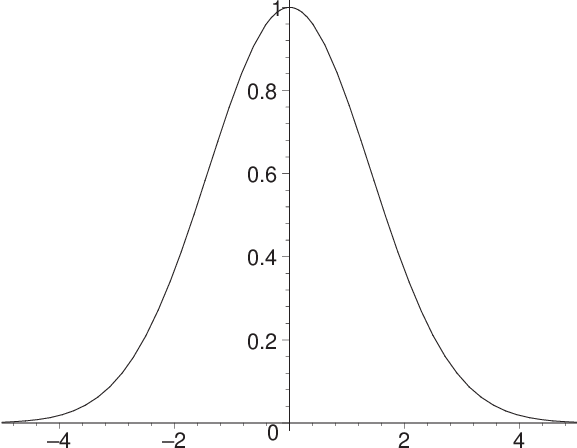}
{Gráfico da condição inicial $f_1$.}
\insertfigure{fig:f2}{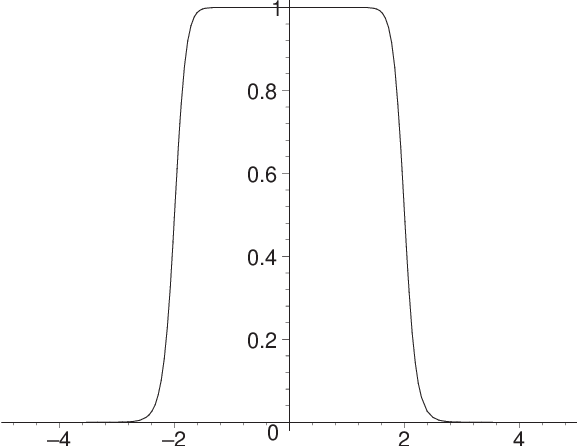}
{Gráfico da condição inicial $f_2$.}
\insertfigure{fig:f3}{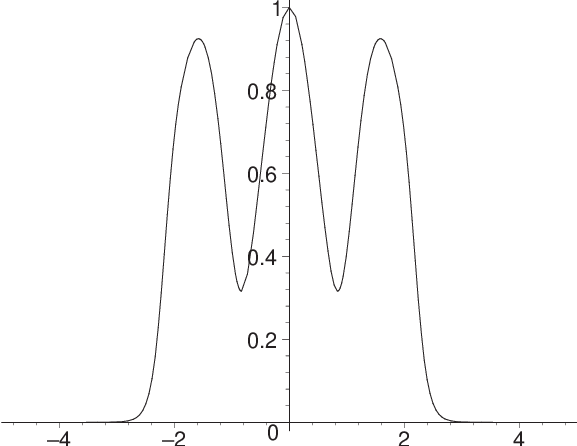}
{Gráfico da condição inicial $f_3$.}

A Figura~\ref{fig:alphaheat} mostra a converg\^encia
do expoente $\alpha_n$, calculado na $n$-ésima iteração do
algoritmo, para o valor teórico $\alpha = 1/2$, com
diferentes condições iniciais. Isso ilustra a universalidade de $\alpha$
com respeito a condições iniciais suficientemente localizadas (por
exemplo, de suporte compacto).
\insertfigure{fig:alphaheat}{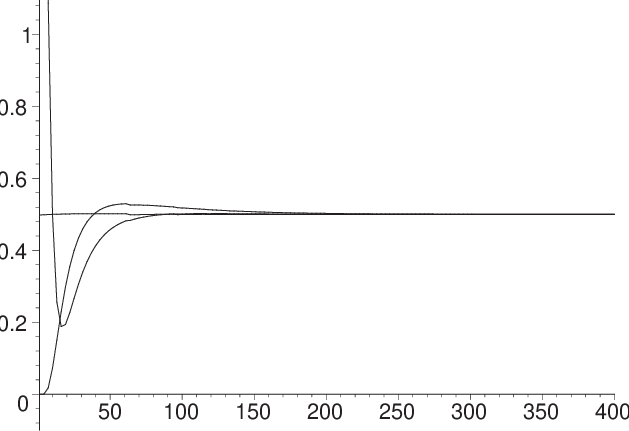}
{Curvas $\alpha_n \times n$ para a equação do calor com
diferentes condições inicias.}

A Figura~\ref{fig:aheat} mostra a converg\^encia dos prefatores $A_n$
para os valores teóricos. 
\insertfigure{fig:aheat}{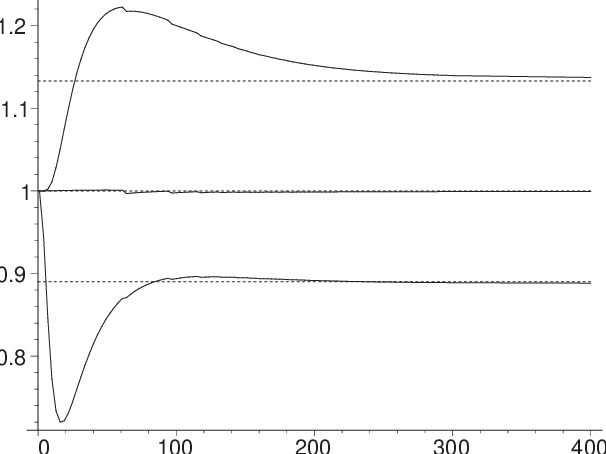}
{Curvas $A_n \times n$ para a equação do calor com
diferentes condições iniciais. As linhas pontilhadas correspondem
aos valores teóricos. Estes valores, dados por
$(4\pi)^{-1/2}\int_\R f(y)dy$, são
1,000, 1,133 e 0,8901.}

A Figura~\ref{fig:profileheat} mostra uma curva parametrizada
cujas componentes são $(x^2/4,-\ln\phi(x))$. O fato de todo ponto
desta curva estar em uma reta de inclinação igual a $1$
significa que $\phi = \phi_*$.
\insertfigure{fig:profileheat}{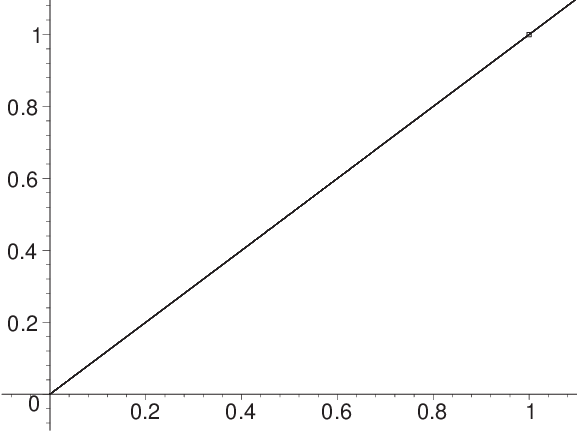}
{Curvas $-\ln\phi \times x^2/4$ para a equação do calor com
diferentes condições iniciais (as curvas se sobrepõem).}

Finalmente, na Figura~\ref{fig:normheat} plotamos a diferença
relativa entre duas funções de perfil consecutivas, $||f_n -
f_{n-1}|| / ||f_n||$, como uma função do número $n$ de
iterações do RG. As normas que utilizamos em todas as simulações são $L^1$
e $L^\infty$. O gráfico exibido na Figura~\ref{fig:normheat} foi
gerado com a norma $L^1(\R)$. A converg\^encia dessa diferença para zero é forte
indicativo da converg\^encia de $(f_n)$ e, conseqüentemente,
de um comportamento assintótico auto-similar. O cálculo da diferença
relativa também pode ser utilizado como um critério de parada
das iterações.
\insertfigure{fig:normheat}{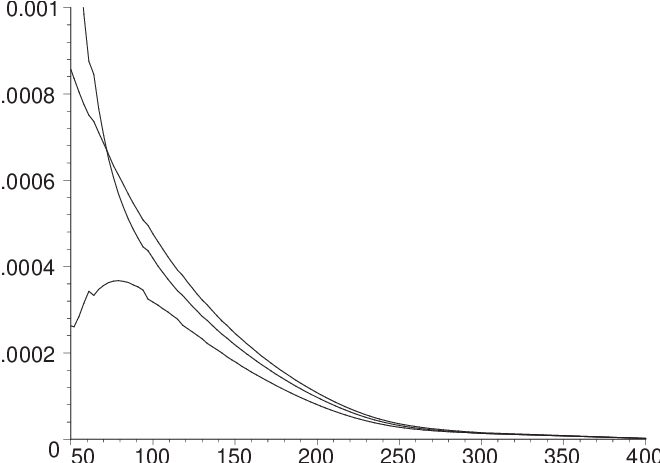}
{Curvas $(||f_n - f_{n-1}||/||f_n||) \times n$ para a equação
do calor com diferentes condições iniciais.}

\clearpage
\subsection{Equação de difusão não-linear}

O Teorema~\ref{teo:nldif} diz que soluções da equação
$$
  u_t = u_{xx} + \lambda u^a u_x^b u_{xx}^c,
$$
com condição inicial  $u(\cdot,0)=f\in\B_4$ e $a+2b+3c-3>0$,
se comportam como
\be
  \label{eq:ncomassnl}
  u(x,t) \approx \frac A{t^\alpha}\ \phi\left(\frac x {
  t^\beta}\right),
\ee
com $\alpha=\beta=1/2$, $\phi=\phi_*$, e $A$ é uma constante positiva
que pode depender de $\lambda$, $a$, $b$, $c$ e $c$.

Fizemos simulações numéricas em que variamos a condição inicial,
o parâmetro $\lambda$ e os expoentes $a$, $b$ e $c$
de maneira a satisfazer à condição  $a+2b+3c-3>0$.
As simulações são descritas pela Tabela~\ref{tab:simnl} e a universalidade de
$\alpha$ é ilustrada na Figura~\ref{fig:alphanl}.

Para todas as
simulações, verificamos que $\phi=\phi_*$, que o prefator converge
para algum número positivo e que a norma relativa da diferença
entre duas funções de perfil consecutivas converge para zero.

\inserttable{tab:simnl}{||c|c|c|c|c||}
{Descrição das simulações feitas no estudo das classes
  de universalidade do comportamento assintótico das soluções
  da equação de difusão não-linear, numeradas de $1$
  a $7$.}
{\hline\hline $n$&$\lambda$&$a$&$b$&$c$ \\ \hline\hline
\input{nltable} \hline \hline}

\insertfigure{fig:alphanl}{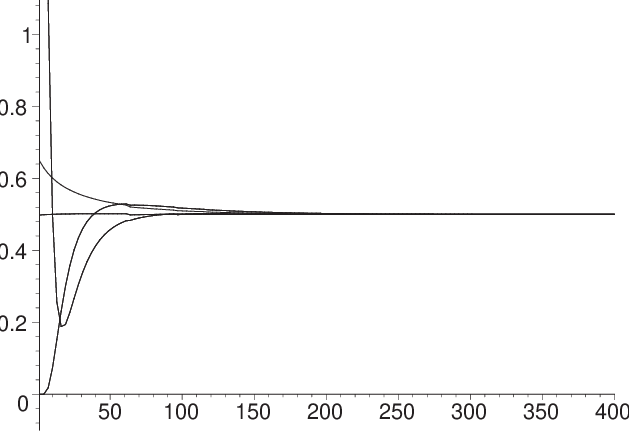}
{Gráfico de $\alpha_n\times n$ para as simulações descritas
  na Tabela~\ref{tab:simnl}, ilustrando a universalidade
  do expoente crítico $\alpha$ para perturbações irrelevantes
  da equação de difusão não-linear.}

\clearpage
\subsection{Equação de Barenblatt}

Estudamos também a equação de Barenblatt
\be
  \label{eq:barenblatt}
  u_t = [1 +\epsilon H(-u_t)]u_{xx},
\ee
onde $H(v)$ é a função de Heaviside com salto em $v=0$.
Resultados teóricos\footnote{
Ver \cite{gold92} e referências lá citadas.}
fornecem
uma aproximação de primeira ordem para o expoente $\alpha$ dada por
$$
  \alpha = \frac 1 2 + \frac \epsilon {\sqrt{2 \pi e}} + O(\epsilon^2).
$$
Simulações foram feitas para diferentes valores de $\epsilon$.
Os valores de $\alpha$ calculados numericamente são comparados
com o resultado teórico acima na Figura~\ref{fig:barenblatt}.

Em todas as simulações, verificamos que o prefator $A_n$ converge
a um número positivo e a norma relativa da diferença de duas funções
de perfil consecutivas converge a zero.

O estudo numérico desta equação também
ilustra uma característica interessante do RG numérico: é possível calcular
os expoentes críticos ainda que os mesmos não sejam conhecidos de forma
exa\-ta.
\insertfigure{fig:barenblatt}{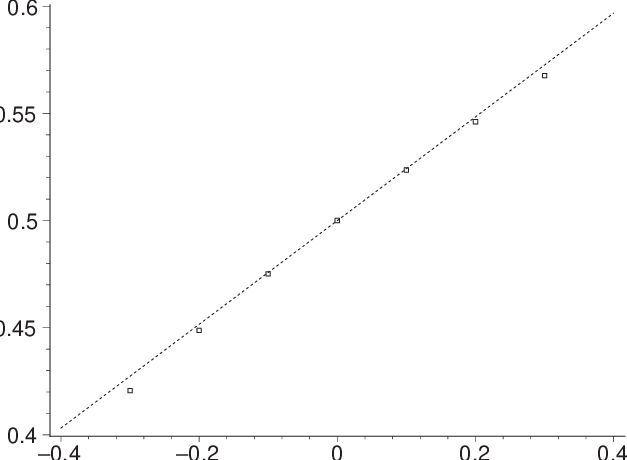}
{Curvas $\alpha\times\epsilon$ para diferentes simulações
  com a equação de Barenblatt, onde $\alpha$ é calculado numericamente
  e a reta pontilhada é dada pela aproximação de primeira ordem.}

\clearpage
\subsection{Refinamento da malha}

Para a equação do calor, fizemos um teste de refinamento de malha
utilizando a condição inicial $f_2$ e uma outra idêntica à $f_2$
porém com uma malha duas vezes mais fina. Os gráficos de $\alpha_n$
e $A_n$ são exibidos nas Figuras~\ref{fig:alpharefine}~e~\ref{fig:arefine},
respectivamente.
Na simulação com a malha refinada foi verificado que $\phi=\phi_*$ e que
a norma relativa da diferença de duas funções de perfil consecutivas
converge a zero.
\insertfigure{fig:alpharefine}{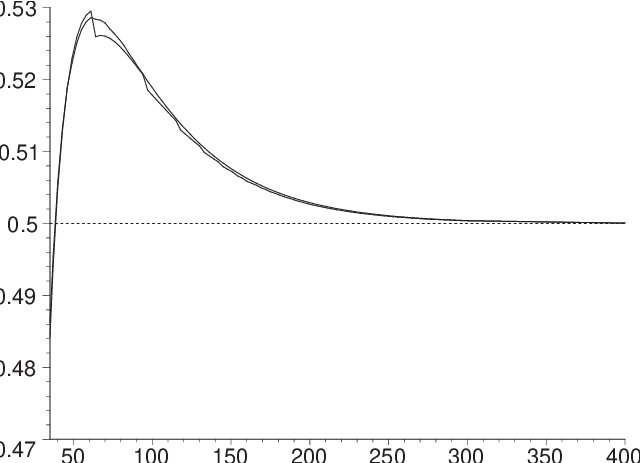}
{Curva $\alpha_n\times n$ para a equação do calor
com condição inicial $f_2$; uma curva foi gerada
com a malha duas vezes mais fina que a outra.}
\insertfigure{fig:arefine}{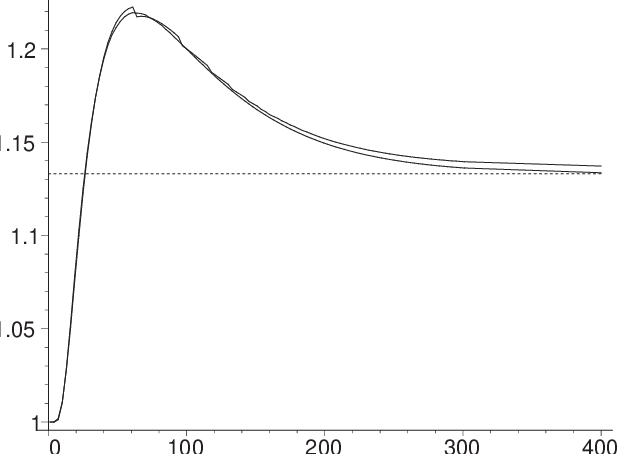}
{Curva $A_n\times n$ para a equação do calor
com condição inicial $f_2$; uma curva foi gerada
com a malha duas vezes mais fina que a outra.}

\clearpage
\section{Difusão com coeficiente periódico}
\label{sec:periodico}

As simulações discutidas nesta seção visam a corroborar a
Conjectura~\ref{conj:oscilante}.

Nas simulações realizadas, observamos a convergência de
$\alpha_n$ para $\alpha = 1/2$. Para verificar a universalidade deste
fato, utilizamos diferentes condições iniciais (Figura~\ref{fig:peralphaic}), perturbações
não-lineares (Figura~\ref{fig:peralphanl}) e valores de $\mu$ (Figura~\ref{fig:peralphamu}).

Em todas as simulações, a convergência da norma relativa da diferença
entre funções de perfil consecutivas convergiu para zero e o prefator
$A_n$ convergiu para um número positivo.

\insertfigure{fig:peralphaic}{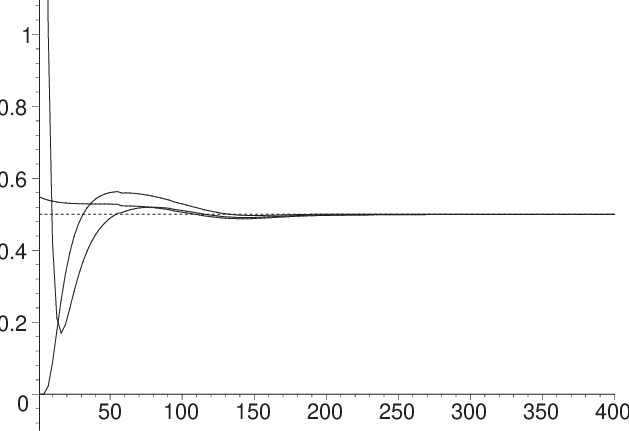}
  {Curvas $\alpha_n\times n$ para o problema de difusão
  com coeficiente periódico e diferentes condições iniciais.}
\insertfigure{fig:peralphamu}{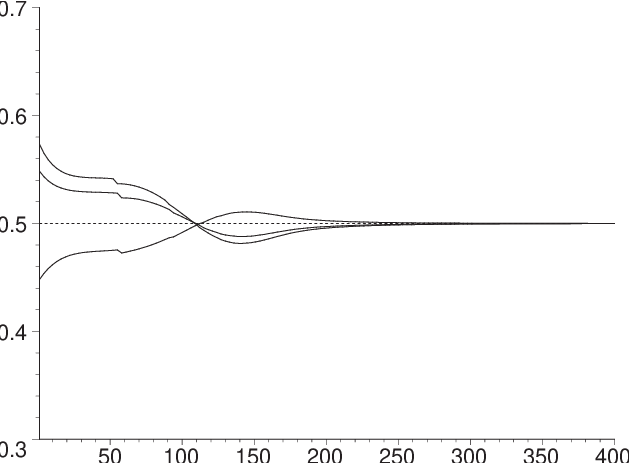}
  {Curvas $\alpha_n\times n$ para o problema de difusão
  com coeficiente periódico e diferentes valores de $\mu$.}
\insertfigure{fig:peralphanl}{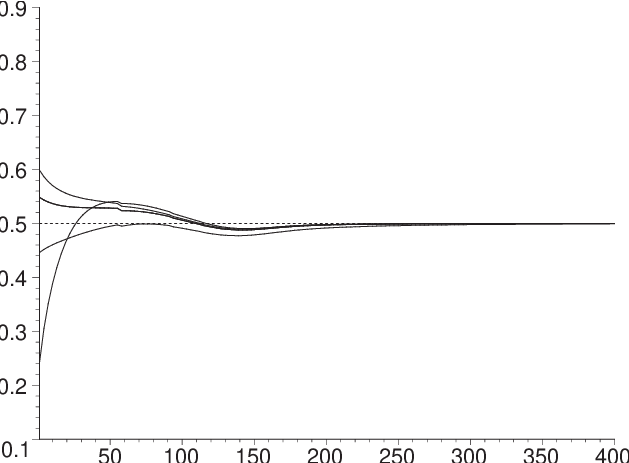}
  {Curvas $\alpha_n\times n$ para o problema de difusão
  com coeficiente periódico e diferentes perturbações
  não-lineares irrelevantes.}

\clearpage
\subsection{Homogeneização do coeficiente de difusão}

A principal diferença desta equação em relação \`a equação
do calor, ou \`a equação de difusão não-linear,
é a mudança da função de perfil causada pelo termo $\mu g(x)$.
Por isso, plotamos as curvas $ -\ln(\phi(x)) \times x^2/4$ para diversos
experimentos. O fato de estas curvas serem retas com inclinações
$\sigma^{-1}$ significa que $\phi = \phi_\sigma$ para cada uma delas
-- veja (\ref{eq:logphisigma}). A Figura~\ref{fig:perlglg} mostra essa
curva para $\mu = 0,8$ e $g=\cos$,
juntamente com o valor de $\sigma=0,6$ esperado teoricamente.

\insertfigure{fig:perlglg}{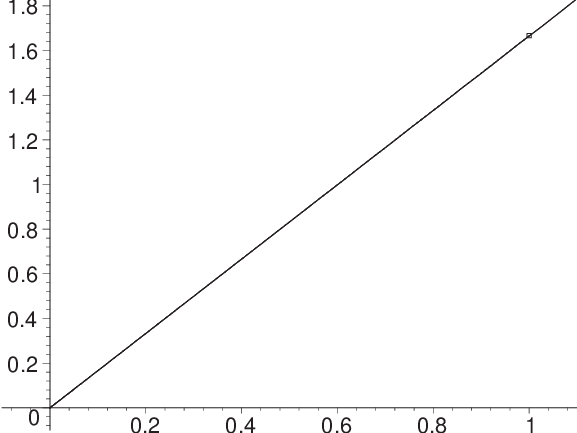}
{Curva parametrizada $(-\ln\phi,x^2/4)$, cujo traço
  é uma reta com inclinação $\sigma^{-1}$.}

\clearpage
\section{Difusão dependente do tempo}
\label{sec:deptempo}

Nesta seção vamos estudar os resultados apresentados na
Seção~\ref{sec:anadeptempo} do Capítulo~\ref{chap:rgana}.
Esses resultados consideram o problema de difusão com
coeficiente dependente do tempo. O caso linear é coberto pelo
Teorema~\ref{teo:compu-0}, enquanto para o caso não-linear
foi formulada a Conjectura~\ref{conj:compassgeral}. Começamos pelo caso
linear.

\subsection{Caso linear}

Considere a equação
$$
  u_t= (t^p + \delta t^r) u_{xx},
$$
onde $r<p$, com condição inicial $u(\cdot,1)=f\in\B_4$.
O Teorema~\ref{teo:compu-0} diz que as soluções deste problema
se comportam assintoticamente como
$$
  u(x,t) \approx \frac A{t^\alpha}\ \phi\left(\frac x {
  t^\beta}\right),
$$
com $\alpha=\beta=(p+1)/2$, $\phi=\phi_p$ e $A=[(p+1)/4\pi]^{1/2}
\int_\R f(y)dy$.

Para estudar numericamente o resultado acima, realizamos um conjunto de
$18$ simulações, que correspondem à combinação de três condições
iniciais diferentes com coeficientes de difusão de seis formas diferentes,
dadas por:
$t^{1/2} + t^{1/4}$; $t^{1/2} + 1$; $t^{1/2}$; $t + t^{3/4}$; $t + 1$; e $t$.

A Figura~\ref{fig:tpalpha} mostra a convergência dos expoentes
$\alpha$ para o valor $(p+1)/2$ e sua universalidade com
respeito ao termo $\delta t^r$ e à condição inicial $f$.
Os valores de $p$ utilizados são $1/2$ e $1$,
determinados pelos coeficientes $K(t)$ listados no parágrafo
anterior.
\insertfigure{fig:tpalpha}{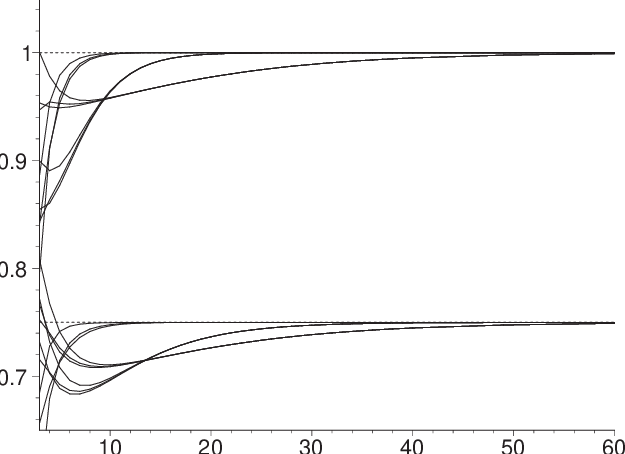}
{Curvas $\alpha_n \times n$ para a equação linear com coeficiente
dependente do tempo para diferentes coeficientes e condições inicias.
Os valores teóricos que correspondem a $p=1/2$ e $p=1$ são,
respectivamente, $\alpha=0,75$ e $\alpha=1,0$.}

O prefator $A$ é também universal com respeito a detalhes
nas condições iniciais e depende apenas da integral das mesmas.
A Figura~\ref{fig:tpprefactor} mostra a convergência
do prefator para os valores teóricos $A = \sqrt{(p+1)/4\pi} \int_\R f(y)dy$.
\insertfigure{fig:tpprefactor}{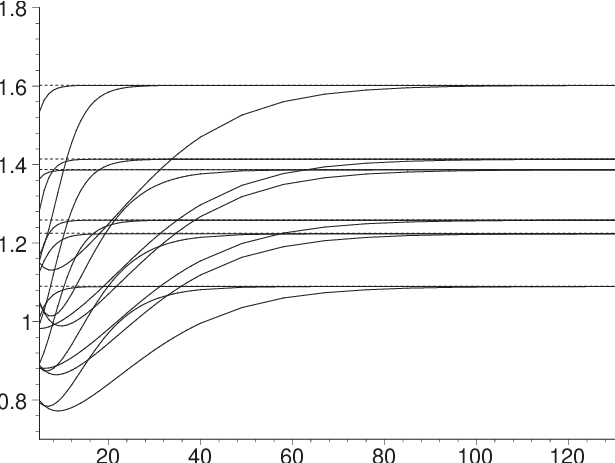}
{Curvas $A_n \times n$ para a equação com coeficiente dependente do
tempo para diferentes coeficientes e condições inicias.
As linhas pontilhadas correspondem
aos valores teóricos.  Estes valores, dados por
$[(p+1)/4\pi]^{1/2}\int_\R f(y)dy$, são
1,090, 1,225, 1,259, 1,388, 1,414 e 1,602.}

Por ultimo, a Figura~\ref{fig:tpprofile} revela que a função de perfil
é, de fato, dada por $\phi=\phi_p$.
\insertfigure{fig:tpprofile}{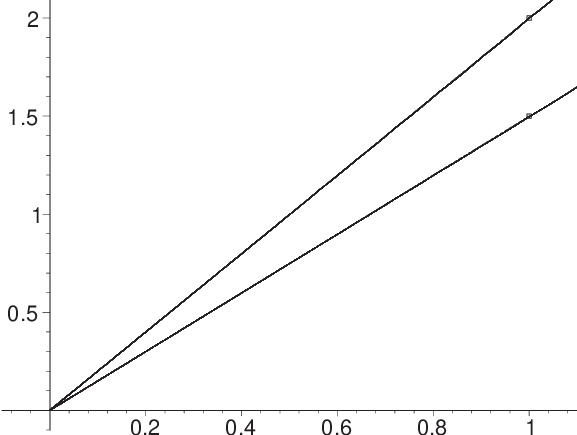}
{Curvas parametrizadas $-\ln\phi \times x^2/4$
para a equação com coeficiente dependente do
tempo com diferentes coeficientes e condições inicias.
A inclinação da reta deve ser dada por $p+1$. Cada curva visível
no gráfico é na verdade a superposição de nove curvas
geradas pelas simulações descritas acima.
As linhas pontilhadas e o ponto marcado $(1,p+1)$ correspondem
aos valores teóricos.}

\clearpage
\subsection{Caso não-linear}

Consideramos a equação
$$
  u_t= (t^p + \delta t^r) u_{xx} + \lambda u^au_x^bu_{xx}^c.
$$
Todas as perturbações não-lineares
que foram classificadas como irrelevantes para a equação de difusão
não-linear, ou seja, aquelas que satisfazem $d_F = a+2b+3c-3 >0$,
ainda são irrelevantes para a equação com coeficiente de difusão
dependente do tempo para $p\geq 0$, pois satisfazem
$\eta_F = a+2b+3c-(3+p)/(1+p) >0$.

No caso em que $p>0$, as perturbações que antes eram classificadas
como marginais, agora passam a ser também irrelevantes.
Essas perturbações podem ser $F=u^3$, $F=uu_x$ ou $F=u_{xx}$.

A Conjectura~\ref{conj:compassgeral} diz que o comportamento
assintótico das soluções da equação acima ainda é o mesmo dado pela equação
linear tratada pelo Teorema~\ref{teo:compu-0}, exceto pelo fato de que
o prefator $A$ é uma constante positiva que depende da condição
inicial, do coeficiente de difusão e da perturbação não-linear,
desde que a perturbação não-linear $F$ satisfaça $\eta_F>0$.

Os gráficos relativos às simulações feitas nesta seção são omitidos
por serem todos análogos aos que já foram apresentados em seções
anteriores. Por esta razão vamos nos limitar a exibir a
Tabela~\ref{tab:simtpnl}, que descreve as simulações realizadas
visando a confirmar a afirmação do parágrafo anterior.

\inserttable{tab:simtpnl}{||c|c|c|c|c|c||}
{Descrição das simulações feitas no estudo das classes
  de universalidade do comportamento assintótico das soluções
  da equação de difusão não-linear com coeficiente de difusão
  dependente do tempo, numeradas de $1$ a $13$.}
{\hline\hline $n$&$p$&$\lambda$&$a$&$b$&$c$ \\ \hline\hline
\input{tpnltable} \hline \hline }

\clearpage
\subsection{Transição de fase}

Vamos estudar numericamente a equação
\be
  \label{eq:phtivp}
  u_t= t^p u_{xx} - u^a,
\ee
onde $a$ e $p$ são parâmetros reais tais que $a>1$
e $p\geq 0$.

Dado $p\geq 0$, o expoente $a$ que corresponde
a uma perturbação marginal é \footnote{
Faça $\eta = 0$ em (\ref{eq:eta}).
}
\be
  \label{eq:mc}
  a = \frac{p+3}{p+1}.
\ee
No plano $(p,a)$, a curva~(\ref{eq:mc}) é formalmente a interface
entre a \emph{região relevante} $ a > (p+3)(p+1)^{-1} $
e a \emph{região irrelevante} $ a < (p+3)(p+1)^{-1}$.
O estudo da função de perfil, do expoente crítico $\alpha$
e dos prefatores, indica que~(\ref{eq:mc}) é de fato
a \emph{curva crítica} para a \emph{transição de fase} entre
os estados relevante e irrelevante.

A Figura~\ref{fig:pht} mostra o espaço de parâmetros $(p,a)$
com a curva crítica (ou separatriz) juntamente com as simulações
que fizemos neste trabalho. A curva $ a_c(p) = (p+3)(p+1)^{-1}$
corresponde aos valores de $(p,a)$ que tornam a equação~(\ref{eq:phtivp})
invariante à renormalização com expoentes
$\beta = (p+1)/2$ e $\alpha = (p+1)/2 =  1/(a-1)$.
\insertfigure{fig:pht}{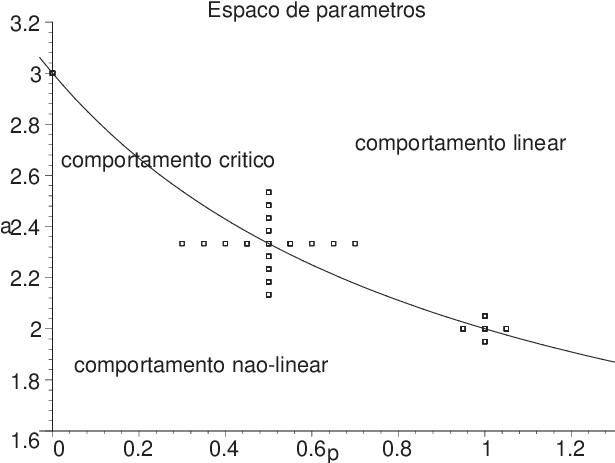}
{Espaço de parâmetros $(p,a)$ para a equação
(\ref{eq:phtivp}). A curva crítica separa os regimes linear e
não-linear da equação. Os quadrados representam os parâmetros
para os quais foram feias as simulações numéricas.}

Para os valores de $(p,a)$ acima da curva crítica, o operador
não-linear tende a zero sob a ação de sucessivas renormalizações
com $\beta = \alpha = (p+1)/2$, enquanto a parte linear permanece
invariante. Neste caso, o regime assintótico das soluções
é essencialmente da equação linear.

Para os valores de $(p,a)$ abaixo da curva crítica,
a equação~(\ref{eq:phtivp})
permanece invariante para $\beta =  (p+1) / 2$ and $\alpha =  1/(a-1)$.
Neste caso o expoente de decaimento
é determinado pela parte não-linear do operador,
enquanto a difusão é devida ao termo linear, e a função
de perfil depende de ambos.

A Figura~\ref{fig:pht2} mostra o gráfico da função
\be
  \label{eq:alphapa}
  \alpha(p,a) = \max \left\{\frac{1+p}2,\frac 1{a-1}\right\}.
\ee
Olhando para as curvas de nível da superfície gerada,
é possível perceber que as singularidades
da superfície $(p,a,\alpha)$ formam uma curva que
corresponde aos parâmetros $p$ e $a$ da curva crítica.
É possível também perceber que de um lado dessa curva
o expoente $\alpha$ é determinado por $p$ enquanto
do outro lado é determinado por $a$.
\insertfigure{fig:pht2}{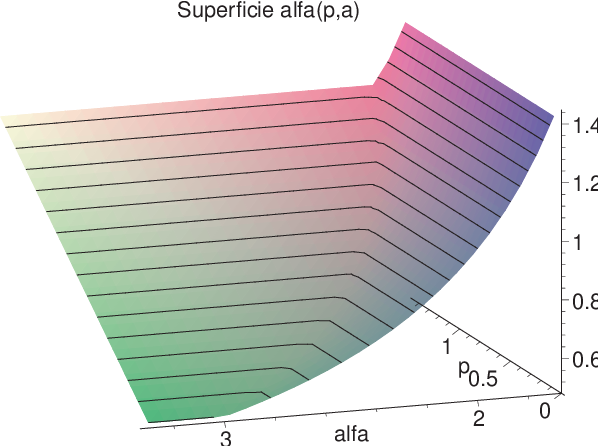}
{Superfície $\alpha(p,a)$ dada por (\ref{eq:alphapa}).
As singularidades da superfície correspondem à curva crítica sobre
a qual ocorre a transição de fase.}

Com o objetivo de verificar a relação~(\ref{eq:alphapa}),
na Figura~\ref{fig:aalphacurve} plotamos
a interseção da superfície (\ref{eq:alphapa}) com o plano
$p = 1/ 2$, juntamente com os expoentes $\alpha$
calculados pelas nossas simulações.
Além disso, a Tabela~\ref{tab:aalpha} compara os valores
teóricos com os valores computados de $\alpha$ como uma
função de $a$ para $p=1/2$ fixo.
A Figura~\ref{fig:palphacurve} e a Tabela~\ref{tab:palpha} são
análogos para $a =  7/ 3$ fixo e $\alpha$ como uma função $p$.
\insertfigure{fig:aalphacurve}{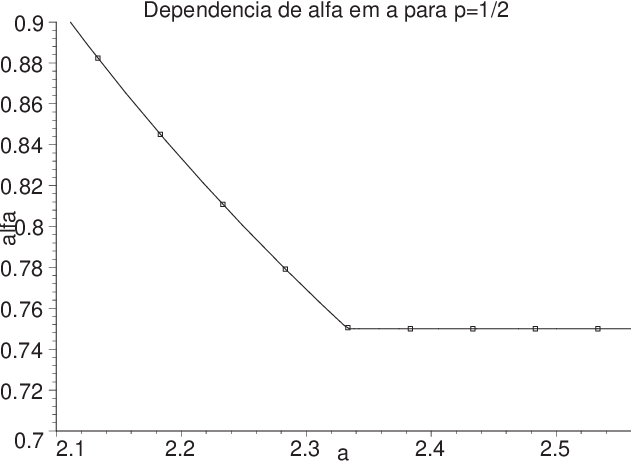}
{Gráfico de $\alpha\times a$ para $p=1/2$ fixo. Os quadrados
representam os valores calculados numericamente e a linha contínua
é dada pela equação (\ref{eq:alphapa}).}
\inserttable{tab:aalpha}{||c|c|c||}
{Valores de $\alpha$ como uma função de
$a$ para $p=1/2$ fixo. O valor calculado numericamente é comparado com
o valor dado pela equação~(\ref{eq:alphapa})}
{\hline\hline $a$&$\alpha$ teórico&$\alpha$ calculado \\ \hline\hline
\input{aalphatable} \hline \hline}
\insertfigure{fig:palphacurve}{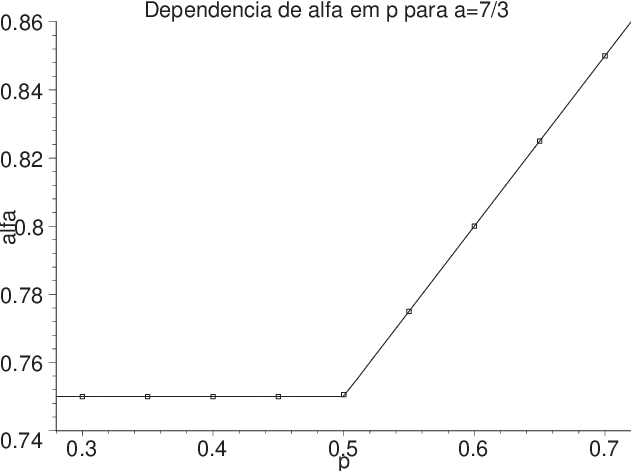}
{Gráfico de $\alpha\times p$ para $a=7/3$ fixo. Os quadrados
representam os valores calculados numericamente e a linha contínua
é dada pela equação (\ref{eq:alphapa}).}
\inserttable{tab:palpha}{||c|c|c||}
{Valores de $\alpha$ como uma função de
$p$ para $a=7/3$ fixo. O valor calculado numericamente é comparado com
o valor dado pela equação~(\ref{eq:alphapa})}
{\hline\hline $p$&$\alpha$ teórico&$\alpha$ calculado \\ \hline\hline
\input{palphatable} \hline \hline}

As Figuras~\ref{fig:alphaiii}~e~\ref{fig:aiii}
ilustram a diferença na convergência de $(\alpha_n)$
e $(A_n)$.
Fazendo $p= 1/ 2$ e escolhendo três valores diferentes
para $a$ (um crítico, um subcrítico e supercrítico),
podemos investigar a diferença dos regimes, no que diz
respeito os principais aspectos qualitativos de seu comportamento
assintótico.

Na Figura~\ref{fig:alphaiii} podemos ver a diferença
no cálculo do expoente $\alpha$. No caso supercrítico
(termo não-linear irrelevante), o expoente de decaimento
é determinado por $\alpha=
[(p+1)/ 2] =  3/ 4$, enquanto que no caso subcrítico este expoente
assume um valor maior. Além disso, quando usamos parâmetros
críticos, a convergência $(\alpha_n)$ para $\alpha= [(p+1)/ 2] =  3/ 4$
é muito mais lenta que nos outros casos (apesar de não aparecer
na figura, após um número muito maior de iterações podemos
ver que a convergência ocorre de fato).

Na Figura~\ref{fig:aiii} vemos a convergência do prefator
nos três regimes. Para parâmetros na região relevante e
irrelevante, o prefator converge rapidamente a um número
positivo, enquanto no caso crítico o prefator
converge para zero quando $n\to\infty$. Este fato é
devido ao aparecimento de uma correção logarítmica
ao comportamento auto-similar, que será discutido
logo adiante.

\insertfigure{fig:alphaiii}{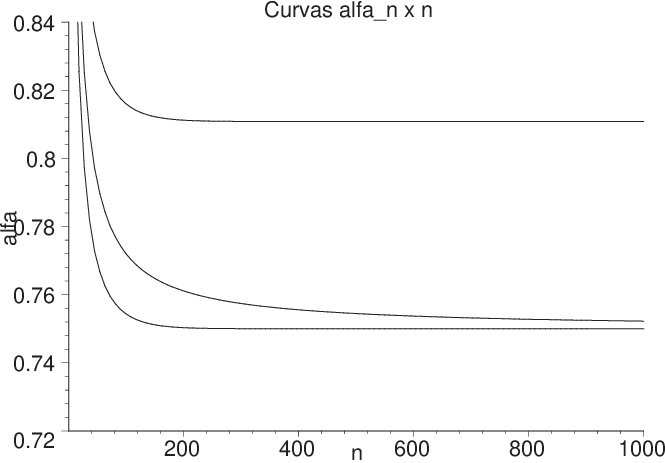}
{Curvas $\alpha_n \times n$ para três regimes distintos.
Fixado o parâmetro $p=1/2$, cada curva foi gerada com um valor de
$a<7/3$, $a=7/3$ e $a>7/3$,  correspondendo aos regimes não-linear,
crítico e linear, respectivamente.}
\insertfigure{fig:aiii}{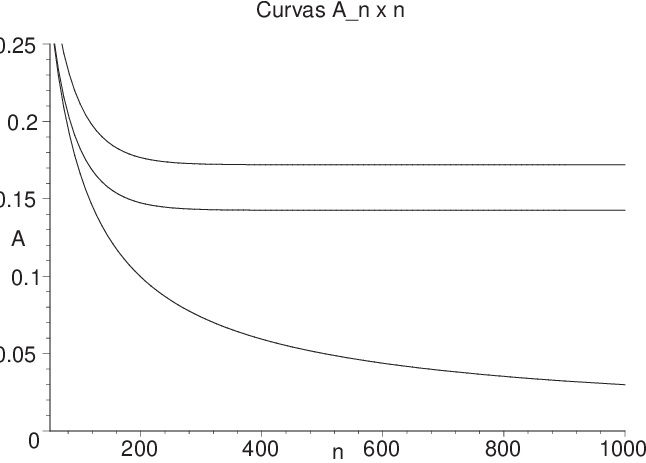}
{Curvas $A_n \times n$ para três regimes distintos.
Fixado o parâmetro $p=1/2$, cada curva foi gerada com um valor de
$a<7/3$, $a=7/3$ e $a>7/3$,  correspondendo aos regimes não-linear,
crítico e linear, respectivamente.}

Na Figura~\ref{fig:nla} ilustramos a universalidade do prefator
com respeito à condição inicial no caso de perturbações
relevantes ($a$ subcrítico). Fixamos $p=1/2$ e escolhemos
três valores distintos de $a$ satisfazendo $a<a_c(p)$.
O prefator obtido como $\lim_{n\to\oo}A_n$ depende apenas
da escolha de $a$ e não das condições iniciais.
\insertfigure{fig:nla}{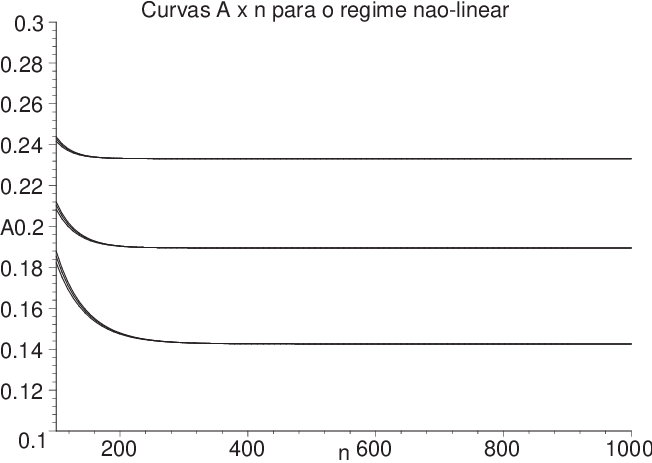}
{Curvas $A_n \times n$ para o regime não-linear,
com três valores distintos de $a$ satisfazendo $a<7/3$ e
três condições iniciais distintas, fixado $p=1/2$. O prefator dado por
$A=\lim_{n\to\oo} A_n$ depende apenas de $a$.}

A Figura \ref{fig:lineara} mostra a não-universalidade
do prefator no caso da perturbação irrelevante. Para
$p=1/2$ fixo, escolhemos três valores de $a>a_c(p)$ e
três condições iniciais distintas. O prefator obtido
é influenciado tanto pelas escolha de $a$ quanto da
condição inicial.
\insertfigure{fig:lineara}{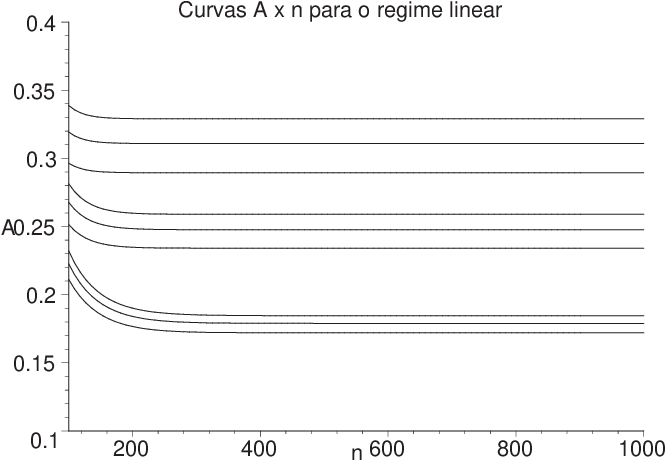}
{Curvas $A_n \times n$ para o regime linear,
com três valores distintos de $a$ satisfazendo $a>7/3$ e
três condições iniciais distintas, fixado $p=1/2$. O prefator dado por
$A=\lim_{n\to\oo} A_n$ depende do parâmetro $a$ e da condição
inicial utilizada.}

\clearpage
\subsection{Correções logarítmicas sobre a curva crítica}

Para a equação~(\ref{eq:phtivp}) com parâmetros sobre a
curva crítica, o comportamento assintótico auto-similar
esperado sofre uma correção logarítmica:
\be
  \label{eq:compasscritlog}
  u(x,t) \approx \frac A {(t \ln t)^\alpha}\ \phi \left(\frac B
  {t^\beta} x \right)
  \mbox{ \ quando \ } t \to \infty,
\ee
onde $\beta=(p+1)/2$ e $\alpha=(p+1)/2=1/(a-1)$.
Além disso, resultados teóricos dizem que o prefator
$A$ não depende da condição inicial, mas apenas da
escolha dos parâmetros $(p,a)$.
Lembrando que o algoritmo faz com que $f_n(0)=1$
no processo construtivo, combinando (\ref{eq:compasscritlog})
com (\ref{eq:prefatores2}) obtemos
\be
\nonumber
\label{eq:lglgrel}
  \ln A_n \approx \ln \frac A {L^\alpha} - \alpha \ln n.
\ee
Essa relação diz que o gráfico de $\ln A_n \times \ln n$
se aproxima, para $n$ suficientemente grande,
de uma linha reta tendo inclinação dada por $-\alpha$
que intercepta o eixo vertical no ponto $\ln A / {L^\alpha}$.
Neste caso, a universalidade de $A$ traduz-se no fato de que
quaisquer condições iniciais dentro de uma classe apropriada
dão origem a curvas que são assintóticas à mesma reta
e não a retas paralelas.

As afirmações do parágrafo acima foram verificadas para
$p=0$, $p= 1/2$, e $p=1$.
A Figura~\ref{fig:loglog1} mostra o gráfico de
$\ln A_n \times \ln n$ para $p = 1 /2$ e três condições
iniciais diferentes. A inclinação da reta assintótica é
$-\alpha$, sendo utilizado o valor de $\alpha$ calculado
pelo RG.
\insertfigure{fig:loglog1}{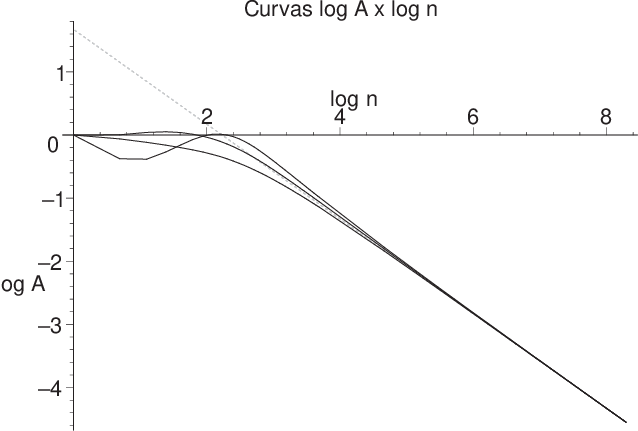}
{Curvas $\ln A_n \times \ln n$ para $p=1/2$
e $a=7/3$, com três condições iniciais distintas. A reta pontilhada
passa pelo ponto $(\ln A_n,\ln n)$ da última iteração do algoritmo
e tem inclinação $-\alpha$, considerando o expoente calculado
na última iteração.}

%% file: nltable.tex
1 & .10 & 3 & 1 & 0 \\
\hline 2 & .20 & 3 & 1 & 0 \\
\hline 3 & -.10 & 1 & 1 & 1 \\
\hline 4 & .10 & 0 & 1 & 1 \\
\hline 5 & .15 & 3 & 1 & 0 \\
\hline 6 & -.20 & 1 & 1 & 1 \\
\hline 7 & .30 & 1 & 0 & 1 \\
\hline 

%% file: tpnltable.tex
1 & .50 & .10 & 1 & 1 & 0 \\
\hline 2 & .50 & .10 & 3 & 0 & 0 \\
\hline 3 & .50 & .10 & 0 & 0 & 1 \\
\hline 4 & .25 & .10 & 1 & 1 & 0 \\
\hline 5 & .25 & .10 & 3 & 0 & 0 \\
\hline 6 & .25 & .10 & 0 & 0 & 1 \\
\hline 7 & .50 & .10 & 3 & 1 & 0 \\
\hline 8 & .50 & .20 & 3 & 1 & 0 \\
\hline 9 & .50 & -.10 & 1 & 1 & 1 \\
\hline 10 & .50 & .10 & 0 & 1 & 1 \\
\hline 11 & .50 & .15 & 3 & 1 & 0 \\
\hline 12 & .50 & -.20 & 1 & 1 & 1 \\
\hline 13 & .50 & .30 & 1 & 0 & 1 \\
\hline 

%% file: aalphatable.tex
2.133333 & .882353 & .882353 \\
\hline 2.183333 & .845070 & .845071 \\
\hline 2.233333 & .810811 & .810811 \\
\hline 2.283333 & .779221 & .779221 \\
\hline 2.333333 & .750000 & .750553 \\
\hline 2.383333 & .750000 & .749991 \\
\hline 2.433333 & .750000 & .749991 \\
\hline 2.483333 & .750000 & .749991 \\
\hline 2.533333 & .750000 & .749991 \\
\hline 

%% file: palphatable.tex
.300000 & .750000 & .750000 \\
\hline .350000 & .750000 & .750000 \\
\hline .400000 & .750000 & .750000 \\
\hline .450000 & .750000 & .750001 \\
\hline .500000 & .750000 & .750553 \\
\hline .550000 & .775000 & .774994 \\
\hline .600000 & .800000 & .799996 \\
\hline .650000 & .825000 & .824987 \\
\hline .700000 & .850000 & .849962 \\
\hline 

%% file: homog.tex
\chapter{Teoria de Homogeneização}
\label{chap:homog}

Neste capítulo vamos apresentar a teoria de homogeneização à qual nos
referimos nos capítulos anteriores. A principal idéia na teoria de
homogeneização é a de tentar descrever um meio não-homogêneo a partir de
seu comportamento \emph{efetivo}. Essa tentativa baseia-se na hipótese de
que existe um meio homogêneo cujas propriedades são próximas daquelas do
meio real quando medidas em grandes escalas. O processo de ``aproximar
pela média'', chamado \emph{homogeneização}, substitui uma estrutura
complicada em escala microscópica por uma estrutura homogênea
assintoticamente equivalente. A teoria que vamos apresentar neste capítulo
diz de uma forma matematicamente precisa que a substituição da equação
original por outra com coeficiente constante é uma aproximação válida em
relação a um certo limite.

Escrevemos esse capítulo com dois objetivos. O primeiro é chegar ao
Teorema~\ref{teo:parab} e ao Co\-ro\-lá\-rio~\ref{cor:parab}, que
justificam a substituição do coeficiente não-homogêneo por um homogêneo no
estudo do Problema~\ref{prob:periodico} da Introdução. O segundo é
apresentar ao leitor um texto auto-contido sobre a teoria de
homogeneização, para equações elípticas e parabólicas lineares, em espaços
de Sobolev. O conteúdo deste capítulo se baseia no artigo de
Papanicolaou~\cite{papa02} e no livro de Jikov~et~al.~\cite{jikov94}.

Vamos fazer uso indiscriminado dos resultado básicos da
teoria de medida e integração em subconjuntos de $\R^m$.
Serão utilizados também alguns elementos de análise funcional como
os conceitos de espaços de Hilbert e de Banach,
ou o Lema de Riez. Contudo, não é necessário que o leitor
tenha feito um curso de análise funcional desde que tenha algum
livro introdutório para consultar.

Na Seção~\ref{sec:unid} vamos apresentar as idéias
básicas da técnica de homogeneização e estudar um caso simples em
que não são necessários muitos resultados matemáticos mas que,
por sua vez, contribuem para uma visão simples e clara do
problema de homogeneização.

Na Seção~\ref{sec:principal} vamos definir alguns espaços de Sobolev
e suas principais propriedades, trabalhar com seus duais e estudar
a equação elíptica linear com condição de contorno periódica. Isso feito,
introduzimos um problema auxiliar, que prepara o terreno para o
Teorema~\ref{teo:parab}, em que uma equação parabólica mais geral
do que o Problema~\ref{prob:periodico} linear é estudada.

\section{Estudo de um problema unidimensional}
\label{sec:unid}

Nesta seção pretendemos apresentar o estudo de um caso bastante
simples com o intuito de dar um enfoque físico e intuitivo
da técnica de homogeneização.
Por esta razão, optamos por deixar em segundo plano as dificuldades
matemáticas (utilizando para isso o recurso de acrescentar mais hipóteses
do que o necessário aos problemas), e demos preferência à tentativa de
justificar e ilustrar a idéia de homogeneizar um meio heterogêneo.

\subsubsection{O problema elíptico linear unidimensional}

Considere a equação abaixo, que modela o estado de equilíbrio de
um condutor (de calor) unidimensional ao qual aplicamos um gradiente
de temperatura $V$:
$$
\left\{
  \ba{l}\displaystyle
  \Der{}x \left[D(x)\Der{u(x)}x\right] = 0, \quad x\in[-N,N] \\
  u(-N)=-VN,\ u(+N)=+VN
  \ea
\right.,
$$
onde $N>0$ é a metade do tamanho do condutor que estamos considerando e
$D(x)$ é uma função periódica (de período $1$), limitada
e satisfaz $\inf_{x\in\R} D(x)>0$.
A função $D(x)$ representa a condutividade térmica no ponto $x$ e
$V$ representa o gradiente de temperatura aplicado ao condutor.

O fluxo em um ponto $x$ é dado por $D(x)\Der{u(x)}x$;
o fluxo médio, por unidade de gradiente de temperatura,
será denotado por $D_N$ e dado pela fórmula
$$
  D_N = \frac 1V \left(\frac 1 {2N}\int_{-N}^{+N}D(x)\Der{u(x)}x dx\right).
$$
No caso unidimensional, que estamos tratando neste momento,
o fluxo $D(x)u_x(x)$ é constante, logo a integral acima é desnecessária.
Entretanto, essa representação tem o significado físico de
estarmos tomando a média do fluxo. Em dimensões maiores,
ao invés de um fluxo constante temos um fluxo com divergente nulo.

Queremos determinar a \emph{condutividade efetiva} $D^*$,
que é definida pela existência do limite
$$
  D^* = \lim_{N\to\oo} D_N.
$$

Ressaltando novamente que estamos tratando do caso unidimensional,
temos uma fórmula explícita para $D_N$. Cada equação abaixo é conseqüência
imediata da equação que a precede:
\begin{eqnarray*}
 D_N &=& \frac1V\, D(x)u_x(x)\quad\forall\ x\in[-N,N]\\
 u_x &=& V\,D_N/D(x)   \\
 {2V\,N} = \int_{-N}^{+N} u_xdx &=& V\,D_N \int_{-N}^{+N} \frac 1{D(x)}dx \\
 D_N &=& \left( \frac 1 {2N} \int_{-N}^{+N}\frac 1{D(x)}dx \right)^{-1}.
\end{eqnarray*}
Como $D$ é periódica, de período $1$, por exemplo,
a última igualdade implica em
$$
  D^* = \lim_{N\to\oo} D_N = \left(\int_0^1\frac1{D(x)}dx\right)^{-1}
  = \av {D^{-1}}^{-1}.
$$
Uma interpretação física para a fórmula acima é que a resistividade
efetiva $(D^*)^{-1}$ é a resistividade média $\av{D^{-1}}$ em um
sistema de pequenos resistores em série. Observamos que o mesmo
raciocínio pode se aplicar ao caso em que a conditividade $D$
não é periódica ou até mesmo quando é um campo estocástico, exceto pela
integral desta última fórmula.

\subsubsection{Uma forma equivalente de ver o mesmo problema}

O mesmo problema pode ser visto sem que se aumente o tamanho
do condutor no limite $N\to\oo$, como fizemos acima. Podemos ajustar
nossa escala de medida ao tamanho do problema e pensar que, nessa escala,
a heterogeneidade é que passa a ter uma oscilação cada vez mais rápida.
Isso é feito através de uma mudança de variáveis. Primeiramente fazemos
$$
  \epsilon = \frac 1 N
$$
e definimos
$$
  D^\epsilon(x) = D(x/\epsilon) \qquad \mbox{e} \qquad u^\epsilon(x) =
  \epsilon u(x/\epsilon).
$$
A nova equação a ser considerada será
$$
\left\{
  \ba{l}
  \left[D^\epsilon(x)u^\epsilon_x(x)\right]_x = 0\quad x\in[-1,1]\\
  u^\epsilon(\pm1) = \pm1
  \ea
\right.
$$
e o coeficiente efetivo dado por
$$
  D^* = \lim_{\epsilon\downarrow0}\frac 12 \int_{-1}^{1}
  D^\epsilon(x) \left[ \Der {u^\epsilon(x)}x \right] dx.
$$
Com essa abordagem obtemos o mesmo valor para $D^*$.

O problema que vamos considerar em seguida ilustra bem a idéia
de que a aproximação do meio heterogêneo pelo meio homogêneo
é uma ``boa'' aproximação.

\subsubsection{Um problema de homogeneização}

Considere $D\in C^1(\R)$ positiva e periódica de período $1$.
Para cada $\epsilon>0$ defina $D^\epsilon(x) =
D(x/\epsilon)$ e seja $u^\epsilon$ a solução de
\be
  \label{eq:dirichlet1d}
  \left\{
  \ba{l}
  -\left[D^\epsilon(x)u^\epsilon_x(x)\right]_x = f,\quad x\in[-1,1]\\
  u^\epsilon(\pm1) = 0
  \ea
\right.,
\ee
onde $f\in C^1([-1,1])$.

Seja $u^0$ a solução do mesmo problema
substituindo-se $D^\epsilon$ por
\be
  \label{eq:d*}
  D^* = \left(\int_0^1\frac1{D(x)}dx\right)^{-1}
  = \av {D^{-1}}^{-1}.
\ee

O restante desta seção é dedicado ao propósito de oferecer uma prova simples
e elegante da seguinte proposição.

\begin{prop} Com $u^\epsilon$ e $u^0$ definidas acima, vale o limite
\be
  \label{eq:lim1d}
  u^\epsilon(x) \toe u^0(x)
\ee
uniforme em $[-1,1]$.
\end{prop}

A equação unidimensional (\ref{eq:dirichlet1d}) pode ser resolvida explicitamente.
Integrando a equação de $0$ a $x$ obtemos
$$
  D^\epsilon(x)u^\epsilon_x(x) = [D^\epsilon(0)u^\epsilon_x(0)] - \int_0^xf(t)dt.
$$
Definindo $F(x) = \int_0^xf(t)dt$, temos
$$
  u^\epsilon_x(x) = \frac1{D^\epsilon(x)}[D^\epsilon(0)u^\epsilon_x(0)]
  - \frac1{D^\epsilon(x)}F(x).
$$
Integrando novamente, obtém-se
$$
  u^\epsilon(x) = u^\epsilon(0) + [D^\epsilon(0)u^\epsilon_x(0)]
  \int_0^x\frac1{D^\epsilon(t)}dt - \int_0^x\frac{F(t)}{D^\epsilon(t)}dt.
$$
Definimos $c_1(\epsilon)=[D^\epsilon(0)u^\epsilon_x(0)]$ e
$c_2(\epsilon)=u^\epsilon(0)$, donde
\be
  \label{eq:ueps}
  u^\epsilon(x) = c_2(\epsilon) + c_1(\epsilon)
  \int_0^x\frac1{D^\epsilon(t)}dt - \int_0^x\frac{F(t)}{D^\epsilon(t)}dt.
\ee
A condição de contorno é $u^\epsilon(\pm1)=0$.
Subtraindo ou somando as equações
$u^\epsilon(-1)=0$ e $u^\epsilon(+1)=0$ obtemos, respectivamente,
$$
  c_1(\epsilon) = \left(\int_{-1}^{+1}\frac1{D^\epsilon(t)}dt\right)^{-1}
  \int_{-1}^{+1}\frac{F(t)}{D^\epsilon(t)}dt
$$
e
$$
  c_2(\epsilon) = \frac12\left(\int_{-1}^{+1}\frac{\sgn(t)}{D^\epsilon(t)}F(t)dt
    - c_1(\epsilon)\int_{-1}^{+1}\frac{\sgn(t)}{D^\epsilon(t)}dt\right),
$$
onde $\sgn(0)=0$ e $\sgn(x) = x/|x|$ para $x\ne0$.

Como $u^0$ satisfaz a mesma equação com $D^\epsilon$ substituído
por $D^*$, fazemos a mesma substituição em (\ref{eq:ueps}) para concluir
que $u^0$ será dado por
\be
  \label{eq:u00}
  u^0(x) = c_2(0) + c_1(0)
  \frac x{D^*} - \frac1{D^*} \int_0^x F(t) dt,
\ee
com
$$
  c_1(0) = \frac12 \int_{-1}^{+1}{F(t)}dt
$$
e
$$
  c_2(0) = \frac1{2D^*} \int_{-1}^{+1}\sgn(t){F(t)}dt.
$$

Valem as seguintes convergências uniformes em $[-1,1]$:
\be
  \label{eq:limsunif}
  \int_0^x \frac1{D^\epsilon(t)}dt \toe \frac x{D^*} \quad\mbox{e}\quad
  \int_0^x \frac{F(t)}{D^\epsilon(t)}dt \toe \frac 1{D^*}\int_0^x {F(t)}dt.
\ee
Para provar essas convergências, vamos utilizar uma propriedade geral de
funções periódicas, dada pelo lema abaixo. Resultados análogos valem sob
condições bem mais gerais, mas com a hipótese de suavidade a demonstração
se torna particularmente simples.
\begin{lema}[Propriedade do Valor Médio]
Seja $F:\R\times[a,b]\to\R$ de classe $C^1$, tal que $F(y+1,x)=F(y,x)$ para
todos $x\in[a,b],\,y\in\R$.
Então, dados $a',b'\in[a,b],\,a'<b'$, vale o seguinte limite
$$
  \lim_{\epsilon\downarrow0}\int_{a'}^{b'} F(x/\epsilon,x)dx = \int_{a'}^{b'} \bar F(x)dx
$$
uniforme em relação a $a'$ e $b'$, onde
$$
  \bar F(x) = \int_0^1 F(y,x)dy.
$$
\end{lema}
\prova
Considere $G:\R\times[a,b]\to\R$ dada por
$$
  G(y,x) = \int_0^y[F(z,x)-\bar F(x)]dz.
$$
Temos $G$ de classe $C^1$, periódica de período $1$ em $y$,
pois é a integral de uma função periódica da qual subtraímos a média;
portanto $G$ é limitada. Temos que
$$
  \partial_1 G(y,x) = F(y,x) - \bar F(x).
$$
Segue que
\begin{eqnarray*}
  \int_{a'}^{b'}[F(x/\epsilon,x)-\bar F(x)]dx
    &=& \int_{a'}^{b'} \partial_1 G(x/\epsilon,x)dx
  \\&=& \int_{a'/\epsilon}^{b'/\epsilon} \partial_1 G(y,\epsilon y) \epsilon dy
  \\&=& \epsilon \left[\int_{a'/\epsilon}^{b'/\epsilon}[\partial_1 G
       (y,\epsilon y) + \epsilon \partial_2 G(y,\epsilon y)]dy -
        \int_{a'/\epsilon}^{b'/\epsilon}  \partial_2 G(y,\epsilon y) \epsilon dy \right]
  \\&=& \epsilon \left[\int_{a'/\epsilon}^{b'/\epsilon}\Der{}y G(y,\epsilon y) dy -
        \int_{a'/\epsilon}^{b'/\epsilon}  \partial_2 G(y,\epsilon y) \epsilon dy \right]
  \\&=& \epsilon \left[ G(b'/\epsilon,b') - G(a'/\epsilon,a') -
        \int_{a'}^{b'}  \partial_2 G(x/\epsilon, x) dx \right].
\end{eqnarray*}
Aplicando a desigualdade triangular, temos
$$
  \left|\int_{a'}^{b'} F(x/\epsilon,x)dx - \int_{a'}^{b'} \bar F(x)dx\right| \leq
  \epsilon \left[\ 2\max_{[0,1]\times[a,b]}|G(y,x)| + (b-a)
  \max_{[0,1]\times[a,b]}|\partial_2 G(y,x)|\ \right],
$$
donde temos o limite uniforme desejado.
\eop

Pela definição de $D^*$ dada por (\ref{eq:d*}),
a Propriedade do Valor Médio tem como conseqüência imediata que vale o limite
$$
  \int_0^x \frac1{D^\epsilon(t)}g(t) dt \toe \frac1{D^*} \int_0^x g(t) dt
$$
uniforme em $[-1,1]$, onde $g$ é qualquer função em $C^1([-1,1])$.

Em particular, temos que valem os limites (\ref{eq:limsunif}), uniformemente
em $[-1,1]$. Também valem os limites $c_1(\epsilon)\toe c_1(0)$
e $c_2(\epsilon) \toe c_2(0)$, separando-se as integrais de $c_2(\epsilon)$
e $c_2(0)$ nos subintervalos $[-1,0)$ e $(0,1]$.

Comparando (\ref{eq:ueps}) e (\ref{eq:u00}), temos o limite (\ref{eq:lim1d})
uniforme em $[-1,1]$ como conseqüência dos limites mencionados acima.

\section{Equações diferenciais com coeficientes periódicos}
\label{sec:principal}

\subsection{Espaços de Sobolev}

Nesta seção $Q$ denotará um conjunto aberto limitado
em $\R^m$.
Os conjuntos de funções reais $L^\alpha(Q),\ 1\leq\alpha\leq+\oo$
são definidos da maneira usual.
O conjunto $C_0^\oo(Q)$ é formado pelas funções reais
infinitamente diferenciáveis com suporte compacto em $Q$.

\begin{defi}[Gradiente Fraco]
Uma função $v:Q\to\R^m$, dada por $v=(v_1,\dots,v_m)$, $v_i\in L^1(Q)$ é
chamada \emph{gradiente} da função $u\in L^1(Q)$ se
\be
  \label{eq:gradiente}
  \int_Q u(x) \dxi \phi(x)dx = - \int_Q v_i(x)
  \phi(x)dx,\quad \forall\phi\in C_0^\oo(Q),\ i=1,\dots,m.
\ee
Denotamos o gradiente por $\grad u$ ou $\dx u$.
\end{defi}
Observe que o gradiente,
quando existe, é único q.t.p. e que $\grad$ é linear, isto é, $\grad(u_1 +
\lambda u_2) = \grad u_1 + \lambda \grad u_2$, sempre que $ \grad u_1$ e $
\grad u_2$ existirem.
Por exemplo, se $u,v:(-1,1)\to\R$ são dadas por $u(x)=|x|$ q.t.p. e
$v(x)=x/|x|$ q.t.p. então $v$ é o gradiente de $u$ no sentido da
definição acima.
Observe também que se $u$ é diferenciável em algum aberto, então
as componentes do gradiente segundo a definição acima coincidem
com as respectivas derivadas parciais em quase todo ponto desse aberto.

Denotamos por $\LL^2(Q)$ o conjunto $\{p:Q\to\R^m \ :\ p_i\in L^2(Q),
i=1,\dots,m\}$ com
produto interno dado por
\be
  \label{eq:ll2inn}
  (\,p\,,\,q\,)_{\LL^2} = \int_Q p\cdot q\ dx
\ee e
norma $\|\,p\,\|_{\LL^2}^2 = {(p,p)}=
\int_Q |p|^2dx =
\|p_1\|_{L^2}^2 + \cdots + \|p_m\|_{L^2}^2$.

\begin{prop} $\LL^2(Q)$ é um espaço de Hilbert.
\end{prop}
\prova
A verificação de que (\ref{eq:ll2inn}) define um produto interno é direta.
Resta então mostrar que $\LL^2(Q)$ é completo.
É claro que $\|p_i\|_{L^2} \leq \|p\|_{\LL^2}$ para $i=1,\dots,m$ e para
todo $p\in\LL^2$. Seja $(p^n)$ uma seqüência de Cauchy em $\LL^2(Q)$.
Temos que $(p_i^n)$ são seqüências de Cauchy para $i=1,\dots,m$ e,
como $L^2(Q)$ é completo, temos que $p_i^n\ton^{L^2}p_i \in L^2(Q)$. Tomando
$p=(p_1,\dots,p_m)\in\LL^2(Q)$, temos que
$\|p^n-p\|_{\LL^2}^2=\|p_1^n-p_1\|_{L^2}^2+\cdots+\|p_m^n-p_m\|_{L^2}^2
\ton 0$, logo $p^n\ton^{\LL^2}p$. Portanto, $\LL^2(Q)$ é completo.
\eop

O \emph{espaço de Sobolev} $H^1(Q)$ é definido por
$$
  H^1(Q) = \left\{u\in L^2(Q): \mbox{existe } \grad u
  \mbox{ e este está em } \LL^2(Q)\right\}.
$$
Observe que faz sentido falar em gradiente de $u$
na definição acima, pois $L^2(Q)\cont L^1(Q)$ uma vez que
$Q$ é limitado.
É claro que $H^1(Q)$ é um espaço vetorial e, sobre  ele,
podemos definir a forma bilinear
\be
  \label{eq:inner}
  (u_1,u_2) = \int_Q (u_1u_2 + \grad u_1 \cdot \grad u_2)\ dx.
\ee
\'E f\'acil verificar que essa forma bilinear define um
produto escalar e que esse produto escalar gera a seguinte
norma
\be
  \|u\|_{H^1} = \sqrt{(u,u)} = \sqrt{\|u\|_0^2 + \|u\|_1^2},
\ee
onde $\|u\|_0^2 = \int_Q u^2 dx = \|u\|_{L^2}^2$ e
$\|u\|_1^2 = \int_Q |\grad u|^2 dx = \|\grad u\|_{\LL^2}^2$.
Mostraremos, a seguir, que
\begin{prop}
  $H^1(Q)$, com produto interno dado por (\ref{eq:inner}), é um espaço de Hilbert.
\end{prop}
\prova
\'E suficiente mostrar que $H^1(Q)$, com a norma $\|u\|_{H^1} = \sqrt{(u,u)}$, é
completo pois, como j\'a observamos,
(\ref{eq:inner}) define de fato um produto interno.
Seja $(u^n)_{n\in\N}$ uma seqüência de Cauchy em $H^1(Q)$.
Ponha $v^n = \grad u^n,\ n\in\N$.
Como $\|u\|_{L^2}\leq\|u\|_{H^1}$ e $\|\grad u\|_{\LL^2}\leq
\|u\|_{H^1}$ para toda $u\in H^1(Q)$,
temos que $(u^n)$ e $(v^n)$ são seqüências
de Cauchy em $L^2(Q)$ e $\LL^2(Q)$, respectivamente.
Mas esses espaços são completos,
logo $u^n\to u\in L^2(Q)$
e $v^n \to v\in\LL^2(Q)$ quando $n\to\infty$.
Basta então mostrar que
$v=\grad u$, pois neste caso teremos $u\in H^1(Q)$ e,
como
$\|u^n-u\|_{H^1}^2 = \|u^n-u\|_{L^2}^2 +
\|v^n - v \|_{\LL^2}^2 \ton 0$, teremos tamb\'em que
$u^n\ton^{H^1}u$.
Seja $\phi\in C_0^\oo(Q)$. Se $\phi=0$ temos (\ref{eq:gradiente})
trivialmente; assuma então que $\phi\ne0$.
Seja $\epsilon>0$. Faça $M=\|\phi\|_{H^1}>0$ e tome $n\in\N$ tal que
$\|u-u^n\|_{L^2} + \|v-v^n\|_{\LL^2} <\epsilon /M$.
Para cada $i=1,\dots,m$, temos por Cauchy-Schwartz
\begin{eqnarray*}
  \left| \int_Q(u-u^n)\dxi\phi dx + \int_Q(v_i - v_i^n) \phi dx \right|
  \leq \|u-u^n\|_{L^2} \|\grad \phi\|_{\LL^2} + \|v-v^n\|_{\LL^2}
  \|\phi\|_{L^2} \leq
  \\
  \leq (\|u-u^n\|_{L^2} + \|v-v^n\|_{\LL^2})M < \epsilon.
\end{eqnarray*}
Agora note que, pela definição de gradiente, vale
$$
  \int_Q u^n\dxi\phi dx + \int_Q  v_i^n \phi dx =0.
$$
Logo,
$$
  \left| \int_Q u\dxi\phi dx + \int_Q v_i  \phi dx \right| < \epsilon.
$$
Como $\epsilon>0$ é arbitrário, temos que $u$ e $v$ satisfazem
(\ref{eq:gradiente}) e, portanto, $v=\grad u$.
\eop

A desigualdade de Poincaré, que enunciamos abaixo\footnote{
A prova da Desigualdade de Poincaré pode ser encontrada em
\cite{evans98}.
},
será usada mais adiante no estudo da equação de  Poisson
com condição de contorno periódica.
\begin{teo}[Desigualdade de Poincaré]
Existe $c_0>0$, que depende apenas de $Q\cont\R^m$, tal que
$$
  \int_Qu^2dx \leq c_0\left[\left(
    \int_Qudx\right)^2 + \int_Q|\grad u|^2dx \right],\quad \forall\ u\in
    H^1(Q).
$$
\end{teo}

Um campo de vetores $v\in\LL^2(Q)$ é chamado de campo \emph{potencial}, ou
campo {\em gradiente}, se $v=\grad u$ para algum $u\in H^1(Q)$.

Dado $p\in\LL^2(Q)$, definimos o \emph{divergente} de
$p$ como sendo o funcional linear dado por
$$
  \div\ p \cdot \phi = (\div\ p,\phi) = -\int_Q(p\cdot \grad \phi)dx,
  \quad \phi\in H^1(Q).
$$

A desigualdade de
Cauchy-Schwartz em $\LL^2(Q)$ implica que $\div\ p$ é contínua, pois
$$
  (\div\ p, \phi)
  \leq \|p\|_{\LL^2} \|\grad\phi\|_{\LL^2}
  \leq \|p\|_{\LL^2} \|\phi\|_{H^1}.
$$
Convém observar que, para funções diferenciáveis em um aberto
de $Q$, as definições de divergente e campo potencial coincidem
com aquelas definidas no cálculo vetorial.

Outro exemplo de funcionais lineares e contínuos em
$H^1(Q)$ são aqueles dados por $u\mapsto\int_Qfudx$, onde
$f$ é um elemento de $L^2(Q)$.

\subsubsection{Espaços de funções periódicas}

Denotamos por $\square$ o paralelepípedo em $\R^m$, dado
por $\square = (0,l_1) \times \cdots \times (0,l_m)$.
Dizemos que uma funç\~ao definida em $\R^m$ \'e
peri\'odica, de período $\square$,  se ela for uma
função peri\'odica em
cada uma das variáveis $x_1,\dots,x_m$, e com períodos $l_1,\dots,l_m$,
respectivamente. Neste caso, tamb\'em dizemos que a função
\'e $\square-$peri\'odica.

Se $g(\cdot) : \R^m \to \R$ é uma função $\square-$periódica e
Lebesgue-integrável em $\square$, definimos o \emph{valor médio} de $g$ por
$$
  \langle g\rangle = \frac 1 {|\square|}\int_\square g(x)dx,
$$
onde $|\square| = l_1l_2\dots l_m$ é o volume do paralelepípedo $\square$.

A seguir, definiremos alguns espaços de funções que estarão relacionados,
de alguma maneira, com as funções de $\R^m$ em $\R$ que sejam
$\square-$periódicas. Assim, o símbolo $\square$ representa não apenas
o domínio dado pelo hipercubo mas também que estamos tratando de
funções periódicas. Em alguns casos, essa diferença vai além
do domínio de definição das funções do espaço em questão,
pois quando tratamos de espaços
que levam em consideração a derivada das funções, estaremos impondo
condições na fronteira deste domínio.

Por exemplo, dado $p\in[1,+\oo)$, denotamos por $L^p(\square)$ às
funções $g(\cdot):\R^m \to \R$ que sejam $\square-$periódicas e
tais que a norma
$$
  \|g\|_p = \langle |g|^p\rangle^{1/p}
$$
seja finita. É claro que $L^2(\square)$ com produto interno
$(u,v)=\av{uv}$ é um espaço de Hilbert. Neste caso, existe
uma bijeção direta entre $L^2(Q)$ e $L^2(\square)$ com
$Q=\square$, ou seja, os espaços são essencialmente os mesmos,
exceto porque as funções do primeiro espaço têm que
ser estendidas
de forma $\square$-periódica para obtermos elementos do
segundo.

Denotamos por $C^\oo(\square)$  o conjunto das funções
de $\R^m$ em $\R$,
$\square$-periódicas  e infinitamente diferenciáveis.
Definimos $H^1(\square)$ como sendo o completamento de
$C^\oo(\square)$ com respeito à norma $\|\cdot\|_{H^1}$.
Em $H^1(\square)$ definimos o produto interno $(u,v)=\av{uv+
\grad u \cdot \grad v}$ e norma $\|u\|_{H^1(\square)}=
\sqrt{(u,u)}$.

Observe que a diferença entre funções de $H^1(Q)$ e $H^1(\square)$ não
reside simplesmente no domínio em que estão definidas, mesmo quando
$Q=\square$.
De fato, se $\phi\in C^\oo(\R^m)$  então claramente $\phi|_{Q}\in C^\oo(Q)$.
Em particular, se $\phi\in C^\oo(\square)$,
vale $\phi|_{Q}\in C^\oo(Q)$. Entretanto, a recíproca
é falsa: se estendemos uma função $\phi\in C^\oo(Q)$
para todo o $\R^m$ de forma $\square$-periódica, em geral
não obtemos uma função de classe $C^\oo$,
pois não temos qualquer informação sequer sobre a continuidade de
$\phi$ na fronteira de $\square$ -- a função e suas derivadas
teriam que ``colar'' nesta fronteira. O ponto aqui é
que as funções de $H^1(Q)$ cujas extensões $\square$-periódicas estão
em $H^1(\square)$ são somente aquelas que podem ser arbitrariamente
aproximadas em $H^1(Q)$ por funções de $C^\oo(\square)$.
Podemos então pensar que $H^1(\square)$ é, em certo sentido,
um subespaço fechado de $H^1(Q)$.

De forma análoga ao que ocorre em $H^1(Q)$, temos que qualquer
$f\in L^2(\square)$ define um funcional linear
contínuo em $H^1(\square)$ dado por $(f,u)=\av {fu}\ \forall\
u\in H^1(\square)$. Da mesma forma,
se $g\in \LL^2(\square)$, então
o funcional $\div\ g$, dado por $(\div\ g,u)= -\av{g\cdot\grad u}$,
é linear e contínuo
em $H^1(\square)$. Uma conseqüência disso é que
 $u\mapsto\av u$  é um funcional
 contínuo em $H^1(\square)$ e $u\mapsto\av{\grad u}$  é uma
aplicação contínua de $H^1(\square)$ em $\R^m$.

Definimos também o espaço
$$
  \LL^2_\pot(\square) = \left\{f:\R^m\to\R^m :
  \exists~p\in\R^m,~u\in H^1(\square)~
   \text{tais que}~ f = p + \grad u  \right\}.
$$
Claramente $\LL^2_\pot(\square)$ é um espaço vetorial. Na verdade
este espaço é também de Hilbert, mas não vamos utilizar este fato.

\begin{prop}
  \label{prop:l2pot}
  Todo campo de vetores $v\in\LL^2_\pot(\square)$ é escrito
  como $v = \av{v} + \grad u,\ u\in H^1(\square)$.
\end{prop}
\prova
Primeiramente escrevemos ${v} = {p} + {\grad u},
\ p\in\R^m,u\in H^1(\square)$, donde $\av{v} = {p} + \av{\grad u}$.
O resultado desejado segue do fato
de que $\av{\grad u}=0\ \forall\ u\in H^1(\square)$.

Com efeito, se $u\in H^1(\square)$, temos que existe uma seqüência
$(\phi^n)$ em $C^\oo(\square)$ tal que $\phi^n\ton^{H^1}u$,
e, como $u\mapsto\av{\grad u}$ é contínua, temos
que $\av {\grad\phi^n} \to \av {\grad u}$.

Agora observe que, para qualquer função $\phi\in C^\infty(\square)$,
vale
$$
  \int_0^{l_i}\dxi \phi\ dx_i =
  \phi(x_1, \dots, x_{i-1},l_i,x_{i+1} \dots x_m)-
  \phi(x_1, \dots, x_{i-1},0,x_{i+1} \dots x_m) = 0
$$
e, pelo
Teorema de Fubini aplicado à integral múltipla em $\av {\dxi \phi}$,
temos que $\av {\dxi \phi}=0$.

Repetindo esse argumento para $i=1,\dots,m$, temos que $\av {\grad\phi} = 0$,
completando a prova. \eop

\subsection{Lema de Lax-Milgram}


Seja $V$ um espaço de Hilbert e $a(u,v):V\times V\to\R$
um funcional bilinear. Suponha que $a$ seja \emph{limitado},
no sentido de que
  \be
    \label{eq:acont}
    a(u,v) \leq \nu_2 \|u\|_V \|v\|_V
  \ee
para todo $u,v\in V$, onde $\nu_2>0$. Suponha também
  que $a$ seja \emph{coercivo}, isto é,
  \be
    \label{eq:acoer}
    a(u,u) \geq \nu_1 \|u\|_V^2,
  \ee
  para todo $u\in V$, onde $\nu_1>0$.

Seja $V^*$ o conjunto dos funcionais lineares e contínuos de $V$
em $\R$. Para cada $u\in V$, a forma linear
$v\mapsto a(u,v)$ é contínua  e, portanto, representa um elemento
de $V^*$, que denotaremos por $Au$. É claro que o operador
$A:V\to V^*$ que leva $u$ em $Au$ é linear e satisfaz
$$
  (Au,v) = a(u,v).
$$
onde $(\cdot,\cdot)$ denota o produto escalar em $V$.
Dada $f\in V^*$, considere a equação em $u$:
\be
  \label{eq:lax}
  a(u,v) = (f,v)\quad \forall\ v\in V.
\ee
isto \'e, dado $f\in V^*$, queremos resolver $Au=f$, com $u\in V$.

\begin{teo}[Lema de Lax-Milgram]
  O problema~(\ref{eq:lax}) tem uma única solução $u$, que sa\-tis\-faz
\be
  \label{eq:isomorf}
  \nu_1 \|u\|_V \leq  \|f\|_{V^*} \leq \nu_2 \|u\|_V.
\ee
Em outras palavras,
uma forma bilinear limitada e coerciva define um
isomorfismo  entre
$V$ e $V^*$.
\end{teo}
\prova
Seja $f\in AV\cont V^*$ e seja $u\in V$ tal que $Au=f$. Por
(\ref{eq:acoer}) temos
$$
  \nu_1 \|u\|_V^2 \leq (Au,u) = (f,u) \leq \|f\|_{V^*}\|u\|_V,
$$
donde segue a primeira desigualdade de (\ref{eq:isomorf}). Daí temos
que $A$ é injetivo, pois $Au=Av$ implica $A(u-v)=0$, logo $\|u-v\|_V \leq 0$.
Da injetividade de $A$ segue a unicidade
de soluç\~oes, quando as mesmas existirem.
Ainda supondo a exist\^encia de soluções,
a segunda desigualdade em (\ref{eq:isomorf}) \'e obtida de (\ref{eq:acont})
pois, dado qualquer $v\in V$, temos
$$
  |(f,v)| = |(Au,v)| = |a(u,v)| \leq \nu_2 \|u\|_V \|v\|_V.
$$

Resta mostrar que $AV=V^*$. Das desigualdades acima segue que
$A:V\to AV\cont V^*$ é um isomorfismo, e seqüências de Cauchy
em $AV$ são imagens de seqüências de Cauchy em $V$, que convergem
porque $V$ é completo. Portanto, $AV$ é fechado em $V^*$.

Agora seja $f\in AV^\perp$. Pelo Lema de Riez,
existe um único $z\in V$ tal que
$(f,v)=(z,v)\ \forall\ v\in V$. Temos que $(Au,z)=(Au,f)=0\ \forall\
u\in V$. Tomando $u=z$ temos $(Az,z)=0$ e, por (\ref{eq:acoer}), $z=0$,
logo $f=0$. Como $f\in AV^\perp$ é arbitrário, temos que $AV^\perp = \{0\}$.
Pelo Teorema da Projeção Ortogonal em Espaços de Hilbert,
$\overline{AV}=V^*$ mas, como $AV$ é fechado, segue que $AV = V^*$.
\eop

\subsection{O problema periódico}

Para tornar o texto mais limpo, vamos assumir, daqui em diante,
que estaremos tratando da soma de $1$ até $m$
sempre que os índices $i$, $j$ ou $k$ aparecerem repetidos,
a não ser que esteja explícito o contrário.

Escrevemos $\A(x) = \left(\,a_{ij}(x)\,\right)_{i,j=1}^m$ para denotar uma
matriz cujos elementos $a_{ij}:\R^m\to\R$ são funções mensuráveis e limitadas.
Fixado $x\in\R^m$ e dados $\xi,\eta\in\R^m$, definimos os produtos
$(\A\xi)_i = a_{ij}\xi_j$, $(\xi\A)_j =  \xi_ia_{ij}$.
Escrevemos também
$\eta \A \xi  = \eta(\A \xi)
= (\eta \A) \xi = \eta_i \A_{ij} \xi_j$.

Definimos
$$
  \|\A\| = \sup_{x\in\R^m} \sup_{\xi,\eta\in\R^m\setminus\{0\}}
  \frac{| \xi \A(x) \eta|}{|\xi||\eta|},
$$
de forma que
$$
  | \xi \A(x) \eta | \leq \|\A\|\ |\xi||\eta| \quad\forall\ x,\xi,\eta\in\R^m.
$$
Uma vez que os elementos de $\A$ satisfazem $|a_{ij}(x)|\leq M_{ij}$, temos que
$\|\A\|<\oo$ e $\sum_{i,j}M_{ij}$ é uma cota superior para $\|\A\|$.

Daqui em diante, vamos assumir que $\A$ sempre satisfaz a
\emph{condição de elipticidade}\label{ind:elipticidade}
$$
  \xi \A(x) \xi
  \geq \nu_1 |\xi|^2,\quad \forall\ x,\xi\in\R^m,
$$
para algum $\nu_1>0$. Dizemos também neste caso que
$\A$ é \emph{uniformemente elíptica}.

Considere a seguinte equação diferencial
$$
  \div (\A\grad u) = -f_0 + \div\ f,
$$
onde $f_0 \in L^2(\square)$ e $f\in\LL^2(\square)$.

Diremos que $u\in H^1(\square)$ é uma \emph{solução fraca} da equação acima
se $u$ satisfizer
\be
  \label{eq:periodica}
  \label{eq:periodica2}
  \av{\grad \phi \cdot \A \grad u} =  \av {f_0\phi} + \av{\grad\phi\cdot f},
  \quad \forall\ \phi\in H^1(\square).
\ee

Mostraremos agora que o problema acima possui solução e a
mesma é única no sentido de que duas soluções diferem por uma
função constante, desde que $\langle f_0\rangle=0$.

\begin{prop}
\label{prop:periodica}
O problema (\ref{eq:periodica}) possui solução $u\in H^1(\square)$
se, e somente se, $\av {f_0} = 0$.
Neste caso, tal solução é única no sentido de que $u-v$ é
constante para qualquer outra
solução $v$. Em particular, existe uma única solução satisfazendo $\av u=0$.
\end{prop}
\prova Suponha que existe uma solução $u\in H^1(\square)$.
Fazendo $\phi(x)=1\ \forall\
x\in\R^m$, temos $\grad\phi=0$ e $\phi\in H^1(\square)$. Substituindo
$\phi$ em (\ref{eq:periodica}), temos $\av {f_0} = 0$.

Reciprocamente, suponha que $\av {f_0} = 0$.
Se $u$ é uma solução então $u+c$ é também solução
para qualquer $c\in\R$, uma vez que somente $\grad u$ aparece em
(\ref{eq:periodica}). Assim, não é possível termos unicidade mais forte
do que a do enunciado acima. Por outro lado, suponha que $u,v$ são
duas soluções do problema e escreva $u = \av u + \tilde u$ e
$v = \av v + \tilde v$. Temos que $\tilde u$ e $\tilde v$ são também soluções
e satisfazem $\av{\tilde u}=\av{\tilde v}=0$. Como $u-v= \av{u-v} + (\tilde u
-\tilde v)$, para que $u-v$ seja constante basta mostrar que
$\tilde u = \tilde v$, ou seja, basta mostrar que existe uma única
solução $u$ satisfazendo $\av u=0$.

Para isso, defina $V = \left\{ u\in H^1(\square):\av u=0 \right\}$,
subespaço de $H^1(\square)$.
$V$ é fechado porque $u\mapsto\av u$ é contínua;
logo $V$ é um espaço de
Hilbert. Considere o funcional bilinear $a:V\times V\to\R$ dado
por $a(u,\phi)=\av{\grad\phi\cdot\A\grad u}$. Vale
$$
  a(u,u) = \av{\grad u\cdot\A\grad u} \geq \av{\nu_1|\grad u|^2}
  = \nu_1 \frac 1{|\square|} \int_\square |\grad u|^2dx,
$$
onde $\nu_1$ é a constante da condição de elipticidade de $\A$.
Da Desigualdade de Poincaré segue que
$$
  {c_0} \int_\square |\grad u|^2dx
  \geq  \int_\square u^2 dx,
$$
uma vez que $\av u=0$. Assim, temos
\begin{eqnarray*}
  a(u,u)
  &\geq& \frac {\nu_1}{(1+c_0)|\square|}  ( 1+ c_0)
  \int_\square |\grad u|^2dx
  \\ &\geq& \frac {\nu_1}{(1+c_0)|\square|} \left(
  \int_\square |\grad u|^2dx + \int_\square u^2 dx \right)
  \\ &=& \frac {\nu_1}{(1+c_0)} \|u\|_V^2,
\end{eqnarray*}
donde concluímos que $a$ é coercivo.

Sejam $u,v\in V$. Temos
$$
  |a(v,u)| = |\av{\grad u\cdot\A\grad v}| \leq
  \av{\,|\grad u\cdot\A\grad v|\,} \leq
  \|\A\| \av{ \, |\grad u|  \, |\grad v| \, },
$$
donde
$$
  |a(v,u)|^2 \leq  {\|\A\|^2 }
  \av{ |\grad u| |\grad v|} ^2
  \leq {\|\A\|^2 }
  \av { |\grad u|^2 } \av { |\grad v|^2 },
$$
logo $|a(v,u)| \leq \|\A\| \|u\|_V \|v\|_V$. Portanto, $a$ é limitado.

Considere o funcional linear $\phi\mapsto\av{f_0\phi}+\av{f\cdot\grad\phi}$,
definido em $V$. Esse funcional é contínuo, pois sua primeira parte
é dada por um elemento $f_0\in L^2(\square)$ e sua segunda parte
é o divergente de $-f\in\LL^2(\square)$. Pelo Lema de
Lax-Milgram, existe um único $u\in V$ tal que
$a(u,\phi)=\av{\grad\phi\cdot\A\grad u}=\av{f_0\phi}+\av{f\cdot\grad\phi}$
para toda $\phi\in V$. Agora seja $\phi\in H^1(\square)$. Temos
$\tilde \phi = \phi - \av \phi \in V$ e $\grad\tilde\phi=\grad\phi$.
Logo, $\av{\grad\phi\cdot\A\grad u}=\av{\grad\tilde\phi\cdot\A\grad u}
=\av{f_0\tilde\phi}+\av{f\cdot\grad\tilde\phi}
=\av{\,f_0(\phi-\av\phi)\,}+\av{f\cdot\grad\phi}
=\av{f_0\phi} - \av{f_0\av\phi}+\av{f\cdot\grad\phi}
=\av{f_0\phi}+\av{f\cdot\grad\phi}$, a última igualdade valendo porque
estamos assumindo que $\av{f_0}=0$. Portanto, (\ref{eq:periodica})
vale para toda $\phi\in H^1(\square)$.
\eop

\subsection{Um problema periódico auxiliar}

Vamos considerar o problema
\be
  \label{eq:auxiliar}
  \div(\A v) =0,\quad v\in\LL^2_\pot(\square),\ \ \av v=\lambda\in\R^m.
\ee
Todo elemento de $\LL^2_\pot(\square)$ é escrito
como $v = \lambda + \grad u,\ u\in H^1(\square)$,
sendo $\lambda=\av v$.
Assim,
obtemos o problema abaixo para $u$:
$$
  \div(\A(\lambda+\grad u))=0,\quad u\in H^1(\square).
$$
Podemos reescrever o problema acima na forma (\ref{eq:periodica2})
com $f=-\A\lambda$ e $f_0=0$. Pela Proposição~\ref{prop:periodica},
existe uma única solução $u$ do problema acima, a menos de
adição de uma constante. Como $v$ depende apenas de $\grad u$,
adição de constante em $u$ não altera a solução $v$, logo $v$ é única.
Em outras palavras, para cada $\lambda\in\R^m$, o problema
(\ref{eq:auxiliar}) possui uma única solução $v\in\LL^2_\pot(\square)$.

Por outro lado, $v\mapsto\div(\A v)$ é uma aplicação linear,
logo seu núcleo é um subespaço vetorial.
Assim, se $v$ e $w$ são duas soluções arbitrárias, isto é,
$\div(\A v)=\div(\A w)=0$, e $c\in\R$ é um escalar arbitrário, temos
$\div(\A(v+cw))=0$ e $\av{v+cw} = \av u + c \av w = \lambda + c\xi$,
onde $\lambda = \av v$ e $\xi = \av w$.
Devido à unicidade da solução de (\ref{eq:auxiliar}) para cada
parâmetro $\lambda$ fixo, podemos tomar o caminho inverso, ou seja,
dados $\lambda,\xi\in\R^m$ e $c\in\R$ arbitrários, temos que
a solução do problema (\ref{eq:auxiliar}) com parâmetro $\lambda + c\xi$
é $v+cw$. Assim, a solução $u$ depende linearmente de $\lambda$.

Portanto, $\av{\A v}$ depende linearmente de $\lambda$ e, conseqüentemente,
pode ser representada por
\be
  \label{eq:a0}
  \av{\A v} = \A^0\lambda,
\ee
onde $\A^0$ é uma matriz real com elementos constantes.

O problema periódico auxiliar e a determinação dos elementos da matriz
$\A^0$ em termos de $\A(x)$ é de importância central
no estudo dos problemas homogeneizados  -- veja as
seções subseqüentes.

Vamos então tentar encontrar os elementos $a_{ij}^0$ de
$\A^0 = (a_{ij}^0)_{i,j=1}^m$. Seja $e^1,\dots,e^m$
a base canônica de $\R^m$. Para cada $k=1,\dots,m$,
a solução de (\ref{eq:auxiliar})
com $\lambda = e^k$ é da forma $v = e^k + \grad N_k,\
\ N_k\in H^1(\square)$. Temos
$$
  a_{ij}^0 = e^i\cdot \A^0 e^j = \av {e^i \cdot \A v} =
  \av { e^i\cdot \A(e^j + \grad N_j) } = \av{a_{ij} + e^i\cdot\A\grad
  N_j},
$$
logo
\be
  \label{eq:aij0}
  a_{ij}^0 = \av{a_{ij}} + \av{a_{ik}\der{N_j}{x_k}}.
\ee
Por definição, $N_k$ é a solução de
$$
  \div(\A(e^k+\grad N_k))=0,\quad N_k\in H^1(\square),
$$
ou, no caso de $\A$ suave:
\be
  \label{eq:nk}
  \dxi\left(a_{ij}\der{N_k}{x_j}\right)=-\dxi(a_{ik}).
\ee

\begin{prop}
A matriz $\A_0$ é elíptica. Além disso, se $\A(x)$ é simétrica para
todo $x$, então $\A^0$ é simétrica.
\end{prop}
\prova
Segue de (\ref{eq:auxiliar}) que
$$
  \grad \phi \perp \A v\quad \forall\phi\in H^1(\square).
$$
Em particular, $\A v \perp (v-\lambda)$, logo $\av {v\A v}=\av{\lambda\A v}$.
Disso segue que
$$
  \lambda\cdot\A^0\lambda = \lambda\cdot\av{\A v} = \av{\lambda\cdot\A v} =
  \av {v\cdot\A v} \geq \nu_1 \av{|v|^2} \geq \nu_1 |\av v|^2 = \nu_1 |\lambda|^2,
$$
ou seja, $\A^0$ é elíptica.

Para provar a segunda parte, consideramos o problema análogo
\be
  \label{eq:auxiliar2}
  \div(w \A ) =0,\quad w\in\LL^2_\pot(\square),\ \ \av w=\xi\in\R^m,
\ee
cuja única solução $w$ depende linearmente de $\xi$, e definimos a matriz
constante $\CC^0$ por
\be
  \label{eq:c0}
  \av{w\A} = \xi\CC^0.
\ee
Do mesmo argumento de ortogonalidade, segue que
$$
  \A v \perp (w-\xi)\quad\mbox{e}\quad w\A\perp(v-\lambda),
$$
e, conseqüentemente,
$$
  \xi\A^0\lambda = \xi\av{\A v} =  \av{\xi\A v} =
  \av{w\A v} = \av{w\A \lambda} = \av{w\A}\lambda= \xi \CC^0 \lambda.
$$
A igualdade acima vale para todos $\lambda,\xi\in\R^m$.
Segue que $\A^0 = \CC^0$. Se definimos $\CC^0$
utilizando $\A^T$ ao invés de $\A$, obtemos
$(\A^0)^T$, implicando na seguinte propriedade: $(\A^0)^T=(\A^T)^0$.
Em particular, quando $\A(x)$ é simétrica para todo $x$,
temos $\A^0$ simétrica.
\eop

\subsection{Homogeneização de equações parabólicas}

Nesta seção oferecemos uma prova construtiva para o
Teorema~\ref{teo:parab} utilizando o método conhecido como
\emph{método de expansão assintótica}.

Este método é comum no estudo teórico de sistemas físicos e
consiste em considerar formalmente a seguinte
aproximação para $u^\epsilon$:
$$
  u^\epsilon(x,t) \approx u^0(x,t) + \epsilon u_1(x,y,t)
  + \epsilon^2 u_2(x,y,t) + \epsilon^3 u_3(x,y,t) + \cdots,
$$
onde $y=\epsilon^{-1}x$, as funções $u_1,u_2,u_3,\dots$ são periódicas em $y$
e não dependem de $\epsilon$ e a função $u_0$ depende apenas de $x$ e $t$.

A idéia do método é fazer com que o truncamento da expansão
acima até o $n$-ésimo termo, que denotaremos por $u_n^\epsilon$,
por um lado satisfaça ``aproximadamente'' a equação diferencial
que estamos considerando e por outro lado seja bem aproximado
pelo seu primeiro termo $u_0(x,t)$, que não depende de $\epsilon$ direta-
nem indiretamente.

Na prova que apresentamos abaixo vamos considerar a aproximação
até segunda ordem em $\epsilon$, dada por
\be
  \label{eq:u2e}
  u_2^\epsilon(x,t) = u^0(x,t) + \epsilon u_1(x,y,t)
  + \epsilon^2 u_2(x,y,t),\quad y=\epsilon^{-1}x.
\ee
Apresentamos abaixo o teorema mencionado e seu corolário.

\begin{teo}
\label{teo:parab}
 Para cada $\epsilon>0$, seja $u^\epsilon$ a solução (clássica)
 do seguinte problema
\be
  \label{eq:parab}
  \left\{\ba{l} \displaystyle
   \rho^\epsilon u^\epsilon_t - \div(\A^\epsilon\grad u^\epsilon)
  = 0,\quad x\in\R^m,\,  t>0,
  \\ \displaystyle
  u^\epsilon(x,0)=\phi(x),\quad\phi\in C_0^\oo(\R^m)
  \ea \right.
\ee
onde $\rho^\epsilon = \rho(\epsilon^{-1}x)$ e $\rho(x)$ é
uma função $C^1$, positiva e periódica que satisfaz $\av \rho = 1$;
$\A^\epsilon = \A(\epsilon^{-1}x)$ onde $\A(x)$ é uma matriz $m\times m$
periódica (de mesmo período de $\rho$), uniformemente elíptica,
cujos elementos são funções $C^1$ de $x$. Então
\be
  \label{eq:convhomog}
  \sup_{x\in\R^m} | u^\epsilon(x,t) - u^0(x,t)| \toe 0, \quad
  \forall\ t \geq 0,
\ee
onde $u_0$ é a solução de
\be
  \label{eq:parabh}
  u^0_t -
  \div(\A^0\grad u^0) = 0,\quad
  u^0(\cdot,0)=\phi
\ee
e $\A^0$ é a matriz constante dada por (\ref{eq:a0}).
\end{teo}
\prova
A prova consiste em estudar a aproximação até segunda ordem para $u^\epsilon$,
dada por (\ref{eq:u2e}).

Primeiramente, definimos o operador
$$
  \Der{}{x_i} = \der{}{x_i} + \epsilon^{-1}\der{}{y_i},\quad i=1,\dots,m,
$$
e observamos que o operador $\div(\A^\epsilon\grad)$ pode
ser reescrito como
$$
  \div(\A^\epsilon\grad) = \Der{}{x_i}
  \left( a_{ij}(\epsilon^{-1}x) \Der{}{x_j} \right).
$$
 Desenvolvendo o lado direito da igualdade acima, obtemos
$$
  \div(\A^\epsilon\grad) = \epsilon^{-2}A_1 + \epsilon^{-1}A_2 + A_3,
$$
onde
\begin{eqnarray*}
   A_1 &=& \der{}{y_i}\left(a_{ij}(y)\der{}{y_j}\right),\\
   A_2 &=& \der{}{y_i}\left(a_{ij}(y)\der{}{x_j}\right) +
           a_{ij}(y)\der{^2}{x_i\partial y_j},\\
   A_3 &=& a_{ij}(y) {\partial^2\over\partial x_i\partial x_j}.
\end{eqnarray*}

Definindo o operador $L_\epsilon$ como
$$
  L_\epsilon  = \rho^\epsilon\der{}t -
  \div(\A^\epsilon\grad ),
$$
nós olhamos para a sua ação sobre $u_2^\epsilon$:
\begin{eqnarray}
  -L_\epsilon u_2^\epsilon
  &=& \left( -\rho(t)\der{}t +
  \epsilon^{-2}A_1 + \epsilon^{-1}A_2 + A_3 \right)
  (u^0 + \epsilon u_1 + \epsilon^2 u_2)
  \nonumber
  \\ &=&
  \epsilon^{-1}(A_1u_1+A_2u^0)+\left(-\rho(y)\der{u^0}t+
  A_1u_2 + A_3u^0 + A_2u_1 \right) +
  \nonumber
  \\ &&
  \label{eq:leu2e}
  +\ \epsilon\left( A_3u_1 + A_2u_2 - \rho(y)\der{u_1}t - \epsilon\rho(y)
  \der{u_2}t + \epsilon A_3u_2\right).
\end{eqnarray}

Queremos determinar condições suficientes em $u^0$,  $u_1$ e
$u_2$ para que $L_\epsilon u_2^\epsilon$ seja arbitrariamente
pequeno.

Começamos igualando o termo que multiplica $\epsilon^{-1}$ a zero. Obtemos
$A_1u_1+A_2u^0=0$, ou seja
$$
  \der{}{y_i}\left(a_{ij}(y)\der{}{y_j}u_1(x,y,t)\right)
  = - \der{}{y_i}\left(a_{ij}(y)\der{u^0(x,t)}{x_j}\right).
$$
Dado $u^0(x,t)$, vamos
olhar para a igualdade acima como  uma equação elíptica com incógnita
$u_1$ e variável independente $y$. Portanto, o problema acima é periódico
em $y$, enquanto $x$ e $t$ são considerados como parâmetros.
Pondo $N_k(y),\,k=1,\dots,m$ como
as soluções de (\ref{eq:nk}), podemos tomar
\be
  \label{eq:u1}
  u_1(x,y,t) = N_k(y)\der{u^0(x,t)}{x_k}
\ee
como solução da equação acima.

Agora vamos igualar a zero o termo que multiplica $\epsilon^0$,
obtendo assim um problema periódico para a incógnita $u_2$
e com variável independente $y$, dado por
\be
  \label{eq:u2}
  \div_y(\grad_y u_2(x,y,t)) = -\left( -\rho(y)\der{u^0(x,t)}t +
  A_3(y)u_0(x,t)+A_2(y)u_1(x,y,t)\right).
\ee
Pela Proposição~\ref{prop:periodica}, a existência de uma solução
$u_2$ para a equação acima depende da condição
\be
  \label{eq:avu2}
  \av{ -\rho(y)\der{u^0}t + A_3u_0+A_2u_1} = 0.
\ee
Mas
$$
  A_3u^0 = a_{ij}(y){\partial^2u^0(x,t)\over\partial x_i\partial x_j}
$$
e
\begin{eqnarray*}
  \av{A_2u_1}
  &=& \av{\der{}{y_i}\left(a_{ij}\der{u_1}{x_j}\right)} +
      \av{a_{ij}\der{^2u_1}{x_i\partial y_j}}
  \\&=& 0 + \av{a_{ij}\der{^2u_1}{x_i\partial y_j}}
  \\&=& \av{a_{ij}\der{^3N_ku^0}{x_i\partial x_k\partial y_j}}
  \\&=& \av{a_{ik}\der{N_j}{y_k}}\der{^2u^0}{x_i\partial x_j},
\end{eqnarray*}
onde na primeira igualdade substituímos a definição de $A_2$;
na segunda, aplicamos o mesmo resultado contido na
prova da Proposição~\ref{prop:l2pot}, de que as derivadas
de funções periódicas têm média zero; na terceira,
substituímos (\ref{eq:u1}) e, na quarta, apenas trocamos
$j$ por $k$.

Substituindo em (\ref{eq:avu2}), temos
$$
  \der{u^0}t =
  \av{ \rho \der{u^0}t} = \av{A_3u_0+A_2u_1} =
  \av{ a_{ij} + a_{ik}\der{N_j}{y_k}}\der{^2u^0}{x_i\partial x_j}
  = a_{ij}^0 \der{^2u^0}{x_i\partial x_j},
$$
onde a última igualdade vem de (\ref{eq:aij0}). Portanto, $u^0$
deve satisfazer
$$
  \der{u^0}t - \div(\A^0\grad u^0) = 0.
$$
Ponha então $u^0$ como solução da equação acima com condição inicial
$u(\cdot,0)=\phi$. Com $u^0$, $u_1$ e $u_2$ dados acima,
acabamos de construir a aproximação (\ref{eq:u2e}).

Usamos a desigualdade triangular para obter
$$
  \dsup_{\R^m\times[0,T]}|u^\epsilon-u^0|
  \leq
  \dsup_{\R^m\times[0,T]}|u^\epsilon-u_2^\epsilon|
  + \dsup_{\R^m\times[0,T]}|u_2^\epsilon-u^0|,
$$
e vamos tentar controlar cada termo separadamente.

O primeiro termo depende essencialmente de $|L_\epsilon u_2^\epsilon(x,t)|$
e de $|u^\epsilon(x,0)-u_2^\epsilon(x,0)|$.
Isso vem do Corolário~\ref{cor:maximo} do Apêndice~\ref{sec:maximo}
aplicado à função $u^\epsilon-u_2^\epsilon$ e do fato de que $u^\epsilon$,
sendo solução, satisfaz $L_\epsilon u^\epsilon =0$:
\be
  \label{eq:maximo}
  \dsup_{\R^m\times[0,T]}|u^\epsilon-u_2^\epsilon| \leq
  \dsup_{\R^m}|u^\epsilon(x,0)-u_2^\epsilon(x,0)| +
  (\min\rho)^{-1} T \dsup_{\R^m\times[0,T]}|L_\epsilon u_2^\epsilon(x,t)|.
\ee

O teorema está provado se mostramos que, dado $T>0$, valem
os limites abaixo, quando $\epsilon\to0$:
\begin{eqnarray}
&  \dsup_{\R^m\times[0,T]}|u_2^\epsilon-u^0|\to 0
\label{eq:limite1}\\
&  \dsup_{\R^m}|u^\epsilon(x,0)-u_2^\epsilon(x,0)|\to 0
\label{eq:limite2}\\
&  \dsup_{\R^m\times[0,T]}|L_\epsilon u_2^\epsilon| =
\label{eq:limite3}  \dsup_{\R^m\times[0,T]}|L_\epsilon (u_2^\epsilon-u^\epsilon)| \to  0.
\end{eqnarray}

Comecemos por (\ref{eq:limite1}). É claro que $u^0$ e suas derivadas são
limitados em $\R^m\times[0,T]$. Considerando que as funções $N_k(y)$ são
periódicas e de classe $C^1$, portanto limitadas, temos que $u_1$, dado por
(\ref{eq:u1}), e $u_2$, dado pela única solução de (\ref{eq:u2})
satisfazendo $\av{u_2}=0$, são limitados em $\R^m\times[0,T]$. Segue então
que $u_2^\epsilon$, dada por (\ref{eq:u2e}), converge uniformemente em
$\R^m\times[0,T]$ para $u^0$ quando $\epsilon\to0$. O limite
(\ref{eq:limite2}) é uma conseqüência direta de (\ref{eq:limite1}) pois
$u^\epsilon(x,0)=u^0(x,0)=\phi(x)$. Para obtermos (\ref{eq:limite3}),
primeiro note que $L_\epsilon u_2^\epsilon$ é dado apenas pelo termo
(\ref{eq:leu2e}), uma vez que $u^0$, $u_1$ e $u_2$ foram definidos de forma
a anular os demais termos. Mas da mesma forma que $u_1$ e $u_2$ são
limitados em $\R^m\times[0,T]$, também o são suas derivadas com respeito a
$x$, $y$ e $t$, de forma que $L_\epsilon u_2^\epsilon$ converge
uniformemente para zero em $\R^m\times[0,T]$.
 Observe que estamos nos baseando no fato de que a solução $u_2$ da
equação elíptica~(\ref{eq:u2}), assim como suas primeiras derivadas, são
uniformemente limitadas, uma vez que $u_0$, sendo solução da equação do
calor, também tem suas derivadas limitadas e $u_1$ é dado
por~(\ref{eq:u1}).

Os limites acima, juntamente com a estimativa (\ref{eq:maximo}),
resultam em
$$
  \dlim_{\epsilon\to 0} \dsup_{\R^m\times[0,T]}|u^\epsilon-u^0| = 0.
$$
Agora note que (\ref{eq:convhomog}) é conseqüência direta deste último
limite e a prova está concluída.
\eop

O corolário abaixo justifica o argumento formal utilizado no
Capítulo~\ref{chap:rgana}.
\begin{cor}
\label{cor:parab}
Para cada $\epsilon>0$, seja $u^\epsilon$ a
solução clássica do seguinte problema
$$
  \left\{\ba{l}
  u^\epsilon_t = 
  [1+\mu g(\epsilon^{-1}    x)] u^\epsilon_{xx}, \quad x\in\R,\,t\geq0,
  \\
  u^\epsilon(x,0)=\phi(x),\quad\phi\in C_0^\oo(\R^m)
  \ea \right.,
$$
onde $g$ é suave e periódica e $\mu$ satisfaz
$[1+\mu g(y)]>0\ \forall\ y\in\R$. Seja $u^0$ solução de
\be
  \label{eq:parabhcor}
  u^0_t = \sigma u^0_{xx},\quad u^0(\cdot,t=0)=\phi,
\ee
onde $\sigma = \av {[1+\mu g(y)]^{-1}}^{-1}$.
Então
\be
  \label{eq:convcor}
  \sup_{x\in\R^m} | u^\epsilon(x,t) - u^0(x,t)| \toe 0, \quad
  \forall\ t \geq 0.
\ee
\end{cor}
\prova
Podemos reescrever a equação para $u^\epsilon$ na forma
(\ref{eq:parab}) com
$$
  \rho(y) = [1+\mu g(y)]^{-1}
$$
e $\A(x)$ constante igual a $1$.
Pelo Teorema~\ref{teo:parab}, vale o limite (\ref{eq:convcor})
sendo $u^0$ a solução de
$$
  \av \rho u^0_t = u^0_{xx},\quad u^0(\cdot,t=0)=\phi.
$$
Mas, pela definição de $\sigma$, temos $\av \rho = \sigma^{-1}$,
portanto $u^0$ é também a solução de (\ref{eq:parabhcor}) e
não há mais nada a provar. \eop

%% file: apend.tex
\chapter{Apêndices}

\section{A transformada de Fourier}
\label{sec:tf}

Neste ap\^endice vamos definir a transformada de Fourier e enunciar as suas
propriedades que são utilizadas no Capítulo~\ref{chap:rgana}. As
demonstrações dos resultados abaixo podem ser encontrados no
livro de Iório \cite{iorio}. Representamos por $dx$ a medida de Lebesgue na reta.

\begin{defi}[Transformada de Fourier]
Seja $f$ um elemento de $L^1(\R)$. Definimos a transformada de
Fourier de $f$  através da fórmula \be \label{eq:tf}
  \F\{f(x)\}  = \int_\R f(x)e^{-ikx}dx = \hat f(k).
\ee
\end{defi}

$\F\{\cdot\}$ atua linearmente em $L^1(\R)$ e além disso, para
qualquer elemento $f\in L^1(\R)$ e $a\in\R$,
\be
  \label{eq:proptf}
  \F \{f(ax)\} = \frac 1 {|a|} \widehat f\left({k} / {a}
  \right).
\ee
Outra propriedade é que se $f\in L^1(\R)$, então $\hat f$ é contínua
e limitada.

Para funções dentro da classe de Schwartz ${\cal S}
(\R)$\footnote { ${\cal S}(\R) =
\{f:\mathbb{R}\rightarrow\mathbb{C}/ |x^mD^nf(x)|\leq M(m,n)\,\,
\forall \,\, m,n\in\mathbb{N}\}$, onde
 $D^nf(x)$
denota a $n$-ésima derivada de $f$ no ponto $x$.}, vale o
seguinte resultado:
\begin{teo} Se $f\in{\cal S}$ então $\hat{f}\in{\cal S}$ e
\be
  \widehat{ f'}(k) = -i k \hat f(k).
\ee Além disso, para quaisquer $f,g \in {\cal S}$: \be
  (\hat f,\hat g) = 2 \pi  ( f,  g),
\ee onde $(\cdot,\cdot)$ denota o produto interno de
$L^2(\R)$. Em particular, para todo $f\in {\cal S}$,
\be
  \label{eq:tfisom}
  \|\hat f\|_{L^2} = \sqrt{ 2 \pi } \,  \| f\|_{L^2}.
\ee
\end{teo}
\obs Como ${\cal S}$ é um
subespaço denso em $L^2(\mathbb{R})$, a definição de $\F$,
pode ser estendida ao espaço de Hilbert  $L^2(\mathbb{R})$
ainda satisfazendo a relação (\ref{eq:tfisom}). A
transformada inversa de Fourier está bem definida
em $L^2(\mathbb{R})$ e a mesma é dada pela fórmula de inversão
abaixo.
\begin{teo}[Transformada inversa de Fourier]
Para todo $f\in L^2(\mathbb{R})\cap L^1(\mathbb{R})$, temos que $\hat{f}\in
L^2(\mathbb{R})$ e, além disso: \be
  f(x) = \frac 1 {2 \pi} \lim_{R\to\oo} \int_{-R}^R  \hat{f}(k)e^{ikx}dk,
\ee
o limite valendo em $L^2(\R)$. Para uma subseqüência $R_n\to\oo$ o
limite também vale q.t.p.
\end{teo}
\begin{defi}[Convolução]Sejam $f,g\in L^1(\mathbb{R})$. Definimos
a convolução de $f$ com $g$ pela integral \be
  (f*g)(x) = \int_\R f(x-y)g(y)dy.
\ee
\end{defi}
A convolução de duas funções $f$ e $g$ satisfaz \`as
seguintes propriedades: $f*g = g*f$  e   $(f*g)*h = f*(g*h)$. O
teorema da convolução nos diz que:
\begin{teo}[Teorema de convolução]
\label{teo:conv} Para todo $f,g \in {\cal S}(\mathbb{R})$: \be
  \F\{ f*g \}(k) = \hat f(k) \hat g(k),
\ee ou seja, \be
  \F^{-1} \{ \hat f \hat g \}(x) = (f*g)(x).
\ee
\end{teo}

Para finalizar, lembremos que a distribuição gaussiana se
transforma da seguinte maneira:
\begin{equation}
  \label{eq:tran-gauss}
  \F \left\{e^{\frac{-x^2} {4\sigma}} \right\}
  = \sqrt{4 \sigma \pi}\, e^{-\sigma k^2},
\end{equation}
onde $\sigma$ é uma constante estritamente positiva.

\section{O princípio do máximo}
\label{sec:maximo}

Neste apêndice vamos estudar propriedades sobre o máximo de funções
que satisfazem uma equação (ou melhor, uma inequação) diferencial
dada por um operador parabólico, que consiste em
tomarmos o divergente do produto de uma matriz
uniformemente elíptica\footnote{Veja a \emph{condição de elipticidade}
definida na página~\pageref{ind:elipticidade}.
} pelo gradiente da função, subtraído de sua derivada no tempo.
Trata-se de uma das várias formulações do
conhecido \emph{Princípio do Máximo} para equações diferenciais
parabólicas.

Vamos enunciar o teorema abaixo sem no entanto demonstrá-lo\footnote{
Uma versão mais geral desse teorema
pode ser encontrada em \cite{friedman83},~p.~43,~\textsl{Theorem 9}.
}.

\begin{teo}
\label{teo:maximo}
Seja $u:\R^m\times[0,\oo)\to\R$ de classe $C^2$ satisfazendo
a seguinte desigualdade
$$
  L[u] = \left( a_{ij}(x,t)\der{^2u}{x_ix_j}\right)
  + \left(b_i(x,t)\der u{x_i}\right) - \der ut \leq 0
$$
em $E_T = \R^m\times(0,T]$.
Suponha que os coeficientes de $L$ são contínuos e limitados
e que a matriz $(a_{ij})$ é uniformemente elíptica.

Nessas condições, se a função $u$ é limitada em $\overline{ E_T}$
e $u(x,0)\geq0\ \forall\ x\in\R^m$ então $u(x,t)\geq 0\ \forall\ (x,t)
\in E_T$.
\end{teo}

\begin{cor}
Seja $\rho:\R^m\to\R$ limitada de classe $C^1$ tal que $\inf_{\R^m}\rho(x)>0$.
Então o Teorema~\ref{teo:maximo} continua válido se definimos
\be
  \label{eq:apoprho}
  L[u] = \der{}{x_i}\left(a_{ij}(x,t)\der u{x_j}\right)
  - \rho(x) \der ut,
\ee
desde que os elementos da matriz $(a_{ij})$ sejam de classe $C^1$
e suas derivadas sejam limitadas.
\end{cor}
\prova
Desenvolvendo a primeira parte de $L$ obtemos
$$
  L[u] = \left( a_{ij}(x,t)\der{^2u}{x_ix_j}\right) +
  \left(\der{a_{ij}(x,t)}{x_i}\der u{x_j}\right)
  - \rho(x) \der ut.
$$
Definimos $\tilde L[u](x,t)=\rho(x)^{-1}L[u](x,t)$,
cujos coeficientes ainda satisfazem as hipóteses do
Teorema~\ref{teo:parab}. O corolário segue então do fato de que
$\tilde L[u]\leq 0 \Leftrightarrow L[u]\leq 0$.
\eop

\begin{cor}
Seja $u:\R^m\times[0,\oo)\to\R$ uma
função  de classe $C^2$ tal que
$L[u]\geq 0$ em $E_T=\R^m\times(0,T]$ e $u$ é limitada em $\overline{E_T}$,
onde $L[\cdot]$ é definido por
(\ref{eq:apoprho}). Então
$$
  \sup_{\R^m\times[0,T]} u(x,t) \leq \sup_{\R^m} u(x,0).
$$
\end{cor}
\prova
Defina $v(x,t) = [\sup_{\R^m} u(x,0)] - u(x,t)$. Claramente
$v$ é de classe $C^2$ e satisfaz $v(x,0)\geq 0$ para
todo $x\in\R^m$. Além disso, $L[v] = -L[u] \leq 0$ e $v$ é limitada
em $\overline{E_T}$. Logo, $v(x,t)\geq 0$ em $E_T$.
Substituindo, obtemos que $u(x,t) \leq [\sup_{\R^m} u(x,0)]$
para todo par $(x,t)\in E_T$.
\eop

\begin{cor}
\label{cor:maximo}
Seja $v:\R^m\times[0,\oo)\to\R$ duas vezes diferenciável
e considere $L[\cdot]$ definido por
(\ref{eq:apoprho}). Então
$$
  \sup_{\R^m\times[0,T]} |v(x,t)| \leq \sup_{\R^m} |v(x,0)| + (\inf_{\R^m}\rho)^{-1}T
  \sup_{\R^m\times[0,T]} |L\,v(x,t)|.
$$
\end{cor}
\prova
Defina
$$
  u(x,t) = - \left((\inf\rho)^{-1} \sup_{\R^m\times[0,T]} |L\,v(x,t)| \right)t + v.
$$
Segue que
\begin{eqnarray*}
  Lu &=& L\left[- \left((\inf\rho)^{-1}
            \sup_{\R^m\times[0,T]} |L\,v(x,t)| \right)t\right] + Lv
  \\  &=& \frac\rho{\inf\rho} \sup_{\R^m\times[0,T]} |L\,v(x,t)|  + Lv
  \\ &\geq& \frac \rho{\inf\rho} \sup_{\R^m\times[0,T]} |L\,v(x,t)| - |Lv|
  \\ &\geq& \sup_{\R^m\times[0,T]} |L\,v(x,t)| - |Lv|
  \\ &\geq& 0.
\end{eqnarray*}
Do corolário anterior segue que
$$
  v(x,t) \leq \sup_{\R^m} |v(x,0)| +
  \left((\inf\rho)^{-1} \sup_{\R^m\times[0,T]} |L\,v(x,t)| \right) t
$$
para todo $t\in[0,T]$.
Aplicamos para $-v$ o resultado que acabamos de obter e concluímos que
$$
  |v(x,t)| \leq \sup_{\R^m} |v(x,0)| +
  \left((\inf\rho)^{-1} \sup_{\R^m\times[0,T]} |L\,v(x,t)| \right) t
$$
para todo $t\in[0,T]$. A estimativa desejada segue direto desta última.
\eop

%% file: errata.tex
\noindent {\Huge \bf Errata}
\addcontentsline{toc}{chapter}{Errata}

\vspace{1cm}

Esta errata contém os erros encontrados até 08/05/2006.

Pede-se comunicar quaisquer outros erros matemáticos, tipográficos
ou gramaticais através do endereço ele\-trô\-ni\-co leorolla@impa.br.

\begin{itemize}

\item Página~\pageref{eq:lglgrel}.
Na equação~(\ref{eq:lglgrel}) deve-se ler
$$
  \ln A_n \approx \ln \left[\frac A {(\ln L)^\alpha}\right] - \alpha \ln n.
$$

\end{itemize}

%% file: teseleo.bbl
\begin{thebibliography}{10}

\bibitem{ames92}
W.~F. Ames.
\newblock {\em Numerical Methods for Partial Differential Equations}.
\newblock Academic Press, 3$^{\rm a}$~edi\c c\~ao, 1992.

\bibitem{papa78}
A.~Bensousan, J.~Lions, G.~Papanicolaou.
\newblock {\em Asymptotic Analysis of Periodic Structure}.
\newblock North Holland, Amsterdam, 1978.

\bibitem{braga04}
G.~A. Braga, F.~Furtado, J.~M. Moreira, L.~T. Rolla.
\newblock Renormalization group analysis of nonlinear diffusion equations with
  time dependent diffusion coefficients.
\newblock Artigo em preparação.

\bibitem{braga03a}
G.~A. Braga, F.~Furtado, J.~M. Moreira, L.~T. Rolla.
\newblock Renormalization group analysis of nonlinear diffusion equations with
  periodic coefficients.
\newblock {\em Multiscale Modeling and Simulation}, 1~(4), pp.~630--644, 2003.

\bibitem{braga03b}
G.~A. Braga, J.~M. Moreira, L.~T. Rolla.
\newblock An{\'a}lise assint{\'o}tica de solu{\c c}{\~o}es de {EDP}'s via
  grupos de renormaliza{\c c}{\~a}o.
\newblock Em {\em Atas do 57$^{\rm o}$ Semin{\'a}rio Brasileiro de
  An{\'a}lise}, pp. 213--280, 2003.

\bibitem{bric-kupa}
J.~Bricmont, A.~Kupiainen.
\newblock Renormalizing partial differential equations.
\newblock {\em Constructive Physics}, 446, pp.~83--115, 1995.

\bibitem{bric-kupa-lin}
J.~Bricmont, A.~Kupiainen, G.~Lin.
\newblock Renormalization group and asymptotics of solutions of nonlinear
  parabolic equations.
\newblock {\em Communications in Pure and Applied Mathematics}, 47,
  pp.~893--922, 1994.

\bibitem{chen}
L.~Chen, N.~Goldenfeld.
\newblock Numerical renormalization group calculations for similarity solutions
  and travelling waves.
\newblock {\em Physical Review E}, 51, pp.~5577--5581, 1995.

\bibitem{dagan89}
G.~Dagan.
\newblock {\em Flow and Transport in Porous Formations}.
\newblock Springer-Verlag, New York, 1989.

\bibitem{oliveira01}
C.~R. de~Oliveira.
\newblock {\em Introdução à análise funcional}.
\newblock Publicações Matemáticas. IMPA, Rio de Janeiro, 2001.

\bibitem{zuazua2000}
G.~Duro, E.~Zuazua.
\newblock Large time behavior for convection-diffusion equations in {R$^n$}
  with periodic coefficients.
\newblock {\em Journal of Differential Equations}, 167, pp.~275--315, 2000.

\bibitem{evans98}
L.~C. Evans.
\newblock {\em Partial Differential Equations}, volume~19 de {\em Graduate
  Studies in Mathematics}.
\newblock AMS, 1998.

\bibitem{friedman83}
A.~Friedman.
\newblock {\em Partial Differential Equations of Parabolic Type}.
\newblock Krieger, 1983.

\bibitem{furt-glim-lind-pere}
F.~Furtado, J.~Glimm, B.~Lindquist, F.~Pereira.
\newblock Multi-length scale calculations of mixing length growth in tracer
  floods.
\newblock Em F.~Kovarik, editor, {\em Proc. of the Emerging Technologies
  Conference}, pp. 251--259. Institute for Improved Oil Recovery, University of
  Houston, Houston TX, 1990.

\bibitem{glim-lind-pere-zhan}
J.~Glimm, B.~Lindquist, F.~Pereira, Q.~Zhang.
\newblock A theory of macrodispersion for the scale up problem.
\newblock {\em Transport in Porous Media}, 13, pp.~97--122, 1993.

\bibitem{gold92}
N.~Goldenfeld.
\newblock {\em Lectures on Phase Transitions and the Renormalization Group}.
\newblock Addison-Wesley, Reading, 1992.

\bibitem{halliday}
D.~Halliday, R.~Resnick, K.~S. Krane.
\newblock {\em Physics}.
\newblock John Wiley, 5$^{\rm a}$~edi\c c\~ao, 2000.

\bibitem{iorio}
R.~I{\'o}rio~Jr., V.~M. I{\'o}rio.
\newblock {\em Equa{\c c}{\~o}es Diferenciais Parciais: Uma Introdu{\c
  c}{\~a}o}.
\newblock IMPA, Rio de Janeiro, 1988.

\bibitem{isaia02}
V.~Isaia.
\newblock {\em Intermediate Asymtotic Behavior of Nonlinear Parabolic {PDE}s
  via a Renormalization Group Approach: a Numerical Study}.
\newblock Tese de doutorado, Department of Mathematics, University of Wyoming,
  Laramie, Wyoming, 2002.

\bibitem{jikov94}
V.~V. Jikov, S.~M. Kozlov, O.~A. Oleinik.
\newblock {\em Homogenization of Differential Operators and Integral
  Functionals}.
\newblock Springer Verlag, Berlin, 1994.

\bibitem{leveque92}
R.~J. LeVeque.
\newblock {\em Numerical Methods for Conservation Laws}.
\newblock Lectures in Mathematics ETH Zürich. Birkhäuser Verlag, Boston,
  2$^{\rm a}$~edi\c c\~ao, 1992.

\bibitem{moreira}
J.~M. Moreira.
\newblock {\em O Comportamento Assintótico de Soluções da Equação do Calor
  não-Linear Via Grupos de Renormalização}.
\newblock Tese de mestrado, UFMG, Belo Horizonte, Minas Gerais, 2002.

\bibitem{naddaf97}
A.~Naddaf, T.~Spencer.
\newblock On homogenization and scaling limit of some gradient perturbations of
  a massless free field.
\newblock {\em Comm. Math. Phys.}, 183, pp.~55--84, 1997.

\bibitem{papa02}
G.~Papanicolaou.
\newblock Diffusion in random media.
\newblock {\em Surveys in Applied Mathematics}, pp. 205--255, 1995.

\bibitem{peaceman77}
D.~W. Peaceman.
\newblock {\em Fundamentals of Numerical Reservoir Simulation}.
\newblock Elsevier, New York, 1977.

\bibitem{tartar77}
L.~Tartar.
\newblock Homog{\'e}n{\'e}isation.
\newblock Cours Pecout ou College de France, 1977.

\end{thebibliography}
